\documentclass[12pt]{article}                
\usepackage{graphicx}
\usepackage{psfrag}
\usepackage{amsmath, amssymb}
\usepackage{mathptmx}
\usepackage{ntheorem}
\usepackage[nottoc,numbib]{tocbibind}
\usepackage{hyperref}

\usepackage{tocloft}

\renewcommand{\cftpartpagefont }{\scshape}

\renewcommand{\cftpartpagefont}[2]{}

\usepackage{setspace}
\makeatletter
\newsavebox\myboxA
\newsavebox\myboxB
\newlength\mylenA
\newcommand*\xoverline[2][0.75]{%
    \sbox{\myboxA}{$\m@th#2$}%
    \setbox\myboxB\null
    \ht\myboxB=\ht\myboxA%
    \dp\myboxB=\dp\myboxA%
    \wd\myboxB=#1\wd\myboxA
    \sbox\myboxB{$\m@th\overline{\copy\myboxB}$}
    \setlength\mylenA{\the\wd\myboxA}
    \addtolength\mylenA{-\the\wd\myboxB}%
    \ifdim\wd\myboxB<\wd\myboxA%
       \rlap{\hskip 0.5\mylenA\usebox\myboxB}{\usebox\myboxA}%
    \else
        \hskip -0.5\mylenA\rlap{\usebox\myboxA}{\hskip 0.5\mylenA\usebox\myboxB}%
    \fi}
\makeatother

\newcommand{\spe}{{\cal S}}

\newcommand{\A}{{\cal A}}

\newcommand{\M}{{\cal M}}

\newcommand{\PP}{{\cal P}}

\newcommand{\I}{{\rm I}}

\newcommand{\eps}{{\varepsilon}}
\newcommand{\ch}{{\mbox{\rm ch}}}

\newcommand{\smsp}{\hspace{0.3mm}}
\newcommand{\e}{\mathbb{E}}

\newcommand{\p}{\mathbb{P}}

\newcommand{\Reals}{\mathbb{R}}
\newcommand{\Natural}{\mathbb{N}}
\newcommand{\la}{\langle}
\newcommand{\ra}{\rangle}
\newcommand\qed{\hfill\hbox{\rlap{$\sqcap$}$\sqcup$}}

\newtheorem{lemma}{Lemma}
\newtheorem{theorem}{Theorem}

\theoremstyle{nonumberplain}
\newcommand\specialref{}

\begin{document}

\title{\vspace{2cm}\textsc{Introduction to the SK model}}
\author{\textsc{Dmitry Panchenko}}
\date{}
\maketitle


\vspace{14cm}
\begin{center}
\textsc{Current Developments in Mathematics 2014}
\end{center}

\thispagestyle{empty}

\setstretch{1.4}

\pagebreak

\tableofcontents

\thispagestyle{empty}

\pagebreak


\section{The Sherrington-Kirkpatrick model} 
\setstretch{1.08}

The field of spin glass models originates in statistical physics. Some of the original models, such as the Edwards-Anderson (EA) model \cite{EA} and the Sherrington-Kirkpatrick (SK) model \cite{SK}, were introduced with the goal of understanding the unusual magnetic properties of some metal alloys with competing ferromagnetic and anti-ferromagnetic interactions, such as copper/magnesium or gold/iron,  in which the magnetic spins of the component atoms are not aligned in a regular pattern. The Hamiltonian of the Ising version of the EA model is given by 
\begin{equation}
H_N(\sigma) = \sum_{i\sim j} g_{ij} \sigma_i \sigma_j,
\label{EA}
\end{equation}
where $\sigma = (\sigma_i)_{i\in \Lambda}$ is a configuration vector of spins $\sigma_i \in \{-1,+1\}$ on a $d$-dimensional lattice $\Lambda = \{1,\ldots, N\}^d$, the notation $i\sim j$ means that $i$ and $j$ are neighbours on this lattice, and $(g_{ij})$ are interaction parameters that can be either positive (ferromagnetic) or negative (anti-ferromagnetic). These interaction parameters are modelled as independent random variables and are collectively called \emph{the disorder} of the model. The SK model was introduced as a mean field simplification of the EA model in which all spins interact and, thus, the lattice is replaced by a complete graph on $N$ sites.  The Hamiltonian of the Sherrington-Kirkpatrick model is given by  
\begin{equation}
H_N(\sigma) = \frac{1}{\sqrt{N}} \sum_{1\leq i<j\leq N} g_{ij}\sigma_i \sigma_j,
\end{equation}
where the normalization by $\sqrt{N}$ is introduced to keep the typical energy per spin of order $1$, and can be viewed as a tuning of the interaction strength.

Another common way to introduce the SK model is to consider the following optimization problem, called the Dean's problem. Suppose we have a group of $N$ people indexed by the elements of $\{1,\ldots, N\}$ and a collection of parameters $g_{ij}$ for $1\leq i < j \leq N$, called the interaction parameters, which describe how much people $i$ and $j$ like or dislike each other. A positive parameter means that they like each other and a negative parameter means that they dislike each other. We will consider all possible ways to separate these $N$ people into two groups and it will be convenient to describe such partitions using vectors of $\pm 1$ labels with the agreement that people with the same label belong to the same group. Therefore, vectors 
$$
\sigma = (\sigma_1,\ldots,\sigma_N) \in \Sigma_N = \{-1,+1\}^N
$$
describe $2^N$ possible such partitions. For a given configuration $\sigma$, let us write $i\sim j$ whenever $\sigma_i \sigma_j=1$ or, in other words, if $i$ and $j$ belong to the same group, and consider the following so-called \emph{comfort function},
\begin{equation}
c(\sigma) = \sum_{i<j} g_{ij} \sigma_i \sigma_j 
=
\sum_{i\sim j} g_{ij} - \sum_{i\not \sim j} g_{ij}.
\label{ch11DeanComfort}
\end{equation}
The Dean's problem is then to maximize the comfort function over all configurations $\sigma$ in $\Sigma_N$. This objective is pretty clear -- we would like to keep positive interactions as much as possible within the same groups and separate negative interactions into different groups. It is an interesting problem to understand how this maximum behaves in a `typical situation' and one natural way to formalize this is to model the interaction parameters $(g_{ij})$ as random variables. The simplest choice is to let the interactions be independent among pairs with the standard Gaussian distribution (although, in some sense, the choice of the distribution is not important). 

One of the central questions in the SK model is to compute the properly scaled maximum, 
$$
\frac{1}{N}\max_{\sigma} H_N(\sigma),
$$ 
in the thermodynamic limit $N\to\infty$.  The normalization $\sqrt{N}$ in the definition of the Hamiltonian was introduced, in some sense, for convenience of notation to ensure that the maximum scales linearly with $N$. The above quantity is random, but due to classical Gaussian concentration inequalities, its limit is almost surely equal to the limit of its expected value,
\begin{equation}
\frac{1}{N}\e \max_{\sigma} H_N(\sigma).
\label{ch11limmax}
\end{equation}
Let us mention right away that the physicists predicted precise asymptotics for this maximum,
$$
\lim_{N\to\infty} \frac{1}{N}\e \max_{\sigma} H_N(\sigma) = 0.7633\ldots,
$$
which is a consequence of the Parisi formula that will be discussed below. We mentioned above that the distribution of the interactions is not important, and Carmona and Hu proved in \cite{CarmonaHu} that the limit will be the same as long as $\e g_{ij} = 0$ and $\e g_{ij}^2 = 1$. For example, one can replace Gaussian interactions with random $\pm 1$ signs.

\medskip
\noindent
\textbf{Example.} If we use the limit $0.76\ldots$ above as an approximation for the case of $N=10,000$, one can check that optimal solution of the Dean's problem will result in a given person having on average $2462$ people they dislike within their group. Compare this to $2500$ in the random solution, when the choice of side is decided by a fair coin. \qed

\medskip
\noindent
A standard approach to this random optimization problem is to think of it as the zero-temperature case of a general family of problems at positive temperature and, instead of dealing with the maximum in (\ref{ch11limmax}) directly, first to try to compute the limit of its `smooth approximation'
\begin{equation}
\frac{1}{N\beta}\smsp \e \log \sum_{\sigma\in\Sigma_N} \exp \beta H_N(\sigma)
\label{ch11limFE}
\end{equation}
for every parameter $\beta=1/T>0,$ which is called \emph{the inverse temperature parameter}. To relate the above two quantities, let us write
\begin{align}
\frac{1}{N}\smsp \e \max_{\sigma\in\Sigma_N} H_N(\sigma)
& \,\leq\,
\frac{1}{N\beta}\smsp \e \log \sum_{\sigma\in\Sigma_N} \exp \beta H_N(\sigma)
\nonumber
\\
& \,\leq\,
\frac{\log 2}{\beta}
+\frac{1}{N}\smsp \e \max_{\sigma\in\Sigma_N} H_N(\sigma),
\label{ch11compare}
\end{align}
where the lower bound follows by keeping only the largest term in the sum inside the logarithm and the upper bound follows by replacing each term by the largest one. This shows that (\ref{ch11limmax}) and (\ref{ch11limFE}) differ by at most $\beta^{-1}\log 2$ and (\ref{ch11limFE}) approximates (\ref{ch11limmax}) when the inverse temperature parameter $\beta$ goes to infinity. Let us introduce the notation
\begin{equation}
F_N(\beta) = \frac{1}{N}\smsp \e \log Z_N(\beta),
\,\,\mbox{ where }\,\,
Z_N(\beta) = \sum_{\sigma\in\Sigma_N} \exp \beta H_N(\sigma).
\label{ch11FE}
\end{equation}
The quantity $Z_N(\beta)$ is called the \emph{partition function} and $F_N(\beta)$ is called the \emph{free energy} of the model. There is a short proof due to Guerra and Toninelli \cite{GuerraToninelli} that the limit of the free energy
$$
F(\beta) = \lim_{N\to\infty} F_N(\beta)
$$ 
exists, although it took a while to find this proof. The limit of (\ref{ch11limFE}) then equals to $\beta^{-1}F(\beta)$. It is easy to see, by H\"older's inequality, that the quantity
$$
\beta^{-1}(F_N(\beta) - \log 2) =
\frac{1}{N\beta}\smsp \e \log \frac{1}{2^{N}}\sum_{\sigma\in\Sigma_N} \exp \beta H_N(\sigma)
$$
is increasing in $\beta$ and, therefore, so is $\beta^{-1}(F(\beta) - \log 2)$, which  implies that the limit $\lim_{\beta\to\infty}\beta^{-1} F(\beta)$ exists. It then follows from (\ref{ch11compare}) that
\begin{equation}
\lim_{N\to\infty} \frac{1}{N}\smsp \e \max_{\sigma\in\Sigma_N} H_N(\sigma)
=
\lim_{\beta\to\infty}\frac{F(\beta)}{\beta}.
\label{ch11limmaxF}
\end{equation}
The problem of computing (\ref{ch11limmax}) was reduced to the problem of computing the limit $F(\beta)$ of the free energy $F_N(\beta)$ at every positive temperature. 

A formula for $F(\beta)$ was proposed by Sherrington and Kirkpatrick in their original paper \cite{SK} based on the so-called replica formalism, see Section \ref{Sec4label} below. At the same time, they observed that their so-called replica symmetric solution exhibits `unphysical behavior' at low temperature, which means that it can only be correct at high temperature. The correct formula for $F(\beta)$ at all temperatures was famously discovered by Parisi several years later in \cite{Parisi79}, \cite{Parisi}. Let us first describe this formula.

\section{The Parisi formula}

From now on, it will be convenient to redefine the Hamiltonian of the SK model by 
\begin{equation}
H_N(\sigma) = \frac{1}{\sqrt{N}} \sum_{i,j =1}^N g_{ij}\sigma_i \sigma_j.
\label{SKH}
\end{equation}
Compared to the definition above where the sum was over $i<j$, this definition essentially differs by a factor of $\sqrt{2}$, since the contribution of the diagonal elements is negligible for large $N$ and the sum of two independent Gaussian random variables $g_{ij}+g_{ji}$ is equal in distribution to $\sqrt{2} g_{ij}$.

\medskip
\noindent
\textbf{Remark.} It is customary to study the Sherrington-Kirkpatrick model in a slightly more general setting, by adding the external field term $h \sum_{i=1}^N \sigma_i$ to the Hamiltonian $H_N(\sigma)$ for some external field parameter $h.$ All the results discussed below can be proved with very minor modifications in this case so, for simplicity of notation, we will omit the external field. 

\begin{figure}[t]
\centering
\psfrag{que_0}{$q_0$}
\psfrag{que_1}{$q_1$}
\psfrag{que_2}{$q_{2}$}
\psfrag{que_3}{$q_{r-1}$}
\psfrag{que_4}{$q_r=1$}
\psfrag{que_5}{$\cdots$}
\psfrag{zeta_0}{$\zeta_0$}
\psfrag{zeta_1}{$\zeta_1$}
\psfrag{zeta_2}{$\zeta_{r-2}$}
\psfrag{zeta_3}{$\zeta_{r-1}$}
\psfrag{zeta_4}{$\zeta_r= 1$}
\includegraphics[width=0.65\textwidth]{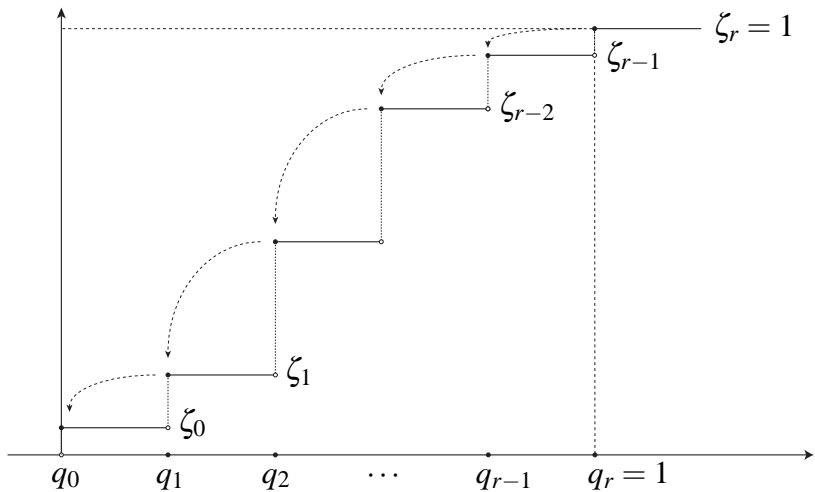}
\caption{\label{Fig0} Solving the Parisi equation (\ref{ParisiEq}) for step function $\zeta$.}
\end{figure}

\medskip
\noindent
The limit of the free energy is given by the following variational formula discovered by Parisi,
\begin{equation}
 \lim_{N\to\infty} F_N(\beta) = F(\beta) = \inf_{\zeta}\, \Bigl(\log 2 + \PP(\zeta) -{\beta^2}\int_{0}^{1}\! \zeta(t)t\, dt\Bigr),
\label{ch30Parisi}
\end{equation}
where the infimum is taken over all cumulative distribution functions $\zeta$ on $[0,1],$  that is, non-decreasing right-continuous functions such that $\zeta(-0)=0$ and $\zeta(1)=0$. The functional $\PP(\zeta)$ is defined as $\PP(\zeta) = f(0,0)$ where $f(t,x)$ is the solution of the parabolic differential equation
\begin{equation}
\frac{\partial f}{\partial t} = -\beta^2\Bigl(
\frac{\partial^2 f}{\partial x^2} + \zeta(t)
\Bigl(\frac{\partial f}{\partial x}\Bigr)^2
\Bigr)
\label{ParisiEq}
\end{equation}
on $[0,1]\times \Reals$ with the boundary condition $f(1,x) = \log \ch(x)$, solved backwards from $t=1$ to $t=0.$ There are several point to be made here.
\begin{enumerate}
\item The parameter $\zeta$ is called \emph{the functional order parameter}, and it has an important physical meaning that will be discussed below.

\item The entire functional being minimized over $\zeta$ is known to be Lipschitz with respect to the $L^1$-distance $\int_{[0,1]} |\zeta_1(x) - \zeta_2(x)|\,dx$ and, therefore, we can minimize over step functions $\zeta$ with finitely many steps.

\item For each step function $\zeta$, the solution of the above equation can be written explicitly. Namely, on any interval where $\zeta$ is constant, the function $\Phi= \exp \zeta f$ satisfies the heat equation 
$$
\frac{\partial \Phi}{\partial t} = - \beta^2 \frac{\partial^2 \Phi}{\partial x^2},
$$
so the function $f(t,x)$ at the beginning of the interval can be written explicitly in terms of the function at the end of the interval. Making the corresponding change of variables on each interval where $\zeta$ is constant, we can write $f(0,x)$ as some recursive formula starting from $f(1,x)$, see Figure \ref{Fig0}.
\end{enumerate}
In fact, in Parisi's original work and in all rigorous proofs, this recursive explicit formula arises naturally from certain computations, so the above equation is simply a convenient compact way to package this recursion.

Let us mention right away that the first proof of the Parisi formula was given by Talagrand in \cite{TPF} building upon an important idea of Guerra in \cite{Guerra}, and another proof was given in Panchenko \cite{PPF} based on the results that we will be reviewed below. We will discuss these proofs briefly in Section \ref{SecLabel-ParisiFormula}. By now, the analogues of the Parisi formula were proved for various modifications of the SK model: spherical SK model in Talagrand \cite{T-sphere} (the formula was discovered by Crisanti and Sommers in \cite{CS}, and another proof was given in Chen \cite{Chen-sphere}); Ghatak-Sherrington model in Panchenko \cite{P-SKgen}; and multi-species SK model in Panchenko \cite{multi-species} (see Section \ref{SecLabel-Problems} for discussion of this and related models).

\section{Beyond the Sherrington-Kirkpatrick model}\label{SecLabel-Beyond}

The original discovery of Parisi generated a lot of activity in the physics community and grew into a far-reaching theory summarized in the classic book \cite{MPV}. We will describe the ideas within the Sherrington-Kirkpatrick model in detail below, but first, to illustrate the success of this theory, let us mention several other models where results analogous to the Parisi formula were obtained.

\medskip
\noindent
\textbf{Minimum matching.} 
Consider the complete graph $K_N$ on $N$ vertices and generate `edge-lengths' $g_{ij}$ independently from some distribution with density $\rho$. The model has \emph{pseudo-dimension} $d$ if 
$$
\lim_{x \downarrow 0} \frac{\rho(x)}{x^{d-1}} =1.
$$
Here, the choice of the constant $1$ on the right hand side is not important and can be adjusted by rescaling the edge-lengths $g_{ij}$. If one picks a point $b$ at random in a neighbourhood of some point $a\in \Reals^d$ then probability that $b\in B_x(a)$ will be proportional to $x^d$, so the density of the distance $\|a-b\|$ will be proportional to $x^{d-1}$ near zero, which explains the name `pseudo-dimension'.

If $N$ is even then a \emph{perfect matching} in $K_N$ is a set of edges with each vertex incident to exactly one of these edges. Let $M_N$ be the minimum total length of edges among all perfect matchings. The following formula was obtained in the work of M\'ezard and Parisi in \cite{MP1, MP3, MP4}:
$$
\lim_{N\to\infty}\frac{\e M_N}{N^{1-1/d}} = d \int_{-\infty}^\infty \, G(x)e^{-G(x)} \,dx,
$$
where $G$ is a solution of the integral equation
$$
G(x) = 2 \int_{-x}^\infty (x+y)^{d-1} e^{-G(y)}\, dy.
$$
When $d=1$, this integral equation has unique solution $G(x) = \log(1+e^{2x})$ and the above limit is equal to $\pi^2/12$, which was first proved by Aldous in \cite{Aldous}. The general case for $d\geq 1$ was solved by W\"astlund in \cite{Wastlund2}. 

\medskip
\noindent
\textbf{Traveling salesman problem.} In the same setting as above, let $L_N$ denote the minimum sum of the edge-lengths of a cycle that visits each vertex exactly once. In this case, M\'ezard-Parisi \cite{MP2, MP3} and M\'ezard-Krauth \cite{KM} obtained that
$$
\lim_{N\to\infty}\frac{\e L_N}{N^{1-1/d}} = \frac{d}{2} \int_{-\infty}^\infty \, G(x)(1+G(x))e^{-G(x)} \,dx,
$$
where $G$ is a solution of the integral equation
$$
G(x) = \int_{-x}^\infty (x+y)^{d-1} (1+G(y))e^{-G(y)}\, dy.
$$
This was proved rigorously for $d= 1$ by W\"astlund in \cite{Wastlund1}. In this case, $\lim_{N\to\infty}\e L_N = 2.04\ldots.$

\medskip
\noindent
\textbf{Diluted spin glass models.} Standard example of diluted models is the diluted version of the Sherrington-Kirkpatrick model corresponding to the Hamiltonian
\begin{equation}
H_N(\sigma) = \sum_{k\leq \pi(\lambda N)} g_k \sigma_{i_k} \sigma_{j_k},
\label{SKdiluted}
\end{equation}
where $\lambda>0$ is called the connectivity parameter, $\pi(\lambda N)$ is a Poisson random variable with the mean $\lambda N$, interaction parameters $(g_k)$ are still independent standard Gaussian, and the indices $i_k$ and $j_k$ are all chosen independently at random from $\{1,\ldots,N\}$. Essentially, this model replaces the complete graph of all possible interactions in the SK model  by the Erd\"os-R\'enyi random graph, where each edge $(i,j)$ is chosen with small probability $p=2\lambda/N.$ The definition using a Poisson number of terms is due to the classical Poisson approximation to the Binomial distribution. 

Another standard example is the random $K$-sat model, with the Hamiltonian 
$$
H_N(\sigma) = - \sum_{k\leq \pi(\lambda N)} \prod_{j=1}^K \frac{1-J_{j,k} \sigma_{i_{j,k}}}{2},
$$
where $\pi(\lambda N)$ is as above, $(J_{j,k})_{j\geq 1}$ are i.i.d. Bernoulli random variables with $\p(J_{j,k}=\pm 1)=1/2$, and indices $i_{j,k}$ are all chosen independently at random from $\{1,\ldots,N\}$. Notice that each term in the sum over $k\leq \pi(\lambda N)$ is equal either to $0$ or $1$, and it is $0$ if and only if one of the variables $\sigma_{i_{j,k}}$ for $j=1,\ldots, K$ is equal to a randomly chosen sign $J_{j,k}.$ As a result, we can think of these terms as representing random Boolean clauses (disjunctions) on $K$ randomly chosen variables, and these disjunctions are generated independently of each other. The Hamiltonian $H_N(\sigma)$ is minus the number of clauses violated by the configuration $\sigma$ and $\max_\sigma H_N(\sigma) = 0$ only if the formula (conjunction of all these clauses) is satisfiable.

In both models, there is an explicit formula for the limit of the free energy  that originates in physics literature in the work of M\'ezard and Parisi in \cite{Mezard}. The formula is too complicated to mention here, and proving it still remains an open problem, although there is a strong indication that all the assumptions made in the derivation of this formula are correct, \cite{AP, Pspins, HEPS, 1RSB, finiteRSB}. In fact, the formula is so complicated that even numerical computations are based on heuristic methods. By considering this formula at zero temperature (letting $\beta\to\infty$) and using numerical computations shows that, for some $\lambda_K>0$, this limit is equal to $0$ for $\lambda\leq \lambda_K$ and it is strictly negative for $\lambda>\lambda_K$. This means that when the number of clauses is not too big, $\lambda\leq \lambda_K$, 
$$
\max_\sigma H_N(\sigma) = - o(N),
$$
so almost all clauses are satisfied. In practice it means that with with high probability all clauses can be satisfies, but it is an important open problem to show that almost satisfiability in the above sense implies satisfiability. For $\lambda >\lambda_K$, a non-trivial proportion of clauses is not satisfied. For example, for $K=3$, the sat-unsat phase transition is predicted to be $\lambda_3 \approx 4.267$ on the basis of the M\'ezard-Parisi formula.  

The value $\lambda_K$ of the satisfiability threshold was described precisely by Mertens, M\'ezard and Zecchina in \cite{MMZ}, who also determined the large $K$ behaviour
$
\lambda_K=2^K \ln 2-\frac{1}{2}(1+\ln 2)+ o_K(1).
$
Another strong indication of the correctness of the physicists' picture is that this asymptotic behaviour was recently proved rigorously by Coja-Oghlan and Panagiotou in \cite{COKSAT} and, soon after, the exact threshold for large enough $K$ (which is a stronger non-asymptotic result) was proved by Ding, Sly and Sun in \cite{DSS}.

Even though there any many questions in these models that remain open even at the heuristic level, one practical success motivated by these analytic developments was the introduction by M\'ezard, Parisi and Zecchina in \cite{MPZ, MZ} of the Survey Propagation algorithm for solving constraint satisfaction problems, such as $K$-sat, which outperforms (at least for small $K$) previously known algorithms on problems with much larger number of variables $N$ and for values of the parameter $\lambda$ (density of clauses per variable) much closer to the satisfiability threshold.

\medskip
\noindent
\textbf{Other models.} There are rigorous results about a number of other models, such as the Hopfield model and perceptron models, and we refer to the manuscript of Talagrand \cite{SG2, SG2-2} for details.

\section{Before and after the Parisi solution}\label{Sec4label}

Before Parisi's discovery, the physicists already had a method to compute the free energy, called \emph{the replica method}. Let us give a brief sketch without going into details of the computations. One starts with the following simple equation,
$$
\e \log Z_N(\beta) = \lim_{n\downarrow 0}\, \frac{1}{n}\log \e Z_N(\beta)^n.
$$
For a fixed \emph{integer} $n\geq 1$, one can use classical large deviations techniques to compute the limit
\begin{equation}
\lim_{N\to\infty} \frac{1}{Nn}\log \e Z_N(\beta)^n = \max_A \Psi_n(A),
\label{RepMet}
\end{equation}
where $\Psi_n$ is some explicit function of a vector $A\in\Reals^{n(n-1)/2}$. For integer $n\geq 1$, the maximum of $\Psi_n(A)$ is achieved on a vector $A$ with all coordinates equal to the same number $q\in [0,1]$; this fact is called \emph{replica symmetry} and was proved in van Hemmen, Palmer \cite{VHP}. In fact, the same replica symmetric formula holds in (\ref{RepMet}) for all $n\geq 1$, as was shown by Talagrand in \cite{TalCortona} (see also Panchenko \cite{PRM}). Then one hopes that the same expression gives the limit in (\ref{RepMet}) for $n\in(0,1)$, crosses ones fingers and interchanges the limits,
$$
\lim_{N\to\infty} \frac{1}{N}\,\e \log Z_N(\beta) = \lim_{n\downarrow 0}\,\max_A \Psi_n(A).
$$
If we pretend that the replica symmetric expression for the maximum holds for $n\in(0,1)$, we get
$$
\lim_{n\downarrow 0}\,\max_A \Psi_n(A)
=
\log 2 + \e\log \ch(\beta g \sqrt{2q}) + \frac{\beta^2}{2}(1-q)^2,
$$
where $g$ is a standard Gaussian random variable. The parameter $q$ here is the largest solution of the equation $q= \e\, \mbox{th}^2(\beta g \sqrt{2q}).$ For $\beta\leq 1/\sqrt{2}$, there will be only one solution $q=0$, but for $\beta> 1/\sqrt{2}$ there will be one more strictly positive solution. This is the formula that Sherrington and Kirkpatrick discovered in \cite{SK}, and it coincides with the infimum of the Parisi formula taken over distributions concentrated on one point $q\in [0,1]$. This formula was known to be incorrect for large values of $\beta$, so there were several attempts to fix it. 

The weakest point of the above heuristic argument was the assumption that (\ref{RepMet}) holds for $n<1$, as well as that the maximum should be replica symmetric. As a result, these attempts focused on \emph{breaking replica symmetry}, which meant looking for various non-constant choices of the vector $A$. However, they were not successful until Parisi introduced a new fundamental principle, namely, that \emph{a maximum of a function of negative number of variables is a minimum}:
$$
\hspace{-2.2cm}\mbox{\textbf{(The Parisi Axiom)}:\hspace{0.8cm} if } \,\mbox{dim}(A)<0 \,\mbox{ then }\,
\max_A \Psi_n(A) \, {\rightsquigarrow} \, \min_A \Psi_n(A).
$$
Since $A$ consists of $n(n-1)/2<0$ variables for $n<1$, we should minimize over the parameters in $A$ for any particular form of replica symmetry breaking before we let $n\downarrow 0$. In part, this allowed Parisi to discover the correct way to break replica symmetry, which was very non-trivial. 

As we mentioned above, the original discovery of Parisi generated a lot of activity in the physics community and one of the main outcomes was that various features of his original computation were reinterpreted via physical properties of \emph{the Gibbs distribution} of the system,
\begin{equation}
G_N(\sigma) = \frac{\exp \beta H_N(\sigma)}{Z_N(\beta)}.
\end{equation}
In statistical physics, $G_N(\sigma)$ represent the chance to observe the system at temperature $T=1/\beta$ in a particular configuration $\sigma$. Intuitively, since the Gibbs weight $G_N(\sigma)$ of the configuration $\sigma$ corresponds to its relative contribution in the partition function $Z_N(\beta)$ and free energy $F_N(\beta)$, discovering a special structure inside the Gibbs distribution should be helpful in computing the free energy. What the physicists did was, in some sense, reverse engineering -- discovering a special structure of the Gibbs distribution from the Parisi formula and the method by which it was obtained. This was done mostly in Parisi \cite{Parisi83}, M\'ezard, Parisi, Sourlas, Toulouse, Virasoro \cite{M1, M2}, and M\'ezard, Parisi, Virasoro \cite{MPV0,MPV1}. In \cite{MPV1}, a new approach to proving the Parisi formula was developed based on these physical assumptions on the organization of the Gibbs distribution. Recent mathematical results were much closer in spirit to this so-called \emph{cavity approach} from \cite{MPV1}, although one had to find a way to prove the physical properties of the Gibbs distribution, whose origin was not understood well even within physics community. The original replica method of Parisi still remains a big mystery.

One way to describe the structure of the Gibbs distribution $G_N$ is in terms of the random geometry of a large sample $\sigma^1,\ldots,\sigma^n$ from $G_N.$ Let us rescale these configurations, 
$$
\tilde{\sigma}^\ell = \frac{\sigma^\ell}{\sqrt{N}},
$$ 
to make their lengths of order one (in fact, equal to one), and let us look at the distances between all these points. Equivalently, we can consider their scalar products
\begin{equation}
R_{\ell,\ell'} =\tilde{\sigma}^\ell \cdot \tilde{\sigma}^{\ell'} 
=\frac{1}{N} \sum_{i=1}^N \sigma_i^\ell  \sigma_i^{\ell'}, 
\label{overlap}
\end{equation}
called the \emph{overlaps} of configuration $\sigma^\ell$ and $\sigma^{\ell'}$. The Gram matrix 
$$
R_N^n = \bigl(R_{\ell,\ell'}\bigr)_{1\leq \ell,\ell'\leq n}
$$ 
encodes the geometry of the sample up to orthogonal transformations. This matrix is random and there are two sources of randomness -- first, the Gibbs distribution $G_N$ itself is random since it depends on the random disorder $(g_{ij})$ and, second, the sample is generated randomly from the Gibbs distribution. The key predictions of the physicists can be summarized as follows. This will be only the first quick reference, and somewhat imprecise. We will discuss these properties in much more detail in the coming sections.  

\begin{enumerate}
\item[] \textbf{1. Functional order parameter.} The distribution of one overlap,
\begin{equation}
\zeta_N(A) = \e G_N^{\otimes 2}\bigl((\sigma^1,\sigma^2) : R_{1,2} = \tilde{\sigma}^1\cdot \tilde{\sigma}^2 \in A\bigr),
\label{fopsummary}
\end{equation}
is called \emph{the functional order parameter} of the model, because, in some sense, everything else is determined by it. This means, for example, that the distribution of the Gram matrix $(R_{\ell,\ell'})_{\ell,\ell'\geq 1}$  is determined by $\zeta_N$. Of course, this should be understood in the asymptotic sense and can hold only approximately for large $N$.

\item[]\textbf{2. Ultrametricity.} The Gibbs distribution is concentrated on an approximately \emph{ultrametric} subset of configurations $\{-1,+1\}^N$ in the sense that if we consider a sample $\sigma^1,\ldots,\sigma^n$ from $G_N$ then, with high probability, it will form an approximately ultrametric set in $\Reals^N$. This means that the ultrametric distance inequality between any three points can be violated by at most $\eps_N$ for some $\eps_N \to 0$ and, with high probability,
$$
\|\tilde{\sigma}^2-\tilde{\sigma}^3\| \leq \max\bigl(\|\tilde{\sigma}^1-\tilde{\sigma}^2\|, \|\tilde{\sigma}^1-\tilde{\sigma}^3\|\bigr) +\eps_N.
$$
Since $\|\tilde{\sigma}^1-\tilde{\sigma}^2\|^2 = 2(1-R_{1,2})$, this can be also rewritten in terms of overlaps,
$$
R_{2,3} \geq \min(R_{1,2},R_{1,3}) - \delta_N.
$$
This geometric ultrametricity property was, perhaps, the most surprising feature that came out of the interpretation of the original Parisi solution, and it is one of the key reasons why everything is determined by the distribution of one overlap asymptotically.

\item[]\textbf{3. Replica symmetry breaking.} It could happen that the distribution (\ref{fopsummary}) of one overlap concentrates around one value. This is called \emph{replica symmetry}, since all overlaps $R_{\ell,\ell'}$ are essentially equal with high probability. In particular, the ultrametricity property becomes trivial in this case. Replica symmetry breaking means that the distribution of one overlap is non-trivial asymptotically, and this is what the physicists predicted in the SK model at low enough temperature.
\end{enumerate}

\noindent
\textbf{Important remark.} We should clarify an important point right away, which will only appear explicitly much later, when we describe technical results. Notice that the SK Hamiltonian (\ref{SKH}) is symmetric, $H_{N}(-\sigma) = H_N(\sigma)$, while the overlap $R_{1,2}$ is anti-symmetric under the map $\sigma^1\to -\sigma^1$, so the distribution of $R_{1,2}$ is symmetric. If we consider three overlaps $(R_{1,2},R_{1,3}, R_{2,3})$, for the same reason, their distribution will be equal to the distribution of $(-R_{1,2},-R_{1,3}, R_{2,3}).$ This shows that ultrametricity can not be taken literally, because for $0<a<b,$ the triple of overlaps $(b,a,a)$ satisfies ultrametricity, while $(-b,-a,a)$ does not. Also, in physics papers the overlap is always assumed to be non-negative, which can be true only if it concentrates at zero. The explanation here is that one can \emph{remove this symmetry} in a way that does not affect the free energy in the limit, and this is what is implicitly assumed in the physics papers. There are different ways to do this. One way is to add an external field term $h\sum_{i=1}^N\sigma_i$ to the Hamiltonian with $h>0$, which creates a preferred direction and forces the overlaps to be non-negative (this is proved by Talagrand in Section 14.12 of \cite{SG2, SG2-2}). A small value of $h$ will result in a small change of free energy, and one can let $h\downarrow 0$ after computing the limiting free energy. We will take a different approach in this paper and add a special random perturbation to our Hamiltonian that will be of a smaller order and will not have any affect on the limit of the free energy at all. This perturbation will serve a more important purpose (in fact, it will be the main tool in the proof of ultrametricity), but it will also force the overlap to be non-negative.
\qed

\medskip
In addition to the geometric property of ultrametricity, another important part of the physicists' picture was a probabilistic property that describes precisely how the functional order parameter determines the joint distribution of overlaps $(R_{\ell,\ell'})$, as well as the free energy. This can be encoded by the so called \emph{Ruelle Probability Cascades} (RPC) that will be defined later.

Let us now mention why the overlaps play such a crucial role in the SK model. After all, the spin configurations are elements in $\{-1,+1\}^N$, and the fact that the coordinates take values $\pm 1$ could conceivably play a role. However, the SK Hamiltonian (\ref{SKH}) is a Gaussian process with the covariance 
\begin{equation}
\e H_N(\sigma^1)H_N(\sigma^2)
=
\frac{1}{N} \sum_{i,j =1}^N \sigma_i^1 \sigma_j^1  \sigma_i^2 \sigma_j^2
=
N\Bigl(
\frac{1}{N} \sum_{i=1}^N \sigma_i^1  \sigma_i^2 
\Bigr)^2
=
N R_{1,2}^2
\label{Cov}
\end{equation}
that depends on the spin configurations $\sigma^1, \sigma^2$ only through their overlap. Since the distribution of a Gaussian process is determined by its covariance, this means that any orthogonal transformation of the space $\{-1,+1\}^N$ does not affect the distribution of the Hamiltonian and the free energy. For this reason, the overlaps between configurations encode all the relevant information in the SK model. 

All the results we will present below can be proved with only minor modifications in a larger class of models -- the so called \emph{mixed $p$-spin models} with the Hamiltonian 
\begin{equation}
H_N(\sigma) = \sum_{p\geq 1} \beta_p H_{N,p}(\sigma)
\label{mixedH}
\end{equation}
given by a linear combination of {pure $p$-spin} Hamiltonians\index{Hamiltonian!$p$-spin model}
\begin{equation}
H_{N,p}(\sigma)
=
\frac{1}{N^{(p-1)/2}}
\sum_{i_1,\ldots,i_p = 1}^N g_{i_1\ldots i_p} \sigma_{i_1}\cdots\sigma_{i_p},
\label{mixedp}
\end{equation}
where the random variables $(g_{i_1\ldots i_p})$ are standard Gaussian, independent for all $p\geq 1$ and all $(i_1,\ldots,i_p)$. The covariance of this Hamiltonian is given by
\begin{equation}
\e H_N(\sigma^1) H_N(\sigma^2) = N\xi(R_{1,2}),
\,\mbox{ where }\, \xi(x)=\sum_{p\geq 1}\beta_p^2 x^p,
\label{Covxi}
\end{equation}
and we assume that the coefficients $(\beta_p)$ decrease fast enough to ensure that the process is well defined when the sum in (\ref{mixedH}) includes infinitely many terms. Again, the covariance is a function of the overlap, which makes these models very similar to the SK model.

In other models, the overlaps might not encode all the information and one needs to understand the structure of the Gibbs distribution $G_N$ in much more detail. For example, in the diluted SK model, the Parisi ansatz for the overlaps still plays the central role, but it is only a starting point of a more detailed description. Nevertheless, physicists predicted that the picture that describes the behaviour of the overlaps in the SK model should be universal and shared by other models, except that some models might be replica symmetric and the picture is trivial. We will see that recent progress on the mathematical side occurred according to this most optimistic scenario and, instead of trying to prove the physicists' predictions in some ad hoc way in each model, one can obtain it using the same soft approach. We should also emphasize that even when the model is replica symmetric, it does not mean that computing the free energy or $\max H_N(\sigma)$ is easy. For example, the minimal matching and traveling salesman problems are replica symmetric according to the physicists, but current rigorous solutions are decidedly not easy.

\section{Asymptotic Gibbs distributions} 

In order to describe the physics predictions in more detail, it will be convenient to define some limiting object that will be the analogue of the Gibbs distribution in the thermodynamic limit $N\to\infty$. Since we already mentioned above that in the SK model all relevant information about the system is encoded by the overlaps $R_{\ell,\ell'}$ between replicas $(\sigma^\ell)$ sampled from the Gibbs measure, we will define the limiting object in terms of the distribution of the overlap array $(R_{\ell,\ell'})_{\ell,\ell'\geq 1}$. We will use the following obvious property:
\begin{equation}
\bigl(R_{\pi(\ell), \pi(\ell')}\bigr)_{\ell,\ell'\geq 1}  \stackrel{d}{=}\bigl(R_{\ell,\ell'}\bigr)_{\ell,\ell'\geq 1}
\label{ch12weakexch}
\end{equation}
for any permutation $\pi$ of finitely many indices, where the equality is in distribution. This may be called \emph{replica symmetry in distribution}. 

If we knew that the overlap array $(R_{\ell,\ell'})_{\ell,\ell'\geq 1}$ converges in distribution as $N\to\infty$, we could use its limiting distribution to define the asymptotic Gibbs distribution as follows. Notice that, in the limit the array $(R_{\ell,\ell'})_{\ell,\ell'\geq 1}$ would still be non-negative definitive and exchangeable in the sense of (\ref{ch12weakexch}). Such arrays are called \emph{Gram-de Finetti arrays} and the Dovbysh-Sudakov representation \cite{DS} guarantees that they can be generated as follows.\
\begin{theorem}[Dovbysh-Sudakov representation]\label{DSintro}
If $(R_{\ell,\ell'})_{\ell,\ell'\geq 1}$ is a Gram-de Finetti array, there exists a random measure $\eta$ on $H\times \Reals^+$, where $H$ is a separable Hilbert space, such that
\begin{equation}
\bigl(R_{\ell,\ell'}\bigr)_{\ell, \ell'\geq 1}\stackrel{d}{=}
\Bigl(\sigma^{\ell}\cdot \sigma^{\ell'} + a_\ell\, \I(\ell =\ell')\Bigr)_{\ell, \ell'\geq 1},
\label{ch12RDS}
\end{equation}
where $(\sigma^\ell, a_\ell)_{\ell\geq 1}$ is an i.i.d. sample from $\eta$. 
\end{theorem}
These days, this result is proved using the so called Aldous-Hoover representation \cite{Aldous0,Hoover} and the most elegant proof is due to Austin \cite{Austin} (see Sections $1.4$ and $1.5$ in \cite{SKmodel}).

When $(R_{\ell,\ell'})$ is the overlap array, the diagonal elements are all equal to one, so we are only interested in the off-diagonal elements. If we let $G$ be the marginal of the random measure $\eta$ on $H$, the Dovbysh-Sudakov representation states that
$$
\bigl(R_{\ell,\ell'}\bigr)_{\ell\not = \ell'}\stackrel{d}{=}
\bigl(\sigma^{\ell}\cdot \sigma^{\ell'} \bigr)_{\ell \not= \ell'},
$$
where $(\sigma^\ell)_{\ell\geq 1}$ is an i.i.d. sample from $G$. This means that, just like for a system of finite size $N$, in the limit, the overlaps $R_{\ell,\ell'}$ can be generated as scalar products of an i.i.d. sample from a random measure $G$, except that now $G$ is defined on a Hilbert space instead of $\Reals^N$. The measure $G$ is a limiting analogue of the Gibbs distribution $G_N$, and we will call it an \emph{asymptotic Gibbs distribution.} Such definition of an asymptotic Gibbs distribution via the Dovbysh-Sudakov representation was first given by Arguin and Aizenman in \cite{AA} (see also \cite{Arguin}).

One problem with this definition is that we do not know how to show that the distribution of the overlap array $(R_{\ell,\ell'})_{\ell,\ell'\geq 1}$ converges. However, this is not a serious problem because we can consider all subsequential limits and try to show that each limit will have properties predicted by the physicists. It turns out that this is sufficient to prove the Parisi formula. Moreover, if each limit is ultrametric then the entire sequence $G_N$ will be approximately ultrametric, so this key geometric property of ultrametricity does not depend on the existence of the limit. With this in mind, we will now describe the predictions of the Parisi solution in more detail in the language of asymptotic Gibbs distributions.

\section{Structure of asymptotic Gibbs distributions} 

First of all, in the limit, the functional order parameter (\ref{fopsummary}) can be written as
\begin{equation}
\zeta(A) = \e G^{\otimes 2}\bigl((\sigma^1,\sigma^2) : R_{1,2} = \sigma^1\cdot \sigma^2 \in A\bigr)
\label{Gfop}
\end{equation}
in terms of an asymptotic Gibbs measure. In some sense, the infimum in the Parisi formula (\ref{ch30Parisi}) is taken over all possible candidates for this distribution in the thermodynamic limit. 

Since we defined asymptotic Gibbs distributions $G$ on a Hilbert space $H$ via distribution of the overlaps, it is natural to write ultrametricity property first in terms of the overlaps,
\begin{equation}
R_{2,3} \geq \min(R_{1,2}, R_{1,3})
\label{ultra2}
\end{equation}
for any three points in the support of $G$. We could also write this in terms of distances,
\begin{equation}
\|{\sigma}^2 - {\sigma}^3\|
\leq
\max\bigl(\|{\sigma}^1 - {\sigma}^2\|,
\|{\sigma}^1 -{\sigma}^3\|\bigr),
\label{ultra1}
\end{equation}
if we knew that the measure $G$ is concentrated on a sphere in $H$. The method that we will use to analyze asymptotic Gibbs distributions will show that this is indeed the case, so we can keep thinking of ultrametricity either in terms of distances or scalar products. It turns out that the radius of this sphere will be non-random, even though the measure $G$ itself is random. In fact, if $q_*$ is the largest point in the support of $\zeta$ in (\ref{Gfop}), then the radius of the sphere will be $\sqrt{q_*}$.

Ultrametricity can be visualized as clustering of the support of $G$, because the ultrametric inequality (\ref{ultra1}) implies that the relation $\sim_d$ defined by
\begin{equation}
\sigma^1 \sim_d \sigma^2 \Longleftrightarrow \|{\sigma}^1 - {\sigma}^2\| \leq d
\label{diameter}
\end{equation}
is an equivalence relation on the support of $G$ for any $d\geq 0$. As we increase $d$, smaller clusters will merge into bigger clusters and the whole process can be visualized by a branching tree. In general, depending on the distribution $\zeta$ of one overlap, merging of clusters can occur at any distance $d$, but it is convenient to discretize the whole process and record merging of clusters at finitely many values $d$. This discretization will allow us to give a very explicit description of the randomness of the asymptotic Gibbs distribution $G$.

It will also be convenient from now on to express everything in terms of scalar products, or overlaps, instead of distances and we can redefine clusters in terms of overlaps by
\begin{equation}
\sigma^1 \sim_q \sigma^2 \Longleftrightarrow R_{1,2} \geq q.
\label{cluster}
\end{equation}
For a given $q$, we will call these equivalence clusters $q$-clusters. To discretize the overlap, let us choose an integer $r\geq 0$ and $r+1$ disjoint intervals of the type
\begin{equation}
I_p = [q_p,q_p') \ \mbox{ or }\ I_p = \{q_p\}
\label{ch31qs}
\end{equation}
for $0\leq p \leq r$ in such a way that
\begin{equation}
{\rm supp}(\zeta)\subseteq \bigcup_{0\leq p\leq r} I_p \ \mbox{ and }\ \zeta(I_p) >0 \ \mbox{ for all }\ 0\leq p\leq r.
\label{intervals}
\end{equation}
We will allow the possibility of $I_p = \{q_p\}$ only when the point $q_p$ is an atom of $\zeta$ isolated from the right, namely, $\zeta(\{q_p\})>0$ and $\zeta((q_p,q_p +\eps))=0$ for some $\eps>0.$ 

Now that we covered the support of $\zeta$ by (small) intervals $I_p$, we will discretize the overlap $R_{1,2}$ between any two points $\sigma^1,\sigma^2$ in the support of the Gibbs distribution $G$ by defining 
\begin{equation}
\hat{R}_{1,2} = q_p \,\,\mbox{ if }\,\, R_{1,2} \in I_p.
\end{equation}
If two clusters (\ref{cluster}) merge somewhere on the interval $I_p$ then we record the merger at $q_p$ or, in other words, we look only at the clusters
\begin{equation}
\sigma^1 \sim_{q_p} \sigma^2 \Longleftrightarrow R_{1,2} \geq q_p 
\Longleftrightarrow \hat{R}_{1,2} \geq q_p.
\end{equation}
The reason we chose the intervals $I_p$ in such a way that $\zeta(I_p) >0$ is to make sure that the merger can occur on $I_p$. If all the intervals $I_p$ are small then $\hat{R}_{1,2}$ will be a good approximation of $R_{1,2}$ for any two points in the support of $G$.

\begin{figure}[t]
\hspace{-0.5cm}
\centering
\psfrag{que_0}{$q_0$}
\psfrag{que_1}{$q_1=q_{\alpha\wedge\beta}$}
\psfrag{que_2}{$q_{r-1}$}
\psfrag{que_3}{$q_r$}
\psfrag{Alpha}{{$\alpha$}}
\psfrag{Beta}{{$\beta$}}
\psfrag{ALB}{{$\alpha$}\hspace{0.6mm}$\wedge$\hspace{0.6mm}{$\beta$}}
\psfrag{a_0}{$*$}
\psfrag{a_1}{$\Natural$}
\psfrag{a_2}{$\Natural^{r-1}$}
\psfrag{AcodeR}{$\Natural^r$}
\psfrag{n_code1}{$\Natural^{1}$}
\psfrag{n_r1}{$\Natural^{r-1}$}
\psfrag{n_r}{$\Natural^r$}
\includegraphics[width=0.75\textwidth]{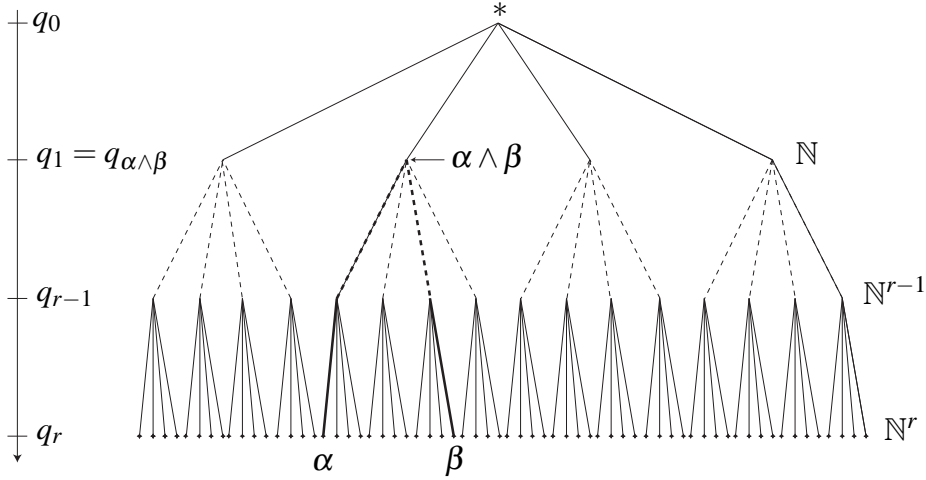}
\caption{\label{Fig1} The leaves $\alpha\in \Natural^r$ index the $q_r$-clusters, called pure states. The figure corresponds to what is called \emph{$r$-step replica symmetry breaking}.}
\end{figure}

We can visualize this discretized clustering process as in Figure \ref{Fig1}, which we will now explain. The leaves of this tree correspond to the smallest clusters, $q_r$-clusters, which will be called the \emph{pure states}. According to the physicists, each $q_p$-cluster contains infinitely many $q_{p+1}$-subclusters, so the clustering tree will be \emph{infinitary}, with each vertex (except the leaves) having infinitely many children. As a result, we can index all the clusters by
\begin{equation}
\A = \Natural^0 \cup \Natural \cup \Natural^2 \cup \ldots \cup \Natural^r,
\label{ch43Atree}
\end{equation}
where $\Natural^0 = \{*\}$, $*$ is the root of the tree that indexes $q_0$-cluster (which is the entire support of the measure $G$) and $q_p$-clusters are indexed by $\Natural^p$. Let us denote by $H_\alpha \subseteq \mbox{supp}(G)$ the equivalence cluster indexed by $\alpha\in\A.$ For any $p\leq r-1$ and $\alpha=(n_1,\ldots,n_p)\in \Natural^{p}$, let us index the $q_{p+1}$-subclusters of $H_\alpha$ by the children of $\alpha$, 
$$
\alpha n : = (n_1,\ldots,n_p,n) \in \Natural^{p+1}
$$
for all $n\in \Natural$. For simplicity of notation, we will write $\alpha n$ instead of $(\alpha,n)$. Each vertex $\alpha$ is connected to the root $*$ by the path
$$
* \to n_1 \to (n_1,n_2) \to\cdots\to (n_1,\ldots,n_p) = \alpha,
$$
and we will denote all the vertices in this path (excluding the root) by
\begin{equation}
p(\alpha) = \bigl\{  n_1, (n_1,n_2),\ldots,(n_1,\ldots,n_p)  \bigr\}.
\label{ch43pathtoleaf}
\end{equation}
Given any two vertices $\alpha, \beta \in\A$, let us denote by
\begin{equation}
\alpha\wedge\beta 
:=
 |p(\alpha) \cap p(\beta)  |
\label{ch43wedge}
\end{equation}
the number of common vertices in the paths from the root $*$ to the leaves $\alpha, \beta$. For two pure states indexed by $\alpha, \beta\in \Natural^r$, if we take $\sigma^1\in H_\alpha$ and $\sigma^2\in H_\beta$ then their overlap $R_{1,2} = \sigma^1\cdot \sigma^2 \in I_{\alpha\wedge\beta }$ and discretized overlap $\hat{R}_{1,2} = q_{\alpha\wedge\beta}$, as shown in Figure \ref{Fig1}.

So far, we have rephrased ultrametricity as clustering and incorporated the predicted infinitary nature of the clustering process into the picture. This corresponds to the the geometric part of the physicists' predictions. The probabilistic part describes the distribution of the Gibbs weights $G(H_\alpha)$ of the pure states indexed by the leaves $\alpha\in\Natural^r$. We need to know this distribution if we want to describe the distribution of discretized overlaps $(\hat{R}_{\ell,\ell'})_{\ell,\ell'\geq 1}$ of the sample $(\sigma^\ell)_{\ell\geq 1}$ from the random measure $G$. This distribution was already understood in M\'ezard, Parisi, Sourlas, Toulouse, Virasoro \cite{M1} and M\'ezard, Parisi, Virasoro \cite{MPV0}. However, a more explicit description emerged in the work of Ruelle \cite{Ruelle} thanks to the connection with Derrida's random energy models. 

In the early eighties, Derrida proposed two simplified models of spin glasses: the random energy model (REM) in \cite{DerridaREM1}, \cite{DerridaREM2}, and the generalized random energy model (GREM) in \cite{DerridaGREM}, \cite{DerridaGREM2}. The Hamiltonian of the REM is given by a vector $(H_N(\sigma))_{\sigma\in\{-1,1\}^N}$ of independent Gaussian random variables with variance $N$, which is a rather classical object. The GREM combines several random energy models in a hierarchical way with the ultrametric structure built into the model from the beginning. Even though these simplified models do not shed light on the Parisi solution of the SK model directly, the behavior of the Gibbs distributions in these models was predicted to be, in some sense, identical to that of the SK model. For example, Derrida and Toulouse showed in \cite{DerridaToulouse} that the Gibbs weights in the REM have the same type of distribution in the thermodynamic limit as the Gibbs weights of the pure states in the SK model, described in \cite{M1}. Later, de Dominicis and Hilhorst \cite{deDH} demonstrated a similar connection between the distribution of the cluster weights in the GREM and the cluster weights in the SK model. Motivated by this connection with the SK model, Ruelle \cite{Ruelle} gave an alternative, much more explicit and illuminating description of the Gibbs distribution of the GREM in the infinite-volume limit  in terms of a certain family of Poisson processes, as follows. 

To describe Ruelle's construction, we need to define one more set of parameters. Let
\begin{equation}
\hat{\zeta}(A) = \e G^{\otimes 2}\bigl((\sigma^1,\sigma^2) : \hat{R}_{1,2} \in A\bigr)
\label{GfopRSB}
\end{equation}
be the distribution of the discretized overlap. If we denote
\begin{equation}
\zeta_p = \hat{\zeta}\bigl(\bigl\{q_0,\ldots,q_p\bigr\}\bigr) = \zeta\bigl(I_0 \cup\ldots\cup I_p\bigr)
\label{ch30zetafop}
\end{equation}
for $p=0,\ldots, r$ then the assumption that $\zeta(I_p)>0$ implies that
\begin{equation}
0 <\zeta_0 <\ldots < \zeta_{r-1} <\zeta_r = 1.
\label{ch31zetas}
\end{equation}
For each vertex $\alpha\in \A$, let us denote by $|\alpha|$ its distance from the root of the tree $*$ or, equivalently, the number of coordinates in $\alpha$, i.e. $\alpha\in \Natural^{|\alpha|}.$ Given the parameters (\ref{ch31zetas}) then, for each $\alpha\in \A\setminus \Natural^r$, let $\Pi_\alpha$ be a Poisson process on $(0,\infty)$ with the mean measure 
\begin{equation}
\mu_{|\alpha|}(dx) = \zeta_{|\alpha|} x^{-1-\zeta_{|\alpha|}} \smsp dx
\label{muzeta}
\end{equation}
and let us generate these processes independently for all such $\alpha$. Let us recall that each Poisson process $\Pi_\alpha$ will be a countable collection of distinct points on $(0,\infty)$ that can be generated in three steps as follows. 
\begin{enumerate}
\item Partition $(0,\infty)=\cup_{m\geq 1}S_m$ into disjoint sets of finite measure (\ref{muzeta}), 
$$
\mu_{|\alpha|}(S_m) = \int_{S_m} \! \zeta_{|\alpha|} x^{-1-\zeta_{|\alpha|}}\smsp dx <\infty,
$$
for example, $S_1 = [1,\infty)$ and $S_m = [1/m,1/(m-1))$ for $m\geq 2$.

\item For each $m\geq 1$, generate a random variable $N_m$ from the Poisson distribution
$$
\p(N_m = k) = \frac{\mu_{|\alpha|}(S_m)^k}{k!} e^{-\mu_{|\alpha|}(S_m)},\,\,\, k=0,1,2,\ldots
$$
with the expected value $\mu_{|\alpha|}(S_m)$. Generate all $(N_m)_{m\geq 1}$ independently of each other.

\item On each set $S_m$, generate $N_m$ points from the probability distribution
$$
\frac{\mu_{|\alpha|}(\,\cdot\,\cap S_m)}{\mu_{|\alpha|}(S_m)}
$$ 
proportional to the measure $\mu_{|\alpha|}$ in (\ref{muzeta}). These points are generated independently over different sets.
\end{enumerate}

\medskip
\noindent
It turns out that statistical properties of the set $\Pi_\alpha$ generated in this way do not depend on the partition $(S_m)$ (see e.g. Kingman \cite{king} or Section 2.1 in Panchenko \cite{SKmodel}). One can now arrange all the points in $\Pi_\alpha$ in the decreasing order,
\begin{equation}
u_{\alpha 1} > u_{ \alpha 2} >\ldots >u_{\alpha n} > \ldots,
\label{ch43us}
\end{equation}
and enumerate them using the children $(\alpha n)_{n\geq 1}$ of the vertex $\alpha$. In other words, parent vertices $\alpha \in \A \setminus \Natural^r$ enumerate independent Poisson processes $\Pi_\alpha$ and child vertices $\alpha n\in \A\setminus \Natural^0$ enumerate individual points $u_{\alpha n}$. Given a vertex $\alpha\in \A\setminus \Natural^0$ and the path $p(\alpha)$ in (\ref{ch43pathtoleaf}), we define
\begin{equation}
w_\alpha = \prod_{\beta \in p(\alpha)} u_{\beta}.
\label{ch43ws}
\end{equation}
Finally, for the leaf vertices $\alpha \in \Natural^r$ we define
\begin{equation}
v_\alpha = \frac{w_\alpha}{\sum_{\beta\in \Natural^r} w_\beta}.
\label{ch43vs}
\end{equation}
One can show that the denominator is finite with probability one, so this is well defined. 

The sequence of these random weights $(v_\alpha)_{\alpha\in \Natural^r}$ is called \emph{the Ruelle Probability Cascades} (RPC) corresponding to the parameters (\ref{ch31zetas}). The probabilistic part of the physicists' picture is saying exactly that the Gibbs weights $G(H_\alpha)$ of the pure states $H_\alpha$, which are $q_r$-clusters indexed by $\alpha\in\Natural^r$, are equal in distribution to the Ruelle Probability Cascades with the parameters defined  in (\ref{ch30zetafop}) in  terms of the functional order parameter $\zeta$. In other words, we can generate them as
\begin{equation}
G(H_\alpha) = v_\alpha.
\end{equation}
In the work of Ruelle, \cite{Ruelle}, it was stated as an evident fact that the Gibbs distributions in the Derrida GREM looks like the weights (\ref{ch43vs}) in the infinite-volume limit, but a detailed proof of this was given later by Bovier and Kurkova in \cite{Bovier}. 

\noindent
\textbf{Summary.}  The core of the physical theory of the Gibbs distribution in the Sherrington-Kirkpatrick model is a combination of a geometric property of ultrametricity and probabilistic description of its randomness via the Ruelle Probability Cascades. The parameters of the RPC depend on the functional order parameter $\zeta$. The power of this description is that it allows us to derive various quantities of interest, most importantly, the limiting free energy given by the Parisi formula. 

\section{Overview of results}

\textbf{First proof of the Parisi formula.} As was mentioned above, the Parisi formula for the free energy was proved by Talagrand in \cite{TPF} building upon a breakthrough idea of Guerra in \cite{Guerra}. Guerra discovered a very natural interpolation that showed that the Parisi formula gives an upper bound on the free energy $F_N(\beta)$ for all $N$,
$$
F_N(\beta) \leq \inf_{\zeta}\Bigl(\log 2 + \PP(\zeta) -{\beta^2}\int_{0}^{1}\! \zeta(t)t\, dt\Bigr).
$$
 Soon after, Talagrand showed how to control the difference in the limit $N\to\infty$ by developing a version of Guerra's interpolation for coupled systems. The results of Guerra and Talagrand were quite stunning considering that even the existence of the limit of the free energy, proved a year earlier by Guerra and Toninelli in \cite{GuerraToninelli}, was considered a big progress. However, this approach managed to prove the Parisi formula without explaining the properties of the Gibbs distribution that we discussed above. 

The predictions of the physicists, unfortunately, came without any heuristic explanation of their possible origin. For example, ultrametricity of the Gibbs distribution was a consequence of the choice of some parameters in the original Parisi computation of the free energy, but where ultrametricity is coming from was a big mystery. One promising approach was suggested by Talagrand in Section $4$ in \cite{TalUltra} based on a generalization of Guerra's interpolation to three coupled copies of the system.  Unfortunately, despite significant effort, it was never completed. On the other hand, the approach was very specific to the SK model and, if it did succeed, it is possible that any further attempts to show how the Parisi ultrametric picture can arise naturally in other models, such as diluted models, would seem hopeless. Another approach which eventually did succeed works equally well in many other models.

\medskip
\noindent
\textbf{The Ghirlanda-Guerra identities.} The first big idea came in the work of Guerra \cite{GuerraGG}, where he proved certain identities on the distribution of the overlaps in the SK model. For example, he showed that for a typical value of the inverse temperature parameter $\beta$, 
\begin{equation}
\e \bigl\la R_{1,2}^2 R_{1,3}^2\bigr\ra \approx \frac{1}{2} \e \bigl\la R_{1,2}^2\bigr\ra  \e \bigl\la R_{1,3}^2\bigr\ra + \frac{1}{2} \e \bigl\la R_{1,2}^4\bigr\ra
\label{GGprelim}
\end{equation}
(the actual statement was on average over $\beta$). Here and below we will always denote by $\la\,\cdot\,\ra$ the average with respect to some Gibbs distribution, either $G_N$ for a finite size system or an asymptotic Gibbs distribution $G.$ For example, in the above equation the term $\la R_{1,2}^4\ra$ can be written as
$$
\bigl\la R_{1,2}^4\bigr\ra = \int \! R_{1,2}^4 \,\, dG_N(\sigma^1)dG_N(\sigma^2)
=  \sum_{\sigma^1,\sigma^2} R_{1,2}^4 \, G_N(\sigma^1)G_N(\sigma^2).
$$
One reason why (\ref{GGprelim}) may hold is that half of the time the overlaps $R_{1,2}$ and $R_{1,3}$ are `generated' independently of each other and half of the time they are set to be equal. More precisely, if the joint distribution of $(R_{1,2},R_{1,3})$  under $\e G_N^{\otimes 2}$ was a mixture of the form
$$
\frac{1}{2}\,\zeta \times \zeta + \frac{1}{2}\, \zeta\circ(x\to(x,x))^{-1}
$$
then we would have
$$
\e \bigl\la f(R_{1,2}) \psi(R_{1,3})\bigr\ra = \frac{1}{2} \e \bigl\la f(R_{1,2})\bigr\ra  \e \bigl\la \psi(R_{1,2})\bigr\ra + \frac{1}{2} \e \bigl\la f(R_{1,2})\psi(R_{1,2})\bigr\ra
$$
for any functions $f$ and $\psi.$ Of course, for finite $N$, this can only be an approximate equality, but in the thermodynamic limit one could hope to have equality. 

The fact that the moments behave as in (\ref{GGprelim}) does not mean that we can expect such stronger statement for joint distributions, but Ghirlanda and Guerra \cite{GG} generalized the original idea of Guerra and showed that (in some generic sense that will be explained below) this stronger statement is correct. In fact, they showed more. Namely, if $\zeta$ is the distribution of one overlap as in (\ref{Gfop}), then conditionally on $(R_{\ell,\ell'})_{\ell,\ell' \leq n}$ the distribution of $R_{1,n+1}$ is given by the mixture 
$$
\frac{1}{n}\, \zeta + \frac{1}{n}\, \sum_{\ell=2}^n \delta_{R_{1,\ell}}.
$$ 
In other words, if we already observed the overlaps $R^n=(R_{\ell,\ell'})_{\ell,\ell' \leq n}$ of $n$ replicas $\sigma^1,\ldots,\sigma^n$ from an asymptotic Gibbs distribution, then the overlap $R_{1,n+1}$ of $\sigma^1$ and a new replica $\sigma^{n+1}$ will take one of the values $R_{1,2},\ldots, R_{1,n}$ with probabilities $1/n$ or, with probability $1/n$, it will be sampled independently according to the distribution $\zeta$. Equivalently,
$$
\e \bigl\la f(R^n) \psi(R_{1,n+1})\bigr\ra = \frac{1}{n} \e \bigl\la f(R^n)\bigr\ra  \e \bigl\la \psi(R_{1,2})\bigr\ra + \frac{1}{n}\sum_{\ell=2}^n \e \bigl\la f(R^n)\psi(R_{1,\ell})\bigr\ra
$$
for any function $\psi$ and any function $f=f(R^n)$. These distributional identities are now called \emph{the Ghirlanda-Guerra identities}. We will explain the main ideas behind these identities below, and explain how one can derive them in a very universal way in many other models as well, using an idea due to Talagrand.

\begin{figure}[t]
\hspace{-0.5cm}
\centering
\includegraphics[width=0.4\textwidth]{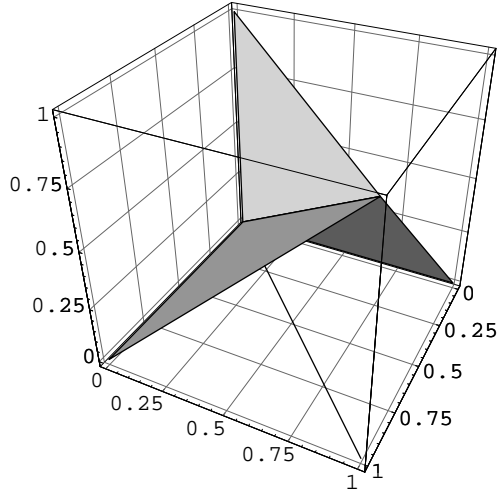}
\caption{\label{Fig10} Ultrametric set of overlaps $(R_{1,2}, R_{1,3}, R_{2,3})$ on $[0,1]^3$. We can assume overlaps to be non-negative by Talagrand's positivity principle.}
\end{figure}

After the discovery of the Ghirlanda-Guerra identities, it was noticed that, if we also assume ultrametricity, the distribution of all overlaps $(R_{\ell,\ell'})_{\ell,\ell' \geq 1}$ can be determined uniquely in terms of the distribution $\zeta$ of one overlap $R_{1,2}$. This is very easy to see geometrically in the case of three overlaps $(R_{1,2}, R_{1,3}, R_{2,3})$, as in Figure \ref{Fig10}. The Ghirlanda-Guerra identities determine all two-dimensional marginal distributions,
$$
\frac{1}{2}\,\zeta \times \zeta + \frac{1}{2}\, \zeta\circ(x\to(x,x))^{-1}.
$$
On the other hand, ultrametricity means than two smallest overlaps are always equal, so the triple lives on the two-dimensional subset of $[0,1]^3$ depicted in Figure \ref{Fig10}. It is clear from the picture  that one can reconstruct the joint distribution of all three overlaps from those marginals. The same idea works for more than three replicas. Every time we add one more replica, the Ghirlanda-Guerra identities describe marginal distributions that involve one new overlap, and ultrametricity allows to reconstruct the joint distribution from those marginals. 

\medskip
\noindent
\textbf{Remark.} Notice that we assumed in Figure \ref{Fig10} that the overlaps are non-negative. We already mentioned in the remark in Section \ref{Sec4label} that we will introduce a small perturbation to the model that will not affect the free energy, but will force the overlap to be non-negative in the limit. In fact, the main purpose of this perturbation will be to prove the Ghirlanda-Guerra identities. At some point, we will prove \emph{Talagrand's positivity principle}, which states that the overlaps are non-negative whenever the Ghirlanda-Guerra identities hold.
\qed

\medskip
\noindent
Since the Ghirlanda-Guerra identities and ultrametricity determine the overlaps $(R_{\ell,\ell'})_{\ell,\ell' \geq 1}$, they essentially determine the randomness of the asymptotic Gibbs distribution $G$, because one can reconstruct $G$ from $(R_{\ell,\ell'})_{\ell,\ell' \geq 1}$ up to isometry. It was believed that, in this sense, the Ghirlanda-Guerra identities and ultrametricity are two \emph{complementary properties} that determine the entire picture in terms of the functional order parameter $\zeta.$ Therefore, if we could prove ultrametricity then, in order to complete the physicists' picture, we would need to answer one more question. If we discretize the support of $G$ as above by discretizing the overlaps and consider the weights $G(H_\alpha)$ of pure states, will these be given by the Ruelle Probability Cascades? The answer is yes, if we can show that the measure defined by the cascades satisfies the Ghirlanda-Guerra identities and, indeed, this was proved by Talagrand \cite{SG} and Bovier and Kurkova \cite{Bovier} (another proof using the Bolthausen-Sznitman \cite{Bolthausen} invariance property for the Ruelle Probability Cascades can be found in \cite{SKmodel}). So understanding where ultrametricity is coming from was now the central problem.

\medskip
\noindent
\textbf{Proving ultrametricity.} It was shown in Panchenko \cite{PUltra} that the Ghirlanda-Guerra identities is not a complementary property to ultrametricity, because they, in fact, imply ultrametricity.  Since the Ghirlanda-Guerra identities can be proved in many models, almost all predictions of the physicists for the Gibbs distribution can be obtained as their consequence (proving replica symmetry breaking utilizes the Parisi formula for the free energy, see Section \ref{SecLabel-PT}). 

The first indication that ultrametricity could possibly be explained by the Ghirlanda-Guerra identities appeared in a seminal work of Arguin and Aizenman \cite{AA}. Instead of the Ghirlanda-Guerra identities they used a closely related property, called the Aizenman-Contucci stochastic stability \cite{AC} (see also \cite{Contucci}), and showed that, under a  technical assumption that the overlap takes only finitely many values in the thermodynamic limit, 
\begin{equation}
R_{1,2} \in \{q_0,\ldots, q_r \},
\label{tecA}
\end{equation}
the stochastic stability implies ultrametricity. Motivated by this development, it was shown in Panchenko \cite{PGG} under the same technical assumption that the Ghirlanda-Guerra identities also imply ultrametricity (an elementary proof was later found in Panchenko \cite{PGG2}). Another approach was given by Talagrand in \cite{Tal-New}. However, according to physicists, at low temperature the overlap it not expected to take finitely many values in the thermodynamic limit, so all these result were not directly applicable to the SK model and could not be used to prove the Parisi formula. Nevertheless, they strongly suggested that this is the right approach and, indeed, it was proved in Panchenko \cite{PUltra} that the Ghirlanda-Guerra identities imply ultrametricity without any technical assumptions. 

We should mention that the proof of the general case in \cite{PUltra}, that will be reproduced below, is very different from all earlier proofs under the assumption (\ref{tecA}). One might hope that the general case can be somehow approximated by the case where (\ref{tecA}) holds, for example, by discretizing the overlap as we did above. However, discretizing the overlap worked well because we assumed ultrametricity. If we consider the array of overlaps $(R_{\ell,\ell'})_{\ell,\ell'\geq 1}$ of a sample from some asymptotic Gibbs distribution $G$ and then discretize the overlaps, the resulting array $(\hat{R}_{\ell,\ell'})_{\ell,\ell'\geq 1}$ might not be positive definite, so it might not correspond to some other Gibbs distribution on a Hilbert space. For this reason, a new approach was needed to solve the general case.

We will start the second part with a few basic techniques used in the proofs. Then, we will explain the idea behind the Ghirlanda-Guerra identities and show how one can ensure their validity in the SK model and other models by introducing a small perturbation of the Hamiltonian. We will prove ultrametricity as a consequence of the Ghirlanda-Guerra identities and show that the distribution of all overlaps can be reconstructed in terms of the distribution of one overlap. Finally, we will sketch the proof of the Parisi formula and describe the phase transition in the SK model.

\pagebreak


\section{Miscellaneous Gaussian techniques} \label{chAppGipSectionLabel}

In this section, we will mention several standard techniques, such as the Gaussian integration by parts, and Gaussian concentration, which will be used in the proofs. We will begin with the Gaussian integration by parts. Let $g$ be a centered Gaussian random variable with variance $v^2$ and let us denote the density function of its distribution by
\begin{equation}
\varphi_v(x) = \frac{1}{\sqrt{2\pi}v}\exp\Bigl(-\frac{x^2}{2v^2}\Bigr).
\end{equation}
Since $x\varphi_v(x) = -v^2 \varphi_v'(x),$ given a continuously differentiable function $F:\Reals \to \Reals$, we can formally integrate by parts,
\begin{align*}
\e g F(g) 
& = 
 \int\! x F(x) \varphi_v(x)  \smsp dx
=
- v^2 F(x) \varphi_v(x) \Bigr|_{-\infty}^{+\infty}
+
v^2 \!\int\! F'(x) \varphi_v(x) \smsp dx
\\
& = 
v^2 \!\int\! F'(x) \varphi_v(x) \smsp dx
=
v^2\e F'(g), 
\end{align*}
if the limits $\lim_{x\to \pm \infty} F(x) \varphi_v(x) = 0$ and the expectations on both sides are finite. Therefore,
\begin{equation}
\e g F(g) = \e g^2 \smsp \e F'(g).
\label{chAppGip}
\end{equation}
This computation can be generalized to Gaussian vectors. Let $g = (g_\ell)_{1\leq \ell\leq n}$ be a vector of jointly Gaussian random variables. Given a continuously differentiable function
$$
F=F((x_\ell)_{1\leq \ell\leq n}): \Reals^n \to \Reals
$$
whose partial derivatives satisfy some mild growth conditions, one can similarly show that
\begin{equation}
\e g_1 F(g)
=
\sum_{\ell\leq n} \e (g_1 g_\ell) \smsp \e \frac{\partial F}{\partial x_\ell} (g).
\label{chAppGIP}
\end{equation}
Typical application of this formula (\ref{chAppGIP}) will be as follows. 

Suppose that we have two jointly Gaussian vectors $(x(\sigma))$ and $(y(\sigma))$ indexed by some finite set of indices $\sigma \in \Sigma$. Let $G$ be a measure on $\Sigma$ and let us define a new (random) measure on $\Sigma$ by the change of density
\begin{equation}
G'(\sigma) = \frac{ \exp y(\sigma) }{Z}G(\sigma)
\,\mbox{ where }\,
Z = \sum_{\sigma\in\Sigma} \exp(y(\sigma)) G(\sigma).
\label{chAppGprimeZ}
\end{equation}
Let us denote by $\la\,\cdot\,\ra$ the average with respect to the product measure $G^{\prime\smsp \otimes \infty}$, which means that for any $n\geq 1$ and any function $f=f(\sigma^1,\ldots, \sigma^n)$,
\begin{equation}
\bigl\la f \bigr\ra = \sum_{\sigma^1,\ldots,\sigma^n \in\Sigma } f(\sigma^1,\ldots, \sigma^n)\smsp
G'(\sigma^1) \cdots G'(\sigma^n).
\label{chAppGprimeave}
\end{equation}
The following is a consequence of the Gaussian integration by parts formula in (\ref{chAppGIP}).
\begin{lemma}\label{chAppLemma1}
If we denote $C(\sigma^1,\sigma^2) = \e x(\sigma^1) y(\sigma^2)$ then
\begin{equation}
\e\bigl\la x(\sigma) \bigr\ra = \e \bigl\la C(\sigma^1,\sigma^1) - C(\sigma^1,\sigma^2)\big\ra.
\label{chAppLem}
\end{equation}
\end{lemma}
\noindent\textbf{Proof.}
Let us consider one term in the sum
$$
\e\bigl\la x(\sigma) \bigr\ra
=
\e \sum_{\sigma^1\in\Sigma}  x(\sigma^1) G'(\sigma^1).
$$
Let us view the function
$$
F = G'(\sigma^1) = \frac{\exp y(\sigma^1)}{Z} G(\sigma^1)
$$ 
as a function of $y(\sigma^1)$ and $(y(\sigma^2))_{\sigma^2\in\Sigma}$, which means that we view a copy of $y(\sigma^1)$ that appears in the denominator $Z$ as a separate variable. Notice that
$$
\frac{\partial F}{\partial y(\sigma^1)} =G'(\sigma^1),\,\,
\frac{\partial F}{\partial y(\sigma^2)}
=
-G'(\sigma^1) G'(\sigma^2),
$$
and, since each factor $G'(\sigma) \in [0,1]$, the function $F$ and all its derivatives are bounded and, therefore, all the conditions in the proof of the Gaussian integration by parts formula  (\ref{chAppGIP}) are satisfied. Then, (\ref{chAppGIP}) implies that
$$
\e x(\sigma^1) G'(\sigma^1)
=
C(\sigma^1,\sigma^1) \e G'(\sigma^1) 
-
\sum_{\sigma^2\in\Sigma} C(\sigma^1,\sigma^2) 
\e G'(\sigma^1) G'(\sigma^2).
$$
If we now sum this equality over $\sigma^1\in\Sigma$, we get (\ref{chAppLem}).
\qed

\medskip
\noindent
Using a similar computation, one can easily generalize Lemma \ref{chAppLemma1} in a couple of ways. First, we can consider finite measures $G$ on a countably infinite set $\Sigma$, under the condition that all the variances 
\begin{equation}
\e x(\sigma)^2, \e y(\sigma)^2 \leq a 
\label{chAppvarbou}
\end{equation}
are uniformly bounded over $\sigma\in\Sigma$. Notice that under this condition $Z<\infty$ and the measure $G'$ in (\ref{chAppGprimeZ}) is well defined, since, by Fubini's theorem,
$$
\e \sum_{\sigma\in\Sigma} \exp(y(\sigma))G(\sigma) 
\leq
e^{a/2}G(\Sigma) <\infty.
$$
The following is then a simple exercise.
\begin{lemma}\label{chAppLemma2}
Suppose that $G$ is a finite measure on a countably infinite set $\Sigma$ and (\ref{chAppvarbou}) holds. If $\Phi = \Phi(\sigma^1,\ldots, \sigma^n)$ is a bounded function of $\sigma^1,\ldots, \sigma^n$  then 
\begin{equation}
\e\bigl\la \Phi \smsp x(\sigma^1) \bigr\ra 
= 
\e \Bigl\la 
\Phi
\Bigl(
\sum_{\ell=1}^n C(\sigma^1,\sigma^{\ell}) - n C(\sigma^1,\sigma^{n+1})
\Bigr)
\Big\ra.
\label{chAppLemExa}
\end{equation}
\end{lemma}
Next, consider a countable set $\Sigma$ and some finite measure $G$ on it, and let
\begin{equation}
X=\log \sum_{\sigma\in\Sigma} \exp(g(\sigma)) G(\sigma),
\label{chAppXFE}
\end{equation}
where $(g(\sigma))_{\sigma\in\Sigma}$ is a Gaussian process such that for some constant $a>0$, 
\begin{equation}
\e g(\sigma)^2 \leq a
\,\mbox{ for all }\,
\sigma\in \Sigma.
\label{chAppvvar}
\end{equation}
The following concentration inequality holds (see e.g. Section 1.2 in \cite{SKmodel}).
\begin{theorem}\label{chAppLextra4}
If (\ref{chAppvvar}) holds then, for all $x\geq 0$,
\begin{equation}
\p\bigl( |X-\e X|\geq x\bigr)\leq 2\exp\Bigl(- \frac{x^2}{4a}\Bigr),
\label{chAppGCieq}
\end{equation}
which implies that $\e(X-\e X)^2 \leq 8 a.$
\end{theorem}

\section{The Ghirlanda-Guerra identities}

Before we go into details, let us first sketch the main idea behind the Ghirlanda-Guerra identities. Let us write the free energy as $F_N(\beta) = \e \hat{F}_N(\beta)$, where
$$
\hat{F}_N(\beta) = \frac{1}{N}\log \sum_{\sigma} \exp \beta H_N(\sigma)
$$
is called \emph{quenched free energy}, which is random since we have not yet averaged in the random interactions $(g_{ij})$. Taking the derivative in $\beta$,
$$
\hat{F}_N'(\beta)= \Bigl\la \frac{H_N(\sigma)}{N} \Bigr\ra_{\hspace{-0.5mm}\beta}
\,\,\mbox{ and }\,\,
{F}_N'(\beta)= \e \Bigl\la \frac{H_N(\sigma)}{N} \Bigr\ra_{\hspace{-0.5mm}\beta},
$$
where $\la\,\cdot\,\ra_\beta$ is the average with respect to the Gibbs distribution $G_N(\sigma)$ corresponding to the inverse temperature $\beta$. We made the dependence of this average on $\beta$ explicit for the moment. There are two very basic properties that $F_N(\beta)$ and $\hat{F}_N(\beta)$ satisfy.
\begin{itemize}
\item[] \textbf{Convexity.} Both $\hat{F}_N(\beta)$ and $F_N(\beta)$ are, obviously, convex in $\beta$.

\item[] \textbf{Concentration.} The random free energy $\hat{F}_N(\beta)$ concentrates around its expectation $F_N(\beta)$,
$$
\e \bigl|\hat{F}_N(\beta) -F_N(\beta) \bigr| \leq \frac{3\beta}{\sqrt{N}},
$$
by Theorem \ref{chAppLextra4}, since $\e (\beta H_N(\sigma))^2 = \beta^2 N.$ 
\end{itemize}
When two convex functions are close to each other, their derivatives are also close, at least on average over intervals. As a result, with just a little bit more work one can show that
\begin{equation}
\int_0^1 \e \Bigl\la \Bigl| \frac{H_N(\sigma)}{N} - \e \Bigl\la \frac{H_N(\sigma)}{N} \Bigr\ra_{\hspace{-0.5mm}\beta} \Bigr| \Bigr\ra_{\hspace{-0.5mm}\beta} \, d\beta \to 0.
\label{GGprelim01}
\end{equation}
Scaling the Hamiltonian $H_N(\sigma)$ by a factor of $N$ makes its typical values of order $O(1)$, and this equation states that, at the right scale, the Hamilitonian is concentrated around its average value, at least for typical values of $\beta$. Another way to rephrase it is to say that the Gibbs measure is concentrated on configurations with nearly constant energy.

The Ghirlanda-Guerra identities then arise by observing this concentration on test functions. If we take any bounded function $f_n=f_n(R^n)$ of overlaps of $n$ replicas $\sigma^1,\ldots,\sigma^n$, we must have that
\begin{equation}
\e \Bigl\la f_n \frac{H_N(\sigma^1)}{N} \Bigr\ra_{\hspace{-0.5mm}\beta}
\approx
\e \bigl\la f_n \bigr\ra_{\hspace{-0.5mm}\beta}
\e \Bigl\la \frac{H_N(\sigma^1)}{N} \Bigr\ra_{\hspace{-0.5mm}\beta}
\label{GGsketcheq1}
\end{equation}
 on average over $\beta$. Using that $\e H_N(\sigma^1)H_N(\sigma^2) = N (R_{1,2})^2$ and the Gaussian integration by parts formula (\ref{chAppLemExa}), this equation can be rewritten as
\begin{equation}
\e \bigl\la f_n R_{1,n+1}^2 \bigr\ra_{\hspace{-0.5mm}\beta}
\approx
\frac{1}{n} \e \bigl\la f_n  \bigr\ra_{\hspace{-0.5mm}\beta}\e \bigl\la R_{1,2}^2 \bigr\ra_{\hspace{-0.5mm}\beta}
+\frac{1}{n} \sum_{\ell=2}^n \e \bigl\la f_n R_{1,\ell}^2 \bigr\ra_{\hspace{-0.5mm}\beta}
\label{GGsketcheq2}
\end{equation}
on average over $\beta\in [0,1]$. This is starting to look exactly like the Ghirlanda-Guerra identities that we described above. However, we would like to strengthen the above argument in two ways.
\begin{itemize}
\item[(i)] We would like to have this for a given $\beta$ instead of on average. This is important, because we want to get the exact Ghirlanda-Guerra identities in the limit for a given model, and not on average over models.

\item[(ii)] In (\ref{GGsketcheq2}), we want to be able to replace $R_{\ell,\ell'}^2$ by any power of the overlap $R_{\ell,\ell'}^p$ for integer $p\geq 1$, since we want to have these identities in distribution and not only for the second moment of the overlaps.
\end{itemize}
It turns out that both of these goals can be achieved by adding a small perturbation term to the Hamiltonian of the model. Namely, for all $p\geq 1$, let us consider
\begin{equation}
g_{p}(\sigma)
=
\frac{1}{N^{p/2}}
\sum_{i_1,\ldots,i_p = 1}^N g_{i_1\ldots i_p}' \sigma_{i_1}\cdots\sigma_{i_p},
\label{ch31mixedppert}
\end{equation}
where the coefficients $(g_{i_1\ldots i_p}')$ are again i.i.d. standard Gaussian random variables independent of all the other random variables, and define
\begin{equation}
g(\sigma) = \sum_{p\geq 1} 2^{-p} x_p\smsp g_{p}(\sigma)
\label{ch31mixedHpert}
\end{equation}
for some parameters $(x_p)_{p\geq 1}$ such that $x_p\in[0,3]$ for all $p\geq 1$. Each term $g_p(\sigma)$ in this sum is very similar to the SK Hamiltonian, only it involves interactions between $p$ spins at a time instead of two. Parameters $x_p$ will play a role of individual inverse temperature parameters for each of these terms. The normalization by $N^{p/2}$ in (\ref{ch31mixedppert}) is chosen so that the covariance
\begin{equation}
\e g_p(\sigma^1) g_p(\sigma^2) = R_{1,2}^p
\end{equation}
is the $p$th power of the overlap. Note that
\begin{equation}
\e g(\sigma^1) g(\sigma^2) = \sum_{p\geq 1} 4^{-p} x_p^2\smsp R_{1,2}^p
\label{ch31Covxipert}
\end{equation}
and $g(\sigma)$ is of a smaller (constant) order than $H_N(\sigma)$ because of the additional factor $1/\sqrt{N}$. If the Hamiltonian of our model was $\sqrt{N}g_p(\sigma)$ then in (\ref{GGsketcheq2}) we would have factors $R_{\ell,\ell'}^p$ instead of $R_{\ell,\ell'}^2$, so including all such terms in the sum (\ref{ch31mixedHpert}) will allow us to extract information about all powers of the overlaps simultaneously. 

Our goal will be to show how the Hamiltonian (\ref{ch31mixedHpert}) used as a perturbation will give rise to the Ghirlanda-Guerra identities in a number of models, so we now consider an arbitrary Hamiltonian $H(\sigma)$ on $\Sigma_N = \{-1,+1\}^N$, either random or non-random, and consider its perturbation
\begin{equation}
H^{\mathrm{pert}}(\sigma) = H(\sigma) + s g(\sigma),
\label{ch31Hpert}
\end{equation} 
for some parameter $s\geq 0.$ Later, we will let $s=s_N$ depend on $N$. First of all, the perturbation should be small enough not to affect the free energy 
$$
\frac{1}{N}\smsp \e\log \sum_{\sigma\in\Sigma_N} \exp  H(\sigma)
$$
in the thermodynamic limit. Notice that we do not write the inverse temperature parameter $\beta$ here, and assume that it is absorbed into the definition of the Hamiltonian $H(\sigma)$. Using (\ref{ch31Covxipert}) and the independence of $g(\sigma)$ and $H(\sigma)$, it is easy to see that
\begin{align}
\frac{1}{N}\smsp \e\log \sum_{\sigma\in\Sigma_N} \exp  H(\sigma)
& \leq 
\frac{1}{N}\smsp \e\log \sum_{\sigma\in\Sigma_N} \exp  \bigl(H(\sigma) + s g(\sigma)\bigr)
\label{ch31FEsamepert}
\\
& \leq 
\frac{1}{N}\smsp \e\log \sum_{\sigma\in\Sigma_N} \exp  H(\sigma)
+ \frac{s^2}{2N} \sum_{p\geq 1} 4^{-p} x_p^2.
\nonumber
\end{align}
Both inequalities follow from Jensen's inequality applied either to the sum or the expectation with respect to $g(\sigma)$ conditionally on $H(\sigma)$. Therefore, if we let $s=s_N$ in (\ref{ch31Hpert}) depend on $N$ and
\begin{equation}
\lim_{N\to\infty}  \frac{s_N^2}{N} = 0
\label{ch31tN}
\end{equation} 
then the limit of the free energy is not affected by the perturbation term $s_N g(\sigma).$ On the other hand, if $s=s_N$ is not too small then we can make the approach of Ghirlanda and Guerra work under some mild assumption on the concentration of the quenched free energy. Consider a  function
\begin{equation}
\varphi = \log Z_N =
\log \sum_{\sigma\in\Sigma_N} \exp \bigl(H(\sigma) + s g(\sigma)\bigr),
\label{ch31theta}
\end{equation}
that will be viewed as a random function $\varphi = \varphi\bigl((x_p)\bigr)$ of the parameters $(x_p)$, and suppose that
\begin{equation}
\sup\Bigl\{ \e |\varphi - \e \varphi | : 0\leq x_p\leq 3, p\geq 1\Bigr\}\leq v_N(s)
\label{ch31vt}
\end{equation}
for some function $v_N(s)$ that describes how well $\varphi((x_p))$ is concentrated around its expected value uniformly over all possible choices of the parameters $(x_p)$ from the interval $[0,3].$ We will make the following assumption about the model.

\medskip
\textbf{Concentration Assumption:} There exists a sequence $s= s_N$ such that
\begin{equation}
\lim_{N\to\infty} s_N=\infty
\,\mbox{ and }\,
\lim_{N\to\infty}  \frac{v_N(s_N)}{s_N^{2}} = 0.
\label{ch31GGassumption}
\end{equation}

\medskip
\noindent
\textbf{Example.} For $H(\sigma)=\beta H_N(\sigma)$ with the SK Hamiltonian $H_N(\sigma)$, we can use Theorem \ref{chAppLextra4} with the counting measure $G$ on $\varSigma_N.$ By (\ref{ch31Covxipert}),
$
\e \bigl(\beta H_N(\sigma) + s g(\sigma) \bigr)^2 \leq \beta^2 N + 3s^2
$ 
if all $0\leq x_p \leq 3$ and Theorem \ref{chAppLextra4} implies that $\e (\varphi - \e \varphi)^2 \leq 8 ( \beta^2 N + 3s^2)$. Hence, we can take $v_N(s) = 5 (\beta^2 N + s^2)^{1/2}$ in (\ref{ch31vt}) and it follows that both (\ref{ch31tN}) and (\ref{ch31GGassumption}) hold if we can take $s_N = N^{\gamma}$ for any $1/4<\gamma < 1/2.$

One can also check that this concentration assumption holds in the Edwards-Anderson model, diluted SK model and random $K$-sat model (see Lemma 1 in \cite{HEPS}), and, probably, in many other models, since there are standard techniques for proving concentration inequalities even when the disorder is not Gaussian.
\qed

\medskip
\noindent
Let us now formulate the main result of this section. Let 
\begin{equation}
G_N(\sigma) = \frac{\exp  H^{\mathrm{pert}}(\sigma)}{Z_N}
\,\mbox{ where }\,
Z_N = \sum_{\sigma\in\Sigma_{N}} \exp  H^{\mathrm{pert}}(\sigma)
\label{ch31GNpert}
\end{equation}
be the Gibbs measure corresponding to the perturbed Hamiltonian (\ref{ch31Hpert}) and, as usual, let $\la\,\cdot\, \ra$ denote the average with respect to $G_N^{\otimes \infty}$. For any $n\geq 2, p\geq 1$ and any bounded function $f$ of the overlaps $R^n=(R_{\ell,\ell'})_{ \ell,\ell'\leq n}$ of $n$ replicas, define
\begin{equation}
\Delta(f,n,p) := 
\Bigl|
\e  \bigl\la f R_{1,n+1}^p \bigr\ra -  \frac{1}{n}\e \bigl\la f \bigr\ra \smsp \e\bigl\la R_{1,2}^p\bigr\ra - \frac{1}{n}\sum_{\ell=2}^{n}\e \bigl\la f R_{1,\ell}^p\bigr\ra
\Bigr|.
\label{ch31GG}
\end{equation}
This is the expression that appeared in the equation (\ref{GGsketcheq2}), only now with the $p$th power of the overlap. The quantity (\ref{ch31GG}) depends on the parameters $(x_p)_{p\geq 1}$ in the perturbation Hamiltonian (\ref{ch31mixedHpert}) and we will show that the Ghirlanda-Guerra identities hold asymptotically on average over these parameters, in the following sense. If we think of $(x_p)_{p\geq 1}$ as a sequence of i.i.d. random variables with the uniform distribution on $[1,2]$ and denote by $\e_x$ the expectation with respect to such sequence then the following holds.\index{Ghirlanda-Guerra identities}
\begin{theorem}[The Ghirlanda-Guerra identities]\label{ch31ThGG} 
If $s=s_N$ in (\ref{ch31Hpert}), and (\ref{ch31GGassumption}) holds, then
\begin{equation}
\lim_{N\to\infty} \e_x \smsp \Delta(f,n,p) = 0
\label{ch31GGxlim}
\end{equation} 
for any $p\geq 1, n\geq 2$ and any bounded measurable function $f=f((R_{\ell,\ell'})_{\ell,\ell'\leq n})$.
\end{theorem}
The fact that these identities hold on average over $x$ is no longer a cause for concern, because these parameters appear in the perturbation term and not in the main Hamiltonian of our model. In particular, we can choose $x^N = (x_p^N)_{p\geq 1}$ varying with $N$ such that
\begin{equation}
\lim_{N\to\infty} \smsp \Delta(f,n,p) = 0
\end{equation} 
with this particular choice of parameters rather than on average. If $G_N(\sigma)$ is the Gibbs distribution corresponding to the perturbed Hamiltonian 
\begin{equation}
H^{\mathrm{pert}}(\sigma) = H(\sigma) + s_N  \sum_{p\geq 1} 2^{-p} x_p^N \smsp g_{p}(\sigma)
\label{ch31HpertNN}
\end{equation} 
with the parameters $x^N = (x_p^N)_{p\geq 1}$ in the perturbation term, then we can define asymptotic Gibbs distributions as before along any converging subsequence of the distribution of the overlap array $(R_{\ell,\ell'})$ generated by a sample from $G_N.$ Such asymptotic Gibbs distributions will satisfy the exact form of the Ghirlanda-Guerra identities
\begin{equation}
\e  \bigl\la f R_{1,n+1}^p \bigr\ra =  \frac{1}{n}\e \bigl\la f \bigr\ra \smsp \e\bigl\la R_{1,2}^p\bigr\ra + \frac{1}{n}\sum_{\ell=2}^{n}\e \bigl\la f R_{1,\ell}^p\bigr\ra.
\label{ch31GGlimit}
\end{equation}
for any $p\geq 1, n\geq 2$ and any bounded measurable function $f=f((R_{\ell,\ell'})_{\ell,\ell'\leq n}).$ This will allow us to characterize all such limits.

In addition to the concentration inequality (\ref{ch31vt}), the proof of Theorem \ref{ch31ThGG} will utilize convexity. The following lemma quantifies the fact that if two convex functions are close to each other then their derivatives are also close. 
\begin{lemma}
If $\varphi(x)$ and $\phi(x)$ are two differentiable convex functions then, for any $y>0,$
\begin{equation}
|\varphi'(x)-\phi'(x)|
\leq 
\phi'(x+y) -\phi'(x-y) + \frac{\delta}{y}, 
\label{ch31conv}
\end{equation}
where 
$
\delta =|\varphi(x+y)-\phi(x+y)|+|\varphi(x-y)-\phi(x-y)| 
+ |\varphi(x) - \phi(x)|.
$
\end{lemma}
\noindent\textbf{Proof.}
The convexity of $\phi$ implies that, for any $y>0$,
$$
\phi(x+y)-\phi(x)\leq y\phi'(x+y) 
\,\mbox{ and }\,
\phi'(x-y) \leq \phi'(x).
$$
Therefore,
$$
 \phi(x+y)-\phi(x)\leq y\bigl( \phi'(x) + \phi'(x+y) -\phi'(x-y) \bigr). 
 $$
The convexity of $\varphi$ and the definition of $\delta$ imply
\begin{align*}
y \varphi'(x)
&\leq 
\varphi(x+y)-\varphi(x)
\leq
\phi(x+y)-\phi(x)+\delta
\\
&\leq
y\bigl( \phi'(x) + \phi'(x+y) -\phi'(x-y)\bigr) +\delta.
\end{align*}
Similarly, one can show that
$$
y \varphi'(x)\geq
y\bigl( \phi'(x) - \phi'(x+y) + \phi'(x-y)\bigr)-\delta
$$
and combining these two inequalities finishes the proof.
\qed

\medskip
The main step in the proof of the Ghirlanda-Guerra identities is the following analogue of (\ref{GGprelim01}) for $p$-spin term in the perturbation Hamiltonian. This is where we utilize convexity and concentration.
\begin{theorem}\label{ch31L1}
For any $p\geq 1$,  if $s>0$ is such that $s^{-2}v_N(s) \leq 4^{-p}$ then
\begin{equation}
\int_1^2 
\e \bigl\la \bigl|g_p(\sigma) -\e\bigl\la g_p(\sigma)\bigr\ra \bigr|\bigr\ra 
\smsp dx_p
\leq 2 + 48 \sqrt{v_N(s)}.
\label{ch31GGbasic}
\end{equation}
\end{theorem}
Of course, the condition (\ref{ch31GGassumption}) will ensure that the assumption $s^{-2}v_N(s) \leq 4^{-p}$ is satisfied for $s=s_N$ and $N$ large enough. 

\noindent\textbf{Proof.} Given the function $\varphi$ in (\ref{ch31theta}), we define $\phi = \e \varphi.$ Let us fix $p\geq 1$ and denote $s_p = s 2^{-p}$. We will think of $\varphi$ and $\phi$ as functions of $x=x_p$ only and work with one term, 
$$
s2^{-p} x_p g_p(\sigma) = xs_p g_p(\sigma),
$$
in the perturbation Hamiltonian (\ref{ch31mixedHpert}). First, let us observe that
\begin{equation}
\varphi'(x)=s_p \bigl\la g_p(\sigma)\bigr\ra
\,\mbox{ and }\,
\phi'(x)=s_p \e\bigl\la g_p(\sigma)\bigr\ra.
\label{ch31derivs}
\end{equation}
Since the covariance $\e g_p(\sigma^1) g_p(\sigma^2) = R_{1,2}^p$, the Gaussian integration by parts in Lemma \ref{chAppLemma1} implies
\begin{equation}
\phi'(x)=s_p \e\bigl\la g_p(\sigma)\bigr\ra
=  x s_p^2 \bigl(1-\e \bigl\la R_{1,2}^p \bigr\ra\bigr)
\in \bigl[0, 2x s_p^2 \bigr].
\label{ch31der}
\end{equation}
Differentiating the derivative $\phi'(x)$ in (\ref{ch31derivs}), it is easy to see that
$$
\phi''(x)= s_p^2 \e \bigl( \bigl\la g_p(\sigma)^2 \bigr\ra - \bigl\la g_p(\sigma)\bigr\ra^2\bigr)
= s_p^2 \e\bigl\la \bigl(g_p(\sigma) - \bigl\la g_p(\sigma)\bigr\ra\bigr)^2\bigr\ra,
$$
and integrating this over the interval $1\leq x\leq 2$ and using (\ref{ch31der}) implies
$$
s_p^2 
\int_{1}^{2}\! \e\bigl\la \bigl(g_p(\sigma) - \bigl\la g_p(\sigma) \bigr\ra\bigr)^2\bigr\ra \smsp dx
=\phi'(2)-\phi'(1)\leq 4 s_p^2.
$$
If we cancel $s_p^2$ on both sides then Jensen's inequality implies that 
\begin{equation}
\int_{1}^{2}\! \e\bigl\la \bigl|g_p(\sigma) - \bigl\la g_p(\sigma)\bigr\ra\bigr| \bigr\ra \smsp dx \leq 2.
\label{ch31step2}
\end{equation}
To prove (\ref{ch31GGbasic}), it remains to approximate $\la g_p(\sigma)\ra$ by $\e\la g_p(\sigma)\ra$ and this is where the convexity plays its role. Since $\varphi(x)$ and $\phi(x)$ are convex and differentiable, we can apply the inequality (\ref{ch31conv}). We will consider $1\leq x\leq 2$ and $0\leq y\leq 1$, in which case $0\leq x-y,x,x+y\leq 3$ and we can use the definition (\ref{ch31vt}) to conclude that $\delta$ in (\ref{ch31conv}) satisfies $\e \delta \leq 3v_N(s)$. Averaging the inequality (\ref{ch31conv}), 
\begin{equation}
\e|\varphi'(x)-\phi'(x)| 
\leq
\phi'(x+y) -\phi'(x-y) 
+ \frac{3v_N(s)}{y}.
\label{ch31etheta3vt}
\end{equation}
By (\ref{ch31der}), $|\phi'(x)| \leq 6s_p^2$ for all $0\leq x\leq 3$ and, by the mean value theorem,
\begin{align*}
\int_{1}^{2} \!(\phi'(x+y) - \phi'(x-y))\smsp dx 
&= 
\phi(2+y)-\phi(2-y)-
\\
& -\phi(1+y)+\phi(1-y)
\leq 
24 y s_p^2.
\end{align*}
Therefore, integrating (\ref{ch31etheta3vt}) and recalling (\ref{ch31derivs}), we get
\begin{equation}
\int_{1}^{2}\!  \e \bigl| \bigl\la g_p(\sigma) \bigr\ra-\e \bigl\la g_p(\sigma) \bigr\ra \bigr| \smsp dx
\leq
24 \Bigl(y s_p  + \frac{v_N(s)}{ys_p}\Bigr).
\label{ch31step1}
\end{equation}
The minimum of the right hand side over $y>0$ is equal to $48 v_N(s)^{1/2}$ and is
achieved at $y=v_N(s)^{1/2}/s_p$. Throughout the argument we assumed that $0\leq y\leq 1$ and this is guaranteed by the condition $s^{-2}v_N(s) \leq 4^{-p}$. Combining (\ref{ch31step1}) with this optimal choice of $y$ and (\ref{ch31step2}) implies (\ref{ch31GGbasic}).
\qed

\medskip
\noindent
The proof of Theorem \ref{ch31ThGG} now follows, essentially, by Gaussian integration by parts.

\noindent
\textbf{Proof of Theorem \ref{ch31ThGG}.}
Let us fix $n\geq 2$ and consider a bounded function $f=f(R^n)$ of the overlaps of $n$ replicas. Without loss of generality, we can assume that $|f|\leq 1$. Then, 
\begin{equation}
\bigl|\e\bigl\la f g_p(\sigma^1) \bigr\ra 
- \e \bigl \la f \bigr\ra \e \bigl\la g_p(\sigma) \bigr\ra \bigr|
\leq
\e \bigl\la \bigl| g_p(\sigma) - \e \bigl\la g_p(\sigma) \bigr\ra\bigr|\bigr\ra. 
\label{ch31GGintegrate}
\end{equation}
We can think of the left hand side as a way to test the concentration of the process $(g_p(\sigma))$ on some function of the overlaps, $f(R^n)$. Using the Gaussian integration by parts formula in Lemma \ref{chAppLemma1} and recalling the factor $s 2^{-p} x_p$ in front of $g_p(\sigma)$ in (\ref{ch31Hpert}),
$$
\e \bigl\la g_p(\sigma) \bigr\ra
=
s 2^{-p} x_p \bigl(1 - \e \bigl\la R_{1,2}^p \bigr\ra \bigr).
$$
For the first term we use Lemma \ref{chAppLemma2}, 
$$
\e\bigl\la f g_p(\sigma^1) \bigr\ra 
=
s 2^{-p} x_p \e \Bigl\la f \Bigl( \sum_{\ell=1}^n R_{1,\ell}^p - nR_{1,n+1}^p \Bigr)\Bigr\ra. 
$$
Therefore, since the self-overlap $R_{1,1} = 1,$ the left hand side of (\ref{ch31GGintegrate}) equals $s 2^{-p} x_p n\Delta(f,n,p),$ where $\Delta(f,n,p)$ was defined in (\ref{ch31GG}). If we now integrate the inequality (\ref{ch31GGintegrate}) over $1\leq x_p\leq 2$ and use Theorem \ref{ch31L1}, we will get  
\begin{equation}
s 2^{-p} n
\int_1^2 \!
\Delta(f,n,p)
\smsp dx_p
\leq 2 + 48 \sqrt{v_N(s)},
\end{equation}
if $s^{-2}v_N(s) \leq 4^{-p}$. If we divide both sides by $s 2^{-p} n$ and then average over all $(x_p)$ on the interval $[1,2]$, by Fubini's theorem, we get 
\begin{equation}
\e_x \Delta(f,n,p)
\leq 
\frac{2^p}{n}
\Bigl(
\frac{2}{s} + 48 \frac{\sqrt{v_N(s)}}{s}
\Bigr),
\label{ch31GGbasicintegrated}
\end{equation}
if $s^{-2}v_N(s) \leq 4^{-p}$.  Finally, using this bound with $s=s_N$ that satisfies the condition (\ref{ch31GGassumption}) implies the Ghirlanda-Guerra identities in (\ref{ch31GGxlim}).
\qed

\section{Ultrametricity}

Let $G$ be any asymptotic Gibbs distribution (or, simply a random measure) on a Hilbert space $H$ and let us denote by $\la\,\cdot\,\ra$ the average with respect to $G^{\otimes \infty}$. As before, let $(\sigma^\ell)_{\ell\geq 1}$ be an i.i.d. sample from the measure $G$ and let 
\begin{equation}
R= \bigl(R_{\ell,\ell'}\bigr)_{\ell,\ell'\geq 1} =\bigl( \sigma^\ell \cdot \sigma^{\ell'}\bigr)_{\ell,\ell'\geq 1}
\label{ch44overlaps}
\end{equation} 
be the overlap array of the sequence $(\sigma^{\ell})$. Suppose that, for any $n\geq 1$ and any bounded measurable functions $f=f((R_{\ell,\ell'})_{\ell\not =\ell'\leq n})$ and $\psi: \Reals \to \Reals$, the Ghirlanda-Guerra identities hold,\vspace{-0.6mm}
\begin{equation}
\e\bigl\la  f \psi(R_{1,n+1})\bigr\ra
=
\frac{1}{n}\e\bigl\la  f\bigr\ra \e\bigl\la \psi(R_{1,2})\bigr\ra
+
\frac{1}{n} \sum_{\ell=2}^n \e\bigl\la  f \psi(R_{1,\ell})\bigr\ra.
\label{ch44GGI}
\end{equation}
It is important to point out that i.i.d. replicas $(\sigma^{\ell})$ play interchangeable roles and the index $1$ in $R_{1,n+1}$ can be replaced by any index $1\leq j\leq n$, in which case the sum on the right hand side will be over $1\leq \ell\leq n$ such that $\ell\not = j$.  Let us also emphasize that in (\ref{ch44GGI}) we consider functions $f$ that depend only on the off-diagonal elements $R_{\ell,\ell'}$ for $\ell\not = \ell'$ and do not depend on the self-overlaps $R_{\ell,\ell}$. In the SK model, the self-overlaps were constant by construction. An asymptotic Gibbs distribution $G$ does not automatically concentrate on a sphere in $H$, but we will prove in Theorem \ref{ch44ThSphere} as a consequence of the Ghirlanda-Guerra identities that, in fact, it does concentrate on a sphere.

Let us denote by $\zeta$ the distribution of the overlap $R_{1,2}$ under the measure $\e G^{\otimes 2}$, 
\begin{equation}
\zeta(A) =\e \bigl\la I\bigl(R_{1,2} \in A\bigr) \bigr\ra
\label{ch44fop}
\end{equation}
for any measurable set $A$ on $\Reals.$ In all applications, random measures $G$ will have bounded support in $H$ so, without loss of generality, we will assume that 
\begin{equation}
G(\sigma : \|\sigma\| \leq 1) = 1
\end{equation}
with probability one, in which case the overlaps $R_{\ell,\ell'}\in [-1,1]$. The main goal of this section is to prove the ultrametricity of the support of $G$.
\begin{theorem}\label{ch44ThUltra}
Suppose that the Ghirlanda-Guerra identities (\ref{ch44GGI}) hold. Then
\begin{equation}
\e \bigl\la I\bigl(R_{1,2}\geq \min(R_{1,3},R_{2,3})\bigr) \bigr\ra=1.
\label{ch44ultra}
\end{equation}
\end{theorem}
We will begin with a couple of basic observations. Our first observation shows that the Ghirlanda-Guerra identities determine the self-overlaps $R_{\ell,\ell}$ and they are, indeed, constant.
 \begin{theorem}\label{ch44ThSphere}
 Suppose that the Ghirlanda-Guerra identities (\ref{ch44GGI}) hold. If $q^*$ is the largest point in the support of $\zeta$ then, with probability one, $G(\|\sigma\|^2 = q^*)=1$. 
\end{theorem} \noindent
Because of this observation, the diagonal elements are non-random, $R_{\ell,\ell} = q^*$, and, if we wish, we can include them in the statement of the Ghirlanda-Guerra identities. The proof of Theorem \ref{ch44ThSphere} will be based on one elementary calculation.
\begin{lemma}\label{LemDS2}
Consider a measurable set $A\subseteq \Reals.$ With probability one over the choice of $G$:
\begin{enumerate}
\item[(a)] if $\zeta(A)>0$ then for $G$-almost all $\sigma^1,$ $G(\sigma^2: \sigma^1\cdot \sigma^2 \in A)>0$,
\item[(b)] if $\zeta(A) = 0$ then for $G$-almost all $\sigma^1,$ $G(\sigma^2: \sigma^1\cdot \sigma^2 \in A)=0$.
\end{enumerate}
\end{lemma}
\noindent\textbf{Proof.} (a) Suppose that $a=\zeta(A^c)<1.$ First of all, using the Ghirlanda-Guerra identities (\ref{ch44GGI}), \setstretch{1}
\begin{align*}
\e \bigl\la I\bigl(R_{1,\ell} \in A^c, 2\leq \ell\leq n+1\bigr) \bigr\ra
& =
\e \bigl\la I\bigl(R_{1,\ell} \in A^c, 2\leq \ell\leq n\bigr) I\bigl(R_{1,n+1} \in A^c\bigr) \bigr\ra
\\
&=
\frac{n-1+a}{n}\smsp \e \bigl \la I\bigl(R_{1,\ell} \in A^c, 2\leq \ell\leq n\bigr) \bigr\ra.
\end{align*}
\setstretch{1.05}
Repeating the same computation, one can show by induction on $n$ that this equals
$$
\frac{(n-1+a)\cdots(1+a) a}{n!}
=
\frac{a(1+a)}{n}\Bigl(1+\frac{a}{2}\Bigr)\cdots\Bigl(1+\frac{a}{n-1}\Bigr).
$$
Using the inequality $1+x\leq e^x$, it is now easy to see that
$$
\e \bigl\la I\bigl(R_{1,\ell} \in A^c, 2\leq \ell\leq n+1\bigr) \bigr\ra
\leq
\frac{a(1+a)}{n} e^{a\log n}
=
\frac{a(1+a)}{n^{1-a}}.
$$
If we rewrite the left hand side using Fubini's theorem then, since $a<1,$ letting $n\to \infty$ implies that 
$$
\lim_{n\to \infty} \e \int\! G( \sigma^2: \sigma^1\cdot \sigma^2 \in A^c)^n \smsp dG(\sigma^1) =0.
$$   
This leads to contradiction if we assume that $G(\sigma^2: \sigma^1\cdot \sigma^2 \in A^c)=1$ with positive probability over the choice of $G$ and the choice of $\sigma^1$, which proves part (a). Part (b) simply follows by Fubini's theorem.
\qed

\medskip
\noindent\textbf{Proof of Theorem \ref{ch44ThSphere}.}
Since $q^*$ is the largest point in the support of $\zeta$, 
$$
\zeta\bigl((q^*,\infty)\bigr) =0 \,\mbox{ and }\, \zeta\bigl([q^*- n^{-1},q^*] \bigr)>0
\,\mbox{ for all }\, 
n\geq 1.
$$ 
Using Lemma \ref{LemDS2}, we get that with probability one, for $G$-almost all $\sigma^1,$
\begin{equation}
G\bigl(\sigma^2 : \sigma^1\cdot \sigma^2 \leq q^* \bigr) = 1
\label{ch44muprop1}
\end{equation}
and, for all $n\geq 1$,
\begin{equation}
G\bigl(\sigma^2 : \sigma^1\cdot \sigma^2 \geq q^*- n^{-1} \bigr) >0.
\label{ch44muprop2}
\end{equation}
The equality in (\ref{ch44muprop1}) implies that $G(\|\sigma\|^2 \leq q^*)=1.$ Otherwise, there exists $\sigma\in H$ with $\|\sigma\|^2> q^*$ such that $G(B_\eps(\sigma))>0$ for any $\eps>0,$ where $B_\eps(\sigma)$ is the ball of radius $\eps$ centered at $\sigma.$ Taking $\eps>0$ small enough, so that $\sigma^1\cdot \sigma^2> q^*$ for all $\sigma^1,\sigma^2\in B_\eps(\sigma)$ contradicts (\ref{ch44muprop1}).  

Next, let us show that $G(\|\sigma\|^2<q^*)=0.$ Otherwise, $G(\|\sigma\|^2< q^* - \eps)>0$ for some small enough $\eps>0$,  while for all $\sigma^1\in \{\|\sigma\|^2< q^* - \eps\}$ and $\sigma^2\in \{\|\sigma\|^2\leq q^*\}$ we have 
$$
\sigma^1\cdot \sigma^2 < \sqrt{q^*(q^*-\eps)} < q^* - n^{-1}
$$
for some large enough $n\geq 1$. Since we already proved that $G(\|\sigma\|^2 \leq q^*)=1,$
this contradicts the fact that (\ref{ch44muprop2}) holds for all $n\geq 1$. 
\qed

\medskip
Our second observation shows that the Ghirlanda-Guerra identities imply that the overlap can take only non-negative values, which is known as the {\it Talagrand positivity principle}.
\begin{theorem}\label{ch44ThPosit} Suppose that the Ghirlanda-Guerra identities (\ref{ch44GGI}) hold. Then the overlap is non-negative, $\zeta([0,\infty))=1.$
 \end{theorem} 
\noindent\textbf{Proof.} Given a set $A,$ consider the event 
$A_n =\{R_{\ell,\ell'} \in A : \ell\not = \ell' \leq n\}$ and notice that 
\begin{equation}
I_{A_{n+1}} 
=  I_{A_n} \prod_{\ell\leq n} \bigl(1- I(R_{\ell,n+1}\not\in A) \bigr)
\geq I_{A_n} - \sum_{\ell\leq n} I_{A_n} I(R_{\ell,n+1}\not\in A).
\label{almostUltra}
\end{equation}
For each $\ell\leq n,$ the Ghirlanda-Guerra identities (\ref{ch44GGI}) imply that 
$$
\e \bigl\la I_{A_n} I(R_{\ell,n+1}\not\in A) \bigr\ra = \frac{\zeta(A^c)}{n} \e \bigl \la I_{A_n} \bigr\ra,
$$ 
and together with (\ref{almostUltra}) this gives 
$$
\e\bigl\la I_{A_{n+1}} \bigr\ra \geq \e\bigl\la I_{A_n} \bigr\ra -\zeta(A^c) \e\bigl\la I_{A_n} \bigr\ra 
= \zeta(A) \e\bigl\la I_{A_n} \bigr\ra  \geq \zeta(A)^{n},
$$ 
by induction on $n$. Therefore, if $\zeta(A)>0$, for any $n\geq 1$, with positive probability over the choice of $G$, one can find $n$ points $\sigma^1,\ldots, \sigma^n$ in the support of $G$ such that their overlaps $R_{\ell,\ell'} \in A.$ If $A=(-\infty,-\eps]$ for some $\eps>0$, this would imply that
$$
0\leq \bigl\|\sum_{\ell\leq n} \sigma^{\ell} \bigr\|^2 = \sum_{\ell,\ell'\leq n} R_{\ell,\ell'}\leq n q^*- n(n-1)\eps <0 
$$
for large $n$, and we can conclude that $\zeta((-\infty,-\eps])=0$ for all $\eps>0$. 
\qed

\medskip
The main idea of the proof of ultrametricity is that, due to the Ghirlanda-Guerra identities, the distribution of the overlaps $(R_{\ell,\ell'})$ is invariant under a large family of changes of density and, as we will show,
this invariance property contains a lot of information about the geometric structure of the measure $G$. The invariance property can be stated as follows. Given $n\geq 1$, we consider $n$ bounded measurable functions $f_1,\ldots, f_n:  \Reals\to\Reals$ and let
\begin{equation}
F(\sigma,\sigma^1,\ldots,\sigma^n) = f_1(\sigma\cdot\sigma^1)+\ldots+f_n(\sigma\cdot\sigma^n).
\label{ch45F1}
\end{equation}
For $1\leq \ell\leq n$, we define
\begin{equation}
F_{\ell}(\sigma,\sigma^1,\ldots,\sigma^n) = F(\sigma,\sigma^1,\ldots,\sigma^n)
 - f_{\ell}( \sigma\cdot\sigma^{\ell})+ \e \bigl \la f_{\ell}(R_{1,2}) \bigr\ra,
\label{ch45F2}
\end{equation}
and, for $\ell\geq n+1$, we define
\begin{equation}
F_{\ell}(\sigma,\sigma^1,\ldots,\sigma^n) =F(\sigma,\sigma^1,\ldots,\sigma^n).
\label{ch45F3}
\end{equation}
The definition (\ref{ch45F3}) for $\ell\geq n+1$ will not be used in the statement, but will appear in the proof of the following invariance property.\index{invariance}
\begin{theorem}\label{ch45Th1}
Suppose that the Ghirlanda-Guerra identities (\ref{ch44GGI}) hold and let $\Phi$ be a bounded measurable function of $R^n=(R_{\ell,\ell'})_{\ell,\ell'\leq n}.$ Then,
\begin{equation}
\e \bigl\la\Phi \bigr \ra =
\e\Bigl\la
\frac{ \Phi \exp \sum_{\ell=1}^{n} F_{\ell}(\sigma^{\ell},\sigma^1,\ldots,\sigma^n)}
{\la\exp F(\sigma,\sigma^1,\ldots,\sigma^n)\ra_{\hspace{-0.3mm}\mathunderscore}^n}
\Bigr\ra,
\label{ch45main}
\end{equation}
where the average $\la\,\cdot\,\ra_{\hspace{-0.3mm}\mathunderscore}$ with respect to $G$ in the denominator is in $\sigma$ only for fixed  $\sigma^1,\ldots, \sigma^n$, and the outside average $\la\,\cdot\, \ra$ of the ratio is in $\sigma^1,\ldots, \sigma^n$.
\end{theorem}
\noindent\textbf{Proof.}
Without loss of generality, let us assume that $|\Phi| \leq 1$ and suppose that $|f_{\ell}|\leq L$ for all $\ell\leq n$ for some large enough $L.$ For $t\geq 0$, let us define
\begin{equation}
\varphi(t) = 
\e\Bigl\la
\frac{\Phi \exp \sum_{\ell=1}^{n} t F_{\ell}(\sigma^{\ell},\sigma^1,\ldots,\sigma^n)}
{\la\exp t F(\sigma,\sigma^1,\ldots,\sigma^n)\ra_{\hspace{-0.3mm}\mathunderscore}^n}
\Bigr\ra.
\label{ch45varphitdefine}
\end{equation}
We will show that the Ghirlanda-Guerra identities (\ref{ch44GGI}) imply that the function $\varphi(t)$ is constant for all $t\geq 0$, proving the statement of the theorem, $\varphi(0)=\varphi(1).$ For $k\geq 1$, let us denote
$$
D_{n+k} = \sum_{\ell=1}^{n+k-1}F_{\ell}(\sigma^{\ell},\sigma^1,\ldots,\sigma^n)
-(n+k-1) F_{n+k}(\sigma^{n+k},\sigma^1,\ldots,\sigma^n).
$$
Recalling (\ref{ch45F3}) and using that the average $\la\,\cdot\,\ra_{\hspace{-0.3mm}\mathunderscore}$ is in $\sigma$ only, one can check by induction on $k$, that
$$
\varphi^{(k)}(t) = 
\e\Bigl\la
\frac{\Phi D_{n+1}\cdots D_{n+k} 
\exp \sum_{\ell=1}^{n+k} t F_{\ell}(\sigma^{\ell},\sigma^1,\ldots,\sigma^n)}
{\la\exp t F(\sigma,\sigma^1,\ldots,\sigma^n)\ra_{\hspace{-0.3mm}\mathunderscore}^{n+k}}
\Bigr\ra.
$$
Next, we will show that $\varphi^{(k)}(0)=0.$ If we introduce the notation 
$$
\Phi' = \Phi D_{n+1}\cdots D_{n+k-1},
$$ 
then $\Phi'$ is a function of the overlaps $(R_{\ell,\ell'})_{\ell,\ell'\leq n+k-1}$ and $\varphi^{(k)}(0)$ equals
\begin{align*}
&
\e\Bigl\la \Phi' \Bigl(
\sum_{\ell=1}^{n+k-1}F_{\ell}(\sigma^{\ell},\sigma^1,\ldots,\sigma^n)
-(n+k-1) F_{n+k}(\sigma^{n+k},\sigma^1,\ldots,\sigma^n)
\Bigr)
\Bigr\ra
\nonumber
\\
& = 
\sum_{j=1}^{n}
\e\Bigl\la \Phi' \Bigl(
\sum_{\ell\not = j,\ell=1}^{n+k-1}f_j(R_{j,l})
+ \e \bigl\la f_j(R_{1,2}) \bigr\ra
-(n+k-1) f_{j}(R_{j,n+k})
\Bigr)
\Bigr\ra
=0,
\end{align*}
by the Ghirlanda-Guerra identities (\ref{ch44GGI}) applied to each term $j$. Furthermore, since $|\Phi|\leq 1$, $|F_{\ell}| \leq L n$ and $|D_{n+k}| \leq 2L (n+k-1)n$, we can bound
\begin{align*}
\bigl|\varphi^{(k)}(t) \bigr|
&\,\leq \,
\Bigl(\prod_{\ell=1}^{k} 2L(n+\ell-1)n\Bigr) \,
\e\Bigl\la
\frac{
\exp \sum_{\ell=1}^{n+k} t F_{\ell}(\sigma^{\ell},\sigma^1,\ldots,\sigma^n)}
{\la\exp t F(\sigma,\sigma^1,\ldots,\sigma^n)\ra_{\hspace{-0.3mm}\mathunderscore}^{n+k}}
\Bigr\ra
\\
&\,=\,
\Bigl(\prod_{\ell=1}^{k} 2L(n+\ell-1)n\Bigr) \,
\e\Bigl\la
\frac{
\exp \sum_{\ell=1}^{n} t F_{\ell}(\sigma^{\ell},\sigma^1,\ldots,\sigma^n)}
{\la\exp t F(\sigma,\sigma^1,\ldots,\sigma^n)\ra_{\hspace{-0.3mm}\mathunderscore}^{n}}
\Bigr\ra,
\end{align*}
where the equality follows from the fact that the denominator depends only on the first $n$ coordinates and, recalling (\ref{ch45F3}), the average of the numerator in $\sigma^{\ell}$ for each $n<\ell\leq n+k$ will cancel exactly one factor in the denominator. Moreover, if we consider an arbitrary $T>0$, using that $|F_{\ell}| \leq L n$, the last ratio can be bounded by $\exp (2L T n^2)$ for $0\leq t\leq T$ and, therefore,
$$
\max_{0\leq t\leq T}\bigl|\varphi^{(k)}(t) \bigr| \leq \exp( 2L T n^2) \frac{(n+k-1)!}{(n-1)!}\,  (2Ln)^k.
$$
Since we proved above that $\varphi^{(k)}(0) = 0$ for all $k\geq 1$, using Taylor's expansion, we can write
$$
\bigl|\varphi(t)-\varphi(0) \bigr| 
\leq 
 \max_{0\leq s\leq t} \frac{|\varphi^{(k)}(s)|}{k!}t^k
\leq
\exp(2L T n^2)  \frac{(n+k-1)! }{k! \,(n-1)!} (2Ln t)^k.
$$
Letting $k\to \infty$ proves that $\varphi(t)=\varphi(0)$ for $0\leq t<(2Ln)^{-1}.$ This implies that for any $t_0<(2Ln)^{-1}$ we have $\varphi^{(k)}(t_0)=0$ for all $k\geq 1$ and, again, by Taylor's expansion for $t_0\leq t \leq T,$
\begin{align*}
\bigl|\varphi(t)-\varphi(t_0)\bigr| 
&\leq 
\max_{t_0\leq s\leq t} \frac{|\varphi^{(k)}(s)|}{k!}(t-t_0)^k
\\
&\leq 
\exp(2L T n^2)  \frac{(n+k-1)! }{k! \,(n-1)!} \bigl(2Ln (t-t_0)\bigr)^k.
\end{align*}
Letting $k\to\infty$ proves that $\varphi(t) = \varphi(0)$ for $0\leq t< 2(2Ln)^{-1}.$ We can proceed in the same fashion to prove this equality for all $t< T$ and note that $T$ was arbitrary.
\qed

\medskip
A special feature of the invariance property (\ref{ch45main}) is that it contains some very useful information not only about the overlaps but also about the Gibbs distribution of the neighborhoods of the replicas $\sigma^1,\ldots,\sigma^n.$ Let us give one simple example. 

\medskip
\noindent
\textbf{Example.}
Recall that the measure $G$ is concentrated on the sphere $\|h\|=\sqrt{q^*}$ and, for $q=q^*-\eps$, let $f_1(x) = t I(x \geq q)$ and $f_2 = \ldots = f_n =0.$ Then
$$
F(\sigma,\sigma^1,\ldots,\sigma^n) = t I(\sigma \cdot \sigma^1 \geq q)
$$
is a scaled indicator of a small neighborhood of $\sigma^1$ on the sphere $\|h\|=\sqrt{q^*}$. If we denote by $W_1 = G(\sigma : \sigma \cdot \sigma^1 \geq q)$ the Gibbs weight of this neighborhood then the average in the denominator in (\ref{ch45main}) is equal to
$$
\la\exp F(\sigma,\sigma^1,\ldots,\sigma^n)\ra_{\hspace{-0.3mm}\mathunderscore}
=
W_1 e^t + 1-W_1.
$$
Suppose now that the function $\Phi = I_A$ is an indicator of the event 
$$
A = \bigl\{(\sigma^1,\ldots, \sigma^n) : \sigma^1 \cdot \sigma^{\ell} <q \mbox{ for } 2\leq \ell\leq n \bigr\}
$$
that the replicas $\sigma^2,\ldots,\sigma^n$ are outside of this neighborhood of $\sigma^1$. Then, it is easy to see that
$$
\sum_{\ell=1}^{n} F_{\ell}(\sigma^{\ell},\sigma^1,\ldots,\sigma^n)
=
t \e \la I(R_{1,2} \geq q)\ra
= : t\gamma
$$
and the equation (\ref{ch45main}) becomes
\begin{equation}
\e \bigl\la I_A \bigr \ra =
\e\Bigl\la
I_A \smsp
\frac{e^{t\gamma}}
{(W_1 e^t + 1-W_1)^n}
\Bigr\ra,
\label{ch45mainEx}
\end{equation}
which may be viewed as a constraint on the weight $W_1$ and the event $A$, since this holds for all $t$. This is just one artificial example, but the idea can be pushed much further.
\qed

\medskip
Let us write down a formal generalization of Theorem \ref{ch45Th1} on which the proof of ultrametricity will be based. Consider a finite index set $\A.$ Given $n\geq 1$ and $\sigma^1,\ldots,\sigma^n\in H,$  let $(B_\alpha)_{\alpha\in\A}$ be some partition of the Hilbert space $H$ such that, for each $\alpha\in\A$, the indicator 
$
I_{B_\alpha}(\sigma) = I(\sigma\in B_\alpha)
$ 
is a measurable function of $R^n = (R_{\ell,\ell'})_{\ell,\ell'\leq n}$ and $(\sigma\cdot\sigma^{\ell})_{\ell\leq n}$. In other words, the sets in the partition are expressed in terms of some conditions on the scalar products between $\sigma,\sigma^1,\ldots,\sigma^n$. Let
\begin{equation}
W_\alpha=W_\alpha(\sigma^1,\ldots,\sigma^n)=G(B_\alpha)
\label{ch45WA}
\end{equation}
be the weights of the sets in this partition with respect to the measure $G$. Let us define a map $T$ by
\begin{equation}
W=(W_\alpha)_{\alpha\in\A}\to T(W) = 
\Bigl(\frac{\la I_{B_\alpha}(\sigma) \exp F(\sigma,\sigma^1,\ldots,\sigma^n )\ra_{\hspace{-0.3mm}\mathunderscore}}
{\la \exp F(\sigma,\sigma^1,\ldots,\sigma^n )\ra_{\hspace{-0.3mm}\mathunderscore}} \Bigr)_{\alpha\in\A}.
\label{ch45TA}
\end{equation}
Then the following holds.
\begin{theorem}\label{ch45Th2}
Suppose that the Ghirlanda-Guerra identities (\ref{ch44GGI}) hold. Then, for any bounded measurable function $\varphi:\Reals^{n\times n}\times\Reals^{|\A|}\to \Reals$,
\begin{equation}
\e \bigl\la  \varphi(R^n, W) \bigr\ra
=
\e\Bigl\la
\frac{ \varphi(R^n,T(W)) \exp \sum_{\ell=1}^{n} F_{\ell}(\sigma^{\ell},\sigma^1,\ldots,\sigma^n)}
{\la\exp F(\sigma,\sigma^1,\ldots,\sigma^n)\ra_{\hspace{-0.3mm}\mathunderscore}^n}
\Bigr\ra.
\label{ch45nA}
\end{equation}
\end{theorem}
\noindent\textbf{Proof.}
Let $n_\alpha\geq 0$ be some integers for $\alpha\in\A$ and let $m=n+\sum_{\alpha\in\A} n_\alpha.$ Let $(S_\alpha)_{\alpha\in\A}$ be any partition of $\{n+1,\ldots,m\}$ such that the cardinalities $|S_\alpha |=n_\alpha.$ Consider a continuous function $\Phi = \Phi(R^n)$ of the overlaps of $n$ replicas and let 
$$
\Phi' = \Phi(R^n) \prod_{\alpha\in\A}\varphi_\alpha,
\,\mbox{ where }\,
\varphi_\alpha = I \bigl(\sigma^{\ell} \in B_\alpha, \forall \ell\in S_\alpha \bigr).
$$
We will apply Theorem \ref{ch45Th1} to the function $\Phi'$, but since it now depends on $m$ coordinates,  we have to choose $m$ bounded measurable functions $f_1,\ldots, f_m$ in the definition (\ref{ch45F1}). We will choose the first $n$ functions to be arbitrary and we let $f_{n+1}=\ldots=f_m=0.$ First of all, integrating out the coordinates $(\sigma^{\ell})_{l>n}$, the left hand side of (\ref{ch45main}) can be written  as
\begin{equation}
\e \bigl\la \Phi'  \bigr\ra 
=
\e \Bigl\la \Phi(R^n) \prod_{\alpha\in\A} \varphi_\alpha \Bigr\ra
=
\e \Bigl\la\Phi(R^n) \prod_{\alpha\in\A} W_\alpha^{n_\alpha}(\sigma^1,\ldots,\sigma^n) \Bigr\ra,
\label{ch45lhscor}
\end{equation}
where $W_\alpha$'s were defined in (\ref{ch45WA}). Let us now compute the right hand side of (\ref{ch45main}). Since $f_{n+1}=\ldots=f_m=0,$ the coordinates $\sigma^{n+1},\ldots,\sigma^{m}$ are not present in all the functions defined in (\ref{ch45F1})--(\ref{ch45F3}) and we will continue to write them as functions of $\sigma,\sigma^1,\ldots,\sigma^n$ only. Then, it is easy to see that the denominator on the right hand side of (\ref{ch45main}) is equal to $\la\exp F(\sigma,\sigma^1,\ldots,\sigma^n)\ra_{\hspace{-0.3mm}\mathunderscore}^m$ and the sum in the numerator equals
$$
\sum_{\ell=1}^{n} F_{\ell}(\sigma^{\ell},\sigma^1,\ldots,\sigma^n)
+
\sum_{\ell=n+1}^{m} F(\sigma^{\ell},\sigma^1,\ldots,\sigma^n).  
$$
Since the function $\Phi$ and the denominator do not depend on  $(\sigma^{\ell})_{l>n}$, integrating the numerator in the coordinate $\sigma^{\ell}$ for $\ell\in S_\alpha$ produces a factor 
$$
\bigl\la I_{B_\alpha}(\sigma) \exp F(\sigma,\sigma^1,\ldots,\sigma^n )  \bigr\ra_{\hspace{-0.3mm}\mathunderscore}.
$$
For each $\alpha\in \A$, we have $|S_\alpha| = n_\alpha$ such coordinates and, therefore, the right hand side of (\ref{ch45main}) 
is equal to
$$
\e\Bigl\la
\frac{\Phi(R^n) \exp \sum_{\ell=1}^{n} F_{\ell}(\sigma^{\ell},\sigma^1,\ldots,\sigma^n)}
{\la\exp F(\sigma,\sigma^1,\ldots,\sigma^n)\ra_{\hspace{-0.3mm}\mathunderscore}^n}
\! \prod_{\alpha\in\A} \hspace{-0.3mm} \Bigl(\hspace{-0.1mm}
\frac{\la I_{B_\alpha} \exp F(\sigma,\sigma^1,\ldots,\sigma^n )\ra_{\hspace{-0.3mm}\mathunderscore}}
{\la \exp F(\sigma,\sigma^1,\ldots,\sigma^n )\ra_{\hspace{-0.3mm}\mathunderscore}}
\Bigr)^{n_\alpha}\Bigr\ra.
$$
Comparing this with (\ref{ch45lhscor}) and recalling the notation (\ref{ch45TA}) proves (\ref{ch45nA}) for
$$
\varphi(R^n, W) =\Phi(R^n) \prod_{\alpha\in \A} W_\alpha^{n_\alpha}.
$$
The general case then follows by approximation. First, we can approximate a continuous function $\phi$ on $[0,1]^{|\A|}$ by polynomials to obtain (\ref{ch45nA}) for products $\Phi(R^n) \phi(W)$. This, of course, implies the result for continuous functions $\varphi(R^n,W)$ and then for arbitrary bounded measurable functions.
\qed

\medskip
\textbf{Duplication property.}
To motivate the rest of the proof, let us recall that the physicist's picture predicts much more than ultrametricity. When we described the clustering process depicted in Figure \ref{Fig1}, we mentioned that the tree is infinitary, meaning that each $q_p$-cluster contains infinitely many $q_{p+1}$-subclusters. This implies the following \emph{duplication property}. Consider $n$ pure states indexed by the leaves $\alpha_1,\ldots, \alpha_n \in\Natural^r$ in Figure \ref{Fig1}, and pick one point $\sigma^\ell \in H_{\alpha_\ell}$ inside each state. Their overlaps can be written as
$
\sigma^{\ell} \cdot \sigma^{\ell'} = q_{\alpha_{\ell}\wedge \alpha_{\ell'}},
$
where $\alpha_{\ell} \wedge \alpha_{\ell'}$ was defined in (\ref{ch43wedge}). Suppose that $\alpha_1 \wedge \alpha_n$ takes the largest value among $\alpha_{\ell}\wedge \alpha_{\ell'}$. Since the tree is infinitary, we can always find another index $\alpha_{n+1}\in \Natural^r$ such that 
$$
\alpha_n \wedge \alpha_{n+1} = \alpha_1 \wedge \alpha_n
\,\mbox{ and }\,
\alpha_{\ell} \wedge \alpha_{n+1} = \alpha_{\ell} \wedge \alpha_n
\,\mbox{ for }\, 
l=1,\ldots, n-1. 
$$
This means that the leaves $\alpha_{n+1}$ and $\alpha_n$ are at exactly the same distance on the tree from all the other points and, moreover, they are not too close to each other, since $\alpha_{n+1}$ is at the same distance from $\alpha_n$ as the closest of the other leaves, in this case, $\alpha_1$. One can think of the pure state $H_{\alpha_{n+1}}$ as a duplicate of $H_{\alpha_{n}}$ and, in some sense, it is a non-trivial duplicate, since they are not too close to each other. Alternatively, if we pick a point $\sigma^{n+1}\in H_{\alpha_{n+1}}$, we can call it a duplicate of $\sigma^n$. It turns out that this possibility of always `duplicating a point' can hold only if the support of $G$ is ultrametric. Our strategy will be to prove this duplication property and show that it implies ultrametricity.

\begin{figure}[t]
\centering
\psfrag{Sigma1}{$\sigma^1$}\psfrag{Sigma2}{$\sigma^2$}\psfrag{Sigman2}{$\sigma^{n-2}$}\psfrag{Sigman1}{$\sigma^{n-1}$}\psfrag{Sigman}{$\sigma^n$}
\includegraphics[width=0.4\textwidth]{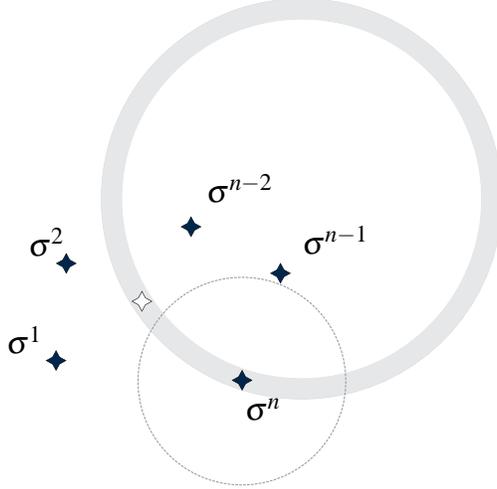}
\caption{\label{Fig2} {\it Duplication property.} The grey area corresponds to all the points on the sphere $\|h\|=c$ which are approximately at the same distance from the first $n-1$ replicas $\sigma^1,\ldots,\sigma^{n-1}$ as the last replica $\sigma^n$. Then the white point is a duplicate $\sigma^{n+1}$ of $\sigma^n$. It is in the grey area, so it is approximately at the same distances from the replicas $\sigma^1,\ldots, \sigma^{n-1}$ as $\sigma^n$, and it is at least as far from $\sigma^n$ as the closest of the first $n-1$ replicas, in this case $\sigma^{n-1}$.}
\end{figure}

Let us recall that, by Theorem \ref{ch44ThSphere}, the measure $G$ is concentrated on the sphere $\|\sigma\|^2 = q^*$ and from now on all $\sigma$'s will be on this sphere. Because of this, whenever we write that the scalar product $\sigma^1\cdot\sigma^2 \geq q$ is larger than some number, this means that the points $\sigma^1,\sigma^2$ are within some distance from each other and, vice versa, if we write that the scalar product $\sigma^1\cdot\sigma^2 \leq q$ is smaller than some number, this means that the points $\sigma^1,\sigma^2$ are separated by a certain distance (this relationship will be very helpful in visualizing the geometric picture since everything will be written in terms of scalar products). Consider a symmetric non-negative definite matrix 
\begin{equation}
A=\bigl(a_{\ell,\ell'}\bigr)_{\ell,\ell'\leq n}
\end{equation}
such that $a_{\ell,\ell} = q^*$ for $\ell\leq n$.
Given $\eps>0$, we will write $x\approx a$ to denote that $a-\eps<x<a+\eps$ and $R^n\approx A$ to denote that 
$R_{\ell,\ell'} \approx a_{\ell,\ell'}$ for all $\ell\not = \ell' \leq n$ and, for simplicity of notation, we will keep the dependence of $\approx$ on $\eps$ implicit.
Below, the matrix $A$ will be used to describe a set of constraints such that the overlaps in $R^n$ can take values close to $A$, 
\begin{equation}
\e\bigl\la I\bigl(R^n \approx A\bigr)\bigr\ra >0,
\label{ch45support}
\end{equation}
for a given $\eps>0$. Let us consider the quantity 
\begin{equation}
a_n^* = \max \bigl(a_{1,n},\ldots, a_{n-1,n}\bigr),
\end{equation}
which describes the constraint on the overlap corresponding to the closest replica among $\sigma^1,\ldots,$ $\sigma^{n-1}$ to the last replica $\sigma^n$. We will only consider the case when $a_n^*<q^*$, because, otherwise, the closest replica essentially coincides with $\sigma^n$. The following is the duplication property described above and also depicted in Figure \ref{Fig2}.
\begin{theorem} \label{ch45ThObs} Suppose that the Ghirlanda-Guerra identities (\ref{ch44GGI}) hold. Given $\eps>0$, if the matrix $A$ satisfies (\ref{ch45support}) and $a_n^* +\eps < q^*$ then
\begin{equation}
\e\bigl\la
I \bigl(
R^n\approx A, R_{\ell,n+1} \approx a_{\ell,n} \mbox{ for }\ell\leq n-1, R_{n,n+1} < a_n^* +\eps
\bigr)
\bigr\ra
>0.
\label{ch45extend}
\end{equation}
\end{theorem}\noindent
This result will be used in the following way. Suppose that $a_n^*<q^*$ and the matrix $A$ is in the support of the distribution of $R^n$ under $\e G^{\otimes \infty}$, which means that (\ref{ch45support}) holds for all $\eps>0$. Since $a_n^*+\eps<q^*$ for small $\eps>0$,  (\ref{ch45extend}) holds for all $\eps>0$. Therefore, the support of the distribution of $R^{n+1}$ under $\e G^{\otimes \infty}$ intersects the event in (\ref{ch45extend}) for every $\eps>0$ and, hence, it contains a point in the set
\begin{equation}
\bigl\{
R^{n+1} :
R^n = A, R_{\ell,n+1} = a_{\ell,n} \mbox{ for }\ell\leq n-1, R_{n,n+1} \leq a_n^*
\bigr\},
\label{ch45Aplus}
\end{equation}
since the support is compact.
\qed

\medskip
Before we prove the duplication property, let us show that it implies ultrametricity.

\medskip
\noindent\textbf{Proof of Theorem \ref{ch44ThUltra}.}
The proof is by contradiction. Suppose that (\ref{ch44ultra}) is violated, in which case there exist $a<b\leq c <q^*$ such that the matrix
\begin{equation}
\left(
\begin{array}{ l c r }
  q^* & a & b \\
  a & q^* & c \\
  b &  c & q^* \\
\end{array}
\right)
\label{ch45violated}
\end{equation}
is in the support of the distribution of $R^3$ under $\e G^{\otimes \infty}$, so it satisfies (\ref{ch45support}) for every $\eps>0$. 

In this case, Theorem  \ref{ch45ThObs} implies the following.  Given any $n_1, n_2, n_3\geq 1$ and $n=n_1+n_2+n_3$, we can find a matrix $A$ in the support of the distribution of  $R^n$ under $\e G^{\otimes \infty}$ such that  for some partition of indices $\{1,\ldots,n\} = I_1\cup I_2\cup I_3$ with $|I_j| = n_j$ we have $j\in I_j$ for $j\leq 3$ and 
\begin{description}
\item[(a)] 
$a_{\ell,\ell'}\leq c$ for all $\ell\not = \ell'\leq n$,

\item[(b)] 
$a_{\ell,\ell'} = a$ if $\ell\in I_1, \ell'\in I_2$, $a_{\ell,\ell'} = b$ if $\ell\in I_1, \ell'\in I_3$ and $a_{\ell,\ell'} = c$ if $\ell\in I_2, \ell'\in I_3$.
\end{description}
This can be proved by induction on $n_1,n_2,n_3.$ First of all, by the choice of the matrix (\ref{ch45violated}), this holds for $n_1=n_2=n_3=1.$ Assuming that the claim holds for some $n_1,n_2$ and $n_3$ with the matrix $A$, let us show how one can increase any of the $n_j$'s by one. For example, let us assume for simplicity of notation that $n\in I_3$ and show that the claim holds with $n_3+1.$ Since $a_n^*\leq c<q^*$, we can use the remark below Theorem \ref{ch45ThObs} to find a matrix $A'$ in the support of  the distribution of $R^{n+1}$ under $\e G^{\otimes \infty}$ that belongs to the set (\ref{ch45Aplus}). Hence, 
$$
\mbox{$a_{\ell,\ell'}' \leq c$ for all $\ell\not = \ell' \leq n+1$ and $a_{\ell,n+1}' = a_{\ell,n}$ for $\ell\leq n-1$,}
$$ 
so, in particular, 
$$
a_{\ell,n+1}' = b \,\mbox{ for }\, \ell\in I_1 \,\mbox{ and }\,  a_{\ell,n+1}' = c \,\mbox{ for }\, \ell\in I_2, 
$$
which means that $A'$ satisfies the conditions (a), (b)  with $I_3$ replaced by $I_3\cup \{n+1\}.$ In a similar fashion, one can increase the cardinality of $I_1$ and $I_2$ which completes the induction. 

Now, let $n_1=n_2=n_3=m$, find the matrix $A$ as above and find $\sigma^1,\ldots, \sigma^n$ on the sphere of radius $\sqrt{q^*}$ such that $R_{\ell,\ell'} = \sigma^{\ell}\cdot\sigma^{\ell'} =a_{\ell,\ell'}$ for all $\ell,\ell'\leq n.$ Let $\bar{\sigma}^j$ be the barycenter of  the set $\{\sigma^{\ell} : \ell\in I_j\}$. The condition (a) implies that  
$$
\|\bar{\sigma}^j\|^2 = \frac{1}{m^2}\sum_{\ell\in I_j} \|\sigma^{\ell}\|^2 + \frac{1}{m^2}\sum_{\ell\not = \ell'\in I_j} R_{\ell,\ell'}
\leq \frac{mq^* + m(m-1) c}{m^2},
$$
and the condition (b) implies that 
\begin{equation}
\bar{\sigma}^1\cdot \bar{\sigma}^2 = a, \bar{\sigma}^1\cdot \bar{\sigma}^3 = b \,\mbox{ and }\, \bar{\sigma}^2\cdot \bar{\sigma}^3 = c.
\end{equation}
Therefore, we can write
$$
\|\bar{\sigma}^2-\bar{\sigma}^3\|^2 = \|\bar{\sigma}^2\|^2 + \|\bar{\sigma}^3\|^2 
- 2\bar{\sigma}^2\cdot \bar{\sigma}^3 \leq \frac{2(q^* -c)}{m}
$$
and $0< b-a = \bar{\sigma}^1 \cdot \bar{\sigma}^3 -\bar{\sigma}^1\cdot \bar{\sigma}^2 \leq K m^{-1/2}.$ We arrive at contradiction by letting $m\to\infty$, which finishes the proof.
\qed

\medskip
It remains to prove the duplication property in Theorem \ref{ch45ThObs}. As we mentioned above, the proof will be based on the invariance property in the form of Theorem \ref{ch45Th2}. 

\begin{figure}[t]
\centering
\psfrag{Sigma1}{$\sigma^1$}\psfrag{Sigma2}{$\sigma^2$}\psfrag{Sigman2}{$\sigma^{n-2}$}\psfrag{Sigman1}{$\sigma^{n-1}$}\psfrag{Sigman}{$\sigma'$}
\includegraphics[width=0.4\textwidth]{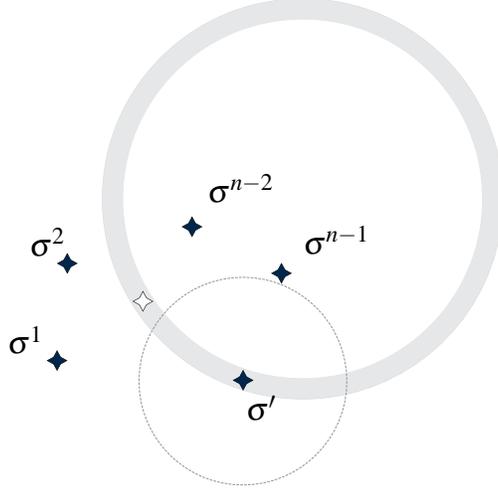}
\caption{\label{Fig3} {\it Proving the duplication property.} Grey circle is the set $\Omega(\sigma^1,\ldots,\sigma^{n-1}).$ The white point is a candidate for the duplicate of $\sigma'$.}
\end{figure}

\noindent\textbf{Proof of Theorem \ref{ch45ThObs}.}
We will prove (\ref{ch45extend}) by contradiction, so suppose that the left  hand side is equal to zero. We will apply Theorem \ref{ch45Th2} with ${\cal A}=\{1,2\}$ and the partition 
$$
B_1 = \bigl\{\sigma : \sigma\cdot \sigma^n \geq a_n^* +\eps \bigr\},\, 
B_2 = B_1^c.
$$ 
Since we assume that $a_n^* +\eps < q^*$, the set $B_1$ contains a small neighborhood of $\sigma^n$ on the sphere of radius $\sqrt{q^*}$ and, on the event $\{R^n\approx A\}$, its complement $B_2$ contains small neighborhoods of $\sigma^1,\ldots,\sigma^{n-1}$, since $R_{\ell,n}<a_{\ell,n}+\eps \leq a_n^* +\eps.$ Therefore, for $\sigma^1,\ldots,\sigma^n$ in the support of $G$, on the event $\{R^n\approx A\}$, the weights  
$$
W_1=G(B_1) \,\mbox{ and }\, W_2 = G(B_2)=1-W_1
$$ 
are strictly positive. Then, (\ref{ch45support}) implies that we can find $0<p<p'<1$ and small $\delta>0$ such that 
\begin{equation}
\e\Bigl\la
I\bigl(
R^n\approx A, W_1 \in (p,p')
\bigr)
\Bigr\ra \geq \delta.
\label{ch45littlec}
\end{equation}
Let us apply Theorem \ref{ch45Th2} with the above partition, the choice of 
\begin{equation}
\varphi(R^n,W) = I\bigl(R^n\approx A, W_1 \in (p,p')\bigr),
\end{equation}
and the choice of functions $f_1=\ldots=f_{n-1}=0$ and $f_{n}(x) = t I(x\geq a_n^*+\eps)$ 
for $t\in\Reals$. The sum in the numerator on the right hand side of (\ref{ch45nA}) will become 
\begin{align}
\sum_{\ell=1}^{n} F_{\ell}(\sigma^{\ell},\sigma^1,\ldots,\sigma^n) 
& = 
\sum_{\ell=1}^{n-1} t  I\bigl(R_{\ell,n} \geq a_n^*+\eps\bigr)
+
t\hspace{0.3mm}  \e \bigl\la I\bigl(R_{1,2}\geq a_n^* + \eps\bigr)\bigr \ra
\nonumber
\\
& = 
t\hspace{0.3mm} \e \bigl\la I\bigl(R_{1,2}\geq a_n^* + \eps\bigr)\bigr \ra =: t\gamma,
\nonumber
\end{align}
since, again, on the event $\{R^n\approx A\}$, the overlaps $R_{\ell,n}< a_{\ell,n}+\eps \leq a_n^*+\eps$ for all $\ell\leq n-1$, and the denominator will become
\begin{align}
\bigl \la \exp F(\sigma,\sigma^1,\ldots,\sigma^n) \bigr \ra_{\hspace{-0.3mm}\mathunderscore} 
&=
\bigl \la \exp t I\bigl(\sigma\cdot \sigma^n \geq a_n^*+\eps\bigr)  \bigr\ra_{\hspace{-0.3mm}\mathunderscore} 
\nonumber
\\
&=
G(B_1) e^t + G(B_2) =W_1 e^t + 1-W_1.
\label{ch45DeltatW}
\end{align}
If we denote $W=(W_1,W_2)$ and $\Delta_t(W) =W_1 e^t + 1-W_1$ then the map $T(W)$ in (\ref{ch45TA}) becomes
\begin{equation}
T_t(W) = \Bigl(\frac{W_1 e^t}{\Delta_t(W) }, \frac{1-W_1}{\Delta_t(W)} \Bigr).
\label{ch45TtW}
\end{equation}
Since $\Delta_t(W)\geq 1$ for $t\geq 0$, the equations (\ref{ch45nA}) and (\ref{ch45littlec}) imply
\begin{align}
\delta 
& \smsp\leq\smsp 
\e\Bigl\la
\frac{I(R^n\approx A, (T_t(W))_1\in (p,p')) \smsp  e^{ t\gamma} }
{ \Delta_t(W)^n}
\Bigr\ra
\nonumber
\\
& \smsp\leq\smsp 
\e\Bigl\la
I\bigl(R^n\approx A, (T_t(W))_1\in (p,p') \bigr)  \smsp e^{ t\gamma} 
\Bigr\ra.
\label{ch45littlec2}
\end{align}
In the average $\la\,\cdot\,\ra$ on the right hand side let us fix $\sigma^1,\ldots, \sigma^{n-1}$ and consider
the average with respect to $\sigma^n$ first. Clearly, on the event $\{R^n\approx A\}$ such average will be taken over the set
\begin{equation}
\Omega(\sigma^1,\ldots,\sigma^{n-1}) = \bigl\{\sigma: \sigma\cdot \sigma^{\ell} \approx a_{\ell,n} \mbox{ for }\ell\leq n-1\bigr\}.
\label{ch45Omega}
\end{equation}
Let us make the following crucial observation about the diameter of this set on the support of $G$. Suppose that with positive probability over the choice of the measure $G$ and replicas  $\sigma^1,\ldots, \sigma^{n-1}$ from $G$ satisfying the constraints in $A$, i.e. $R_{\ell,\ell'}\approx a_{\ell,\ell'}$ for $\ell,\ell'\leq n-1$,  we can find two points $\sigma'$ and $\sigma''$ in the support of $G$ that belong to the set $\Omega(\sigma^1,\ldots,\sigma^{n-1})$ and such that  $\sigma'\cdot \sigma'' < a_n^* + \eps.$ In other words, the matrix of overlaps of $\sigma^1,\ldots,\sigma^{n-1},\sigma'$ is approximately $A$, and $\sigma''$ is a candidate for a duplicate of $\sigma'$, as in Figure \ref{Fig3}. This would then imply
\begin{equation}
\e\bigl\la
I \bigl(
R^n\approx A, R_{\ell,n+1} \approx a_{\ell,n} \mbox{ for }\ell\leq n-1, R_{n,n+1} < a_n^* +\eps
\bigr)
\bigr\ra
>0,
\label{ch45extendAG}
\end{equation}
because for $(\sigma^n,\sigma^{n+1})$ in a small neighborhood of $(\sigma',\sigma'')$ the vector 
$(\sigma^1,\ldots,\sigma^n,\sigma^{n+1})$ would belong to the event on the left hand side,
$$
\bigl\{
R^n\approx A, R_{\ell,n+1} \approx a_{\ell,n} \mbox{ for }\ell\leq n-1, R_{n,n+1} < a_n^* +\eps
\bigr\}.
$$
Since we assumed that the left hand side of (\ref{ch45extendAG}) is equal to zero, we must have that, for almost all choices of the measure $G$ and replicas $\sigma^1,\ldots, \sigma^{n-1}$ satisfying the constraints in $A$, any two points $\sigma',\sigma''$ in the support of $G$ that belong to the set $\Omega(\sigma^1,\ldots,\sigma^{n-1})$ satisfy $\sigma'\cdot \sigma'' \geq a_n^* + \eps.$ In other words, given a point $\sigma'$ in Figure \ref{Fig3}, we can not find the white point $\sigma''$ in the support of $G$ such that  $\sigma'\cdot \sigma'' < a_n^* + \eps.$

Let us also recall that in (\ref{ch45littlec2}) we are averaging over $\sigma^n$ that satisfy the condition $(T_t(W))_1 \in (p,p').$ This means that if we fix any such $\sigma'$ in the support of $G$ that satisfies this condition and belongs to the set (\ref{ch45Omega}) then the Gibbs average in $\sigma^n$ will be taken over its neighborhood 
$$
B_1 = B_1(\sigma') = \bigl\{\sigma'': \sigma'\cdot \sigma'' \geq a_n^* +\eps \bigr\}
$$
of measure $W_1 = W_1(\sigma') = G(B_1(\sigma'))$ that satisfies $(T_t(W))_1\in (p,p').$ It is easy to check that the inverse of the map in (\ref{ch45TtW}) is $T_t^{-1} = T_{-t}$ and, using this for $(T_t(W))_1\in (p,p')$, implies that 
$$
W_1(\sigma') \in \Bigl\{\frac{q e^{-t}}{qe^{-t} + 1-q} : q\in (p, p') \Bigr\}
$$
and, therefore, $W_1(\sigma') \leq (1-p')^{-1}e^{-t}$. This means that the average on the right hand side of (\ref{ch45littlec2}) over $\sigma^n$ for fixed $\sigma^1,\ldots, \sigma^{n-1}$ is bounded by $(1-p')^{-1} e^{-t} e^{t\gamma }$ and, thus, for $t\geq 0$,
$$
0<\delta  
\leq 
\e\Bigl\la
I\bigl (R^n\approx A, (T_t(W))_1\in (p,p') \bigr)  e^{ t\gamma} 
\Bigr\ra 
\leq 
(1-p')^{-1} e^{-t(1-\gamma)}.
$$
Since the constraint matrix $A$ satisfies (\ref{ch45support}), $1-\gamma = \e\la I(R_{1,2}< a_n^* + \eps)\ra > 0$, and letting $t\to+\infty$ we arrive at contradiction.
\qed

\section{Reconstructing the limit}

We mentioned above that the Ghirlanda-Guerra identities and ultrametricity together determine the distribution of the overlap array $R=(R_{\ell,\ell'})_{\ell,\ell'\geq 1}$ uniquely in terms of the functional order parameter $\zeta$ and explained the case of three overlaps in Figure \ref{Fig10}. Now, we will prove the general case.

\begin{theorem}\label{ch44Th1}
Suppose that the Ghirlanda-Guerra identities (\ref{ch44GGI}) hold. Then the distribution of the entire overlap array $R=(R_{\ell,\ell'})_{\ell,\ell'\geq 1}$ under $\e G^{\otimes \infty}$ is uniquely determined by $\zeta$ in (\ref{ch44fop}). 
\end{theorem}
\noindent\textbf{Proof.}
By Theorem \ref{ch44ThSphere}, the diagonal elements $R_{\ell,\ell}$ are determined by $\zeta$ and we only need to consider the off-diagonal elements. By Talagrand's positivity principle in Theorem \ref{ch44ThPosit}, we can assume that all the overlaps are non-negative. Let us begin with the discrete case and suppose for a moment that the overlaps take only finitely many values 
\begin{equation}
0=q_0<q_1<\ldots<q_{r-1}<q_r,
\label{ch44qsagain}
\end{equation}
with probabilities 
\begin{equation}
\zeta(\{q_p\})
=
\e \bigl\la I(R_{1,2} = q_p) \bigr\ra
=
\delta_{p} 
\label{ch44overlapdist}
\end{equation}
for $p=0,\ldots, r$ and some $\delta_p\geq 0$ such that $\sum_{p=0}^r \delta_p= 1$. Some $\delta_p$ here can be equal to zero. In this case, we only need to show how to compute the probability of any particular configuration, 
\begin{equation}
\e \bigl\la I\bigl(R_{\ell,\ell'} = q_{\ell,\ell'} : \ell\not = \ell'\leq n+1\bigr) \bigr\ra,
\label{ch44ultraq1}
\end{equation}
for any $n\geq 1$ and any $q_{\ell,\ell'}\in \{q_0,\ldots, q_r\}$. Because of the ultrametricity property (\ref{ch44ultra}), we only have to consider $(q_{\ell,\ell'})$ that have this property, since, otherwise, (\ref{ch44ultraq1}) is equal to zero. Let us find the largest elements among $q_{\ell,\ell'}$ for $\ell\not = \ell'$ (say, equal to $q_p$) and suppose that $q_{1,n+1}$ is one of them. Let us consider the set of indices $2\leq \ell\leq n$ such that $q_{1,\ell} = q_{p}$ and, without loss of generality, suppose that this holds for $2\leq \ell\leq m$. Then the ultrametricity implies that (\ref{ch44ultraq1}) equals
\begin{align}
&
\e \bigl\la I\bigl(R_{\ell,\ell'} = q_{\ell,\ell'} : \ell,\ell'\leq n\bigr) I\bigl(R_{1,n+1} = q_{p}\bigr) \bigr\ra
\nonumber
\\
& - \sum_{\ell=2}^m \e \bigl\la I\bigl(R_{\ell,\ell'} = q_{\ell,\ell'} : \ell,\ell'\leq n\bigr) I\bigl(R_{\ell,n+1} >q_{p}\bigr) \bigr\ra,
\label{ch44ultraq20}
\end{align}
as follows. Let us consider the set $\{R_{\ell,\ell'} = q_{\ell,\ell'} : \ell,\ell'\leq n \}\cap \{R_{1,n+1} = q_{p}\}$ in the first term. Since we assumed that $q_{1,\ell}=q_{p}$ for $2\leq \ell\leq m,$ by ultrametricity, either one of the overlaps $R_{\ell,n+1}$ for such $\ell$ is strictly bigger than $q_{p}$ or all are equal to $q_{p}$, and these possibilities are disjoint. The second case is exactly the set in (\ref{ch44ultraq1}) and the first case can be written as $\{R_{\ell,\ell'} = q_{\ell,\ell'} : \ell,\ell'\leq n \}\cap \{R_{\ell,n+1} > q_{p}\}$, because the condition $R_{\ell,n+1} > q_{p}$ automatically implies that $R_{1,n+1} = R_{1,\ell}=q_p.$ This proves that (\ref{ch44ultraq1}) is equal to (\ref{ch44ultraq20}). By the Ghirlanda-Guerra identities (\ref{ch44GGI}) and (\ref{ch44overlapdist}), the first term in (\ref{ch44ultraq20}) equals 
$$
\frac{1}{n}\zeta\bigl(\{q_{p}\}\bigr)\e \bigl\la I\bigl(R_{\ell,\ell'} = q_{\ell,\ell'} : \ell,\ell'\leq n\bigr) \bigr\ra
+
\frac{1}{n}\sum_{\ell=2}^m \e \bigl\la I\bigl(R_{\ell,\ell'} = q_{\ell,\ell'} : \ell,\ell'\leq n \bigr) \bigr\ra
$$
and, similarly, each term in the sum over $2\leq \ell\leq m$ equals
$$
\frac{1}{n}\zeta\bigl((q_{p}, 1]\bigr)\e \bigl\la I\bigl(R_{\ell,\ell'} = q_{\ell,\ell'} : \ell,\ell'\leq n\bigr) \bigr\ra.
$$
Notice that all the terms now involve the set $\{R_{\ell,\ell'} = q_{\ell,\ell'} : \ell,\ell'\leq n\}$ that depends only on the indices $\ell\leq n$, so we can continue this computation recursively over $n$ and, in the end, (\ref{ch44ultraq1}) will be expressed completely in terms of the distribution of one overlap (\ref{ch44overlapdist}). 

It remains to show how the general case can be approximated by discrete cases as above. Given $r\geq 1$, let us consider a sequence of points as in (\ref{ch44qsagain}) and a function $\kappa(q)$ on $[0,1]$ such that 
\begin{equation}
\kappa(q) = q_p  \,\mbox{ if }\, q_p \leq q< q_{p+1}
\label{ch44finite}
\end{equation}
for $p=0,\ldots, r-1$ and $\kappa(q_r)=q_r=1$. Notice that, for $0\leq q\leq 1,$
\begin{equation}
| \kappa(q) - q | \leq \Delta_r := \max_{1\leq p\leq r} (q_p-q_{p-1}). 
\label{ch44Deltar}
\end{equation}
We will show that the discrete approximation\index{approximation} 
\begin{equation}
R^\kappa = \bigl(\kappa(R_{\ell,\ell'}) \bigr)_{\ell,\ell'\geq 1}
\label{ch44Rkappa}
\end{equation} 
of the matrix $R= (R_{\ell,\ell'})$ satisfies all the properties needed to apply the first part of the proof and, therefore, its distribution is uniquely determined by $\zeta.$ First of all, since the function $\kappa$ is non-decreasing, (\ref{ch44ultra}) obviously implies that $R^\kappa$ is also ultrametric. Moreover, ultrametricity implies that, for any $q$, the inequality $q\leq R_{\ell,\ell'}$ defines an equivalence relation $\ell\sim \ell'$, since $q\leq R_{\ell_1,\ell_2}$ and $q\leq R_{\ell_1,\ell_3}$ imply that $q\leq R_{\ell_2,\ell_3}$. Therefore, the array $(I(q\leq R_{\ell,\ell'}))_{\ell,\ell'\geq 1}$ is non-negative definite, since it is block-diagonal with blocks consisting of all elements equal to one. This implies that $R^\kappa$ is non-negative definite, since it can be written as a convex combination of such arrays,
\begin{equation}
\kappa\bigl(R_{\ell,\ell'}\bigr) = \sum_{p=1}^{r} (q_p-q_{p-1}) I\bigl(q_p \leq R_{\ell,\ell'}\bigr).
\label{ch44kappa1}
\end{equation}
It is obvious that $R^\kappa$ is replica symmetric in distribution (or exchangeable) under $\e G^{\otimes \infty}$ in the sense of the definition (\ref{ch12weakexch}), and the Dovbysh-Sudakov representation in Theorem \ref{DSintro} yields that (we use here that the diagonal elements of $R^\kappa$ are equal to $\kappa(q^*)$)
\begin{equation}
R^\kappa \stackrel{d}{=}
Q^\kappa =
\Bigl(
\rho^{\ell} \cdot \rho^{\ell'} + \delta_{\ell,\ell'} \bigl(\kappa(q^*) - \|\rho^{\ell}\|^2 \bigr)
\Bigr)_{\ell,\ell'\geq 1},
\label{ch44RkappaG}
\end{equation}
where $(\rho^{\ell})$ is an i.i.d. sample from some random measure $G'$ on a separable Hilbert space $H$. Let us denote by $\la\,\cdot\,\ra'$ the average with respect to ${G}^{\prime\smsp \otimes \infty}$ and, for any function $f$ of the overlaps $(R_{\ell,\ell'})_{\ell\not = \ell'}$, let us denote 
$
f_\kappa\bigl((R_{\ell,\ell'}) \bigr) = f\bigl((\kappa(R_{\ell,\ell'}))\bigr).
$ 
Since
$$
\e\bigl \la 
f\bigl((Q_{\ell,\ell'}^\kappa) \bigr)
\bigr\ra'
=
\e\bigl \la 
f\bigl((\kappa(R_{\ell,\ell'}))\bigr)
\bigr\ra
=
\e\bigl \la 
f_\kappa \bigl((R_{\ell,\ell'}) \bigr)
\bigr\ra,
$$
the Ghirlanda-Guerra identities (\ref{ch44GGI}) for the asymptotic Gibbs distribution $G$ imply that the measure $G'$ also satisfies these identities. The discrete case considered in the first part of the proof implies that the distribution of $Q^\kappa$ (and $R^\kappa$) is uniquely determined by the distribution of one element $\kappa(R_{1,2})$, which is given by the image measure $\zeta\circ \kappa^{-1}$. If we choose sequences (\ref{ch44qsagain}) in such a way that $\lim_{r\to\infty}\Delta_r = 0$ then, by (\ref{ch44Deltar}), $R^\kappa$ converges to $R$ almost surely and in distribution, which means that the distribution of the original overlap array is uniquely determined by $\zeta$. 
\qed

\medskip
It is not difficult to see that, given the overlap array $R=(R_{\ell,\ell'})_{\ell,\ell'\geq 1}$ of an infinite sample from $G$, we can reconstruct $G$ up to isometries (see Lemma 1.7 in \cite{SKmodel}).  This means that, in this sense, the functional order parameter $\zeta$ determines the randomness of $G$ and, as a result, it determines the randomness of Gibbs weights $G(H_\alpha)$ of all pure states described in Figure \ref{Fig1}. To complete the physicists' picture, one needs to show that these weights can be generated by the Ruelle Probability Cascades defined in (\ref{ch43vs}). This can be proved using the properties of Poisson processes in the construction of the Ruelle Probability Cascades, by showing that such Gibbs distribution will also satisfy the Ghirlanda-Guerra identities and will have the same functional order parameter, so it must produce the same overlap array $R=(R_{\ell,\ell'})_{\ell,\ell'\geq 1}$ in distribution. This is explained in detail in Chapter 2 in \cite{SKmodel}.

\section{Sketch of proof of the Parisi formula}\label{SecLabel-ParisiFormula}

\textbf{Guerra's replica symmetry breaking bound.} As we mentioned above, the fact that the Parisi formula in (\ref{ch30Parisi}) gives an upper bound on the free energy was proved in a breakthrough work of Guerra, \cite{Guerra}. This was the starting point of the proof of the Parisi formula by Talagrand \cite{TPF}.

The original argument of Guerra was simplified by Aizenman, Sims and Starr in \cite{AS2}, using some properties of the Ruelle Probability Cascades. The essence of Guerra's result is the interpolation between the SK model and the Ruelle Probability Cascades that we will sketch below.  Let us consider a sequence of parameters 
\begin{equation}
0=\zeta_{-1}< \zeta_0 <\ldots < \zeta_{r-1} < \zeta_r = 1
\end{equation}
as in (\ref{ch31zetas}), only now unrelated to any asymptotic Gibbs distribution. Let $(v_\alpha)_{\alpha\in\Natural^r}$ be the Ruelle Probability Cascades defined in (\ref{ch43vs}) in terms of these parameters. Consider a sequence
\begin{equation}
0=q_0< q_1 <\ldots <q_{r-1}< q_r =1,
\end{equation}
as in the definition of the clustering process in Figure \ref{Fig1} but, again, unrelated to any asymptotic Gibbs distribution. In particular, for some technical reason, we choose $q_r=1$. To these two sequences, we can associate the functional order parameter (a distribution function $\zeta$ on $[0,1]$) via
\begin{equation}
\zeta\bigl(\bigl\{q_p\bigr\}\bigr) = \zeta_{p} - \zeta_{p-1}
\,\,\mbox{ for }\,\,
p=0,\ldots, r.
\label{zetaGRSB}
\end{equation}
We will now sketch the proof of Guerra's upper bound.
\begin{theorem}
For any $\zeta$ defined as in (\ref{zetaGRSB}),
\begin{equation}
F_N(\beta) \leq \log 2 + \PP(\zeta) -{\beta^2}\int_{0}^{1}\! \zeta(t)t\, dt.
\label{FNPzeta2}
\end{equation}
\end{theorem}
\textbf{Sketch of proof.} Consider two Gaussian processes $(z(\alpha))$ and $(y(\alpha))$ indexed by the leaves ${\alpha\in \Natural^r}$ of the tree in Figure \ref{Fig1} with the covariance
\begin{equation}
\e z(\alpha) z(\beta) = 2 \smsp q_{\alpha\wedge \beta},\,\,
\e y(\alpha) y(\beta) = q_{\alpha\wedge\beta}^2.
\label{ch12CovzyG}
\end{equation}
Such processes are very easy to construct explicitly. Namely, let $(\eta_\alpha)_{\alpha \in \A}$ be a sequence of i.i.d. standard Gaussian random variables and, for each $p\geq 1$, let us define a family of Gaussian random variables indexed by $\alpha\in \Natural^r$,
\begin{equation}
g_p(\alpha) = \sum_{\beta\in p(\alpha)} \eta_\beta (q_{|\beta|}^p - q_{|\beta|-1}^p)^{1/2}.
\label{ch43etas}
\end{equation}
Recalling the notation (\ref{ch43wedge}), it is obvious that the covariance of this process is 
\begin{equation}
\e g_p(\alpha)g_p(\beta) = q_{\alpha\wedge \beta}^p,
\label{ch43Cov}
\end{equation}
so we can take $z(\alpha) = \sqrt{2} g_1(\alpha)$ and $y(\alpha) = g_2(\alpha)$. Where these processes are coming from will become clear when we discuss the lower bound.

Let  $(z_i(\alpha))$ and $(y_i(\alpha))$ for $i\geq 1$ be independent copies of the processes $(z(\alpha))$ and $(y(\alpha))$ in (\ref{ch12CovzyG}) and, for $0\leq t\leq 1$, let us consider the Hamiltonian
\begin{equation}
H_{N,t}(\sigma,\alpha) = 
\sqrt{t} H_N(\sigma) + \sqrt{1-t}\smsp \sum_{i=1}^N z_{i}(\alpha) \sigma_i 
+\sqrt{t}\smsp  \sum_{i=1}^N y_{i}(\alpha)
\label{ch33Hta}
\end{equation}
indexed by vectors $(\sigma,\alpha)$, where $H_N(\sigma)$ is the SK Hamiltonian. To the Hamiltonian $H_{N,t}(\sigma,\alpha) $ one can associate the free energy
\begin{equation}
\varphi(t)=\frac{1}{N}\smsp \e\log \sum_{\sigma,\alpha} v_{\alpha} \exp \beta H_{N, t}(\sigma,\alpha).
\label{ch33Gint}
\end{equation}
Let us denote by $\la\, \cdot\,\ra_t$ the average with respect to the Gibbs distribution $\Gamma_t(\sigma,\alpha)$ on $\Sigma_N\times \Natural^r$ defined by
$$
\Gamma_t(\sigma,\alpha) \sim v_\alpha \exp \beta H_{N, t}(\sigma,\alpha).
$$
Then, obviously, for $0<t<1$,
$$
\varphi'(t) = \frac{\beta}{N}\e \Bigl\la \frac{\partial H_{N, t}(\sigma,\alpha)}{\partial t} \Bigr\ra_t.
$$
It is easy to check from the above definitions that 
\begin{align*}
\frac{1}{N}\, \e \frac{\partial H_{N, t}(\sigma^1,\alpha^1)}{\partial t} H_{N, t}(\sigma^2,\alpha^2)
& = 
\frac{1}{2}
\bigl(
R_{1,2}^2 -2 R_{1,2} q_{\alpha^1\wedge\alpha^2} + q_{\alpha^1\wedge\alpha^2}^2
\bigr)
\\
& =
\frac{1}{2}
\bigl(
R_{1,2} - q_{\alpha^1\wedge\alpha^2} 
\bigr)^2.
\end{align*}
When $(\sigma^1,\alpha^1)=(\sigma^2,\alpha^2)$, this is equal to $0$, because $R_{1,1}=q_r = 1$. Using the 
Gaussian integration by parts in Lemma \ref{chAppLemma1} implies that
$$
\varphi'(t) = -\frac{\beta^2}{2}\e \bigl\la \bigl(
R_{1,2} - q_{\alpha^1\wedge\alpha^2} 
\bigr)^2 \bigr\ra_t \leq 0
$$
and, therefore, $\varphi(1)\leq \varphi(0).$  It is easy to see that
$$
\varphi(0)
=
\log 2 + 
\frac{1}{N}\smsp \e\log \sum_{\alpha\in\Natural^r} v_{\alpha} \prod_{i\leq N}
\ch\bigl(\beta z_{i}(\alpha) \bigr)
$$
and
$$
\varphi(1)
=
F_N(\beta)
+
\frac{1}{N}\smsp \e\log \sum_{\alpha\in\Natural^r} v_{\alpha} \prod_{i\leq N} 
\exp \bigl(\beta y_{i}(\alpha) \bigr),
$$
which implies that
$$
F_N(\beta)
\leq
\log 2 + 
\frac{1}{N}\smsp \e\log \sum_{\alpha\in\Natural^r} v_{\alpha} \prod_{i\leq N}
\ch\bigl(\beta z_{i}(\alpha) \bigr)
-
\frac{1}{N}\smsp \e\log \sum_{\alpha\in\Natural^r} v_{\alpha} \prod_{i\leq N} 
\exp \bigl(\beta y_{i}(\alpha) \bigr).
$$
The rest of the proof utilizes the properties of the Poisson processes involved in the construction of the Ruelle Probability Cascades. First, one can show that the independent copies for $i\leq N$ can be decoupled, so the right hand side above does not depend on $N$,
\begin{equation}
F_N(\beta)
\leq
\log 2 + \e\log \sum_{\alpha\in\Natural^r} v_{\alpha} \ch\bigl(\beta z(\alpha) \bigr)
-
\e\log \sum_{\alpha\in\Natural^r} v_{\alpha} \exp \bigl(\beta y(\alpha) \bigr).
\label{ch12AS2reprprelim}
\end{equation}
Then, one can show that the right hand side is exactly 
$$
\log 2 + \PP(\zeta) - {\beta^2}\int_{0}^{1}\! \zeta(t)t\, dt,
$$
and this finishes the proof.
\qed

\medskip
\noindent
\textbf{Lower bound via cavity method.} The first proof of the matching lower bound by Talagrand \cite{TPF} was quite technical, and here we sketch another proof from Panchenko \cite{PPF}, which is more conceptual. It is based on the properties of asymptotic Gibbs distributions that we proved above. In this approach, the lower bound can be proved by what is sometimes called the Aizenman-Sims-Starr scheme \cite{AS2}. In fact, the cavity computation itself was well known much earlier (see e.g. \cite{Guerra95}), and the main new ideas in \cite{AS2} were to turn it into a general variational principle, as well as using the Ruelle Probability Cascades to prove Guerra's upper bound as we did above.

Let us recall the definition of the partition function $Z_N=Z_N(\beta)$ in (\ref{ch11FE}) and, for $j\geq 0$, denote 
\begin{equation}
A_j = \e \log Z_{j+1} - \e \log Z_{j},
\end{equation}
with the convention that $Z_0 = 1$. Then we can rewrite the free energy as follows,
\begin{equation}
F_N = \frac{1}{N}\smsp \e \log Z_N
= \frac{1}{N} \sum_{j=0}^{N-1} A_j.
\label{ch12AS2liminf}
\end{equation}
Clearly, this representation implies that if the sequence $A_N$ converges then its limit is also the limit of the free energy $F_N$. Unfortunately, it is usually difficult to prove that the limit of $A_N$ exists and, therefore, this representation is used only to obtain a lower bound on the free energy, 
\begin{equation}
\liminf_{N\to\infty} F_N \geq \liminf_{N\to\infty} A_N.
\label{FNAN}
\end{equation}
Let us compare the partition functions $Z_N$ and $Z_{N+1}$ and see what they have in common and what makes them different. If we denote $\rho = (\sigma,\eps)\in \Sigma_{N+1}$ for $\sigma\in\Sigma_N$ and $\eps\in\{-1,+1\}$ then we can write
\begin{equation}
H_{N+1}(\rho) = H_N'(\sigma) + \eps z_N(\sigma),
\label{ch12decomp1}
\end{equation}
where
\begin{equation}
H_N'(\sigma) = \frac{1}{\sqrt{N+1}} \sum_{i,j =1}^N g_{ij}\sigma_i \sigma_j
\label{ch12commonH}
\end{equation}
and
\begin{equation}
z_N(\sigma) =  \frac{1}{\sqrt{N+1}} \sum_{i=1}^N \bigl(g_{i(N+1)} + g_{(N+1)i} \bigr)\sigma_i. 
\end{equation}
One the other hand, the part (\ref{ch12commonH}) of the Hamiltonian $H_{N+1}(\rho)$ is, in some sense, also a part of the Hamiltonian $H_N(\sigma)$ since, in distribution, the Gaussian process $H_N(\sigma)$ can be decomposed into a sum of two independent Gaussian processes
\begin{equation}
H_N(\sigma) \stackrel{d}{=}
H_N'(\sigma) + y_N(\sigma),
\label{ch12commonH2}
\end{equation}
where
\begin{equation}
y_N(\sigma) = 
 \frac{1}{\sqrt{N(N+1)}} \sum_{i,j =1}^N g_{ij}'\sigma_i \sigma_j
\end{equation}
for some independent array $(g_{ij}')$ of standard Gaussian random variables. Using the decompositions (\ref{ch12decomp1}) and (\ref{ch12commonH2}), we can write
\begin{equation}
\e\log Z_{N+1} 
= \log 2 + 
\e \log \sum_{\sigma\in\Sigma_N} \ch \bigl(\beta z_N(\sigma) \bigr) \exp \beta H_{N}'(\sigma)
\label{ch12ZN1}
\end{equation}
and
\begin{equation}
\e\log Z_{N} 
=
\e \log \sum_{\sigma\in\Sigma_N} \exp\bigl(\beta y_N(\sigma) \bigr) \exp \beta H_{N}'(\sigma).
\label{ch12ZN}
\end{equation}
Finally, if we consider the Gibbs distribution on $\Sigma_N$ corresponding to the Hamiltonian $H_N'(\sigma)$,
\begin{equation}
G_N'(\sigma) = \frac{\exp \beta H_N'(\sigma)}{Z_N'}
\,\mbox{ where }\,
Z_N' = \sum_{\sigma\in\Sigma_{N}} \exp \beta H_N'(\sigma),
\label{ch12MeasureGNprime}
\end{equation}
then (\ref{ch12ZN1}), (\ref{ch12ZN}) can be combined to give the representation
\begin{equation}
A_N= \log 2 +
\e \log \sum_{\sigma\in\Sigma_{N}} \ch\bigl(\beta z_N(\sigma)\bigr) G_N'(\sigma)
-
\e \log \sum_{\sigma\in\Sigma_{N}} \exp\bigl(\beta y_N(\sigma)\bigr) G_N'(\sigma).
\label{ch12AS2repr}
\end{equation}
Notice that the Gaussian processes $(z_N(\sigma))$ and $(y_N(\sigma))$ are independent of the randomness of the measure $G_N'$ and have the covariance
\begin{equation}
\e z_N(\sigma^1) z_N(\sigma^2) = 2 R_{1,2} + O(N^{-1}),\,\,
\e y_N(\sigma^1) y_N(\sigma^2) = R_{1,2}^2 + O(N^{-1}).
\label{ch12Covzy}
\end{equation}
Notice that (\ref{ch12AS2repr}) looks very similar to the right hand side of (\ref{ch12AS2reprprelim}), and the covariance of Gaussian processes (\ref{ch12Covzy}) is formally similar to the covariance of Gaussian processes (\ref{ch12CovzyG}). If we can show that the limit of $A_N$ over any subsequence can be approximated by expressions of the type 
\begin{equation}
\log 2 + \e\log \sum_{\alpha\in\Natural^r} v_{\alpha} \ch\bigl(\beta z(\alpha) \bigr)
-
\e\log \sum_{\alpha\in\Natural^r} v_{\alpha} \exp \bigl(\beta y(\alpha) \bigr)
\label{equpthere}
\end{equation}
as in (\ref{ch12AS2reprprelim}), we would prove the matching lower bound to Guerra's upper bound. This turns out to be rather straightforward using the theory we developed above (we refer to \cite{SKmodel} for details). Namely, when we define $A_N$, we can introduce the same perturbation term as in the proof of the Ghirlanda-Guerra identities. This would imply that the Gibbs distribution $G_N'$ in (\ref{ch12AS2repr}) satisfies the Ghirlanda-Guerra identities and any asymptotic Gibbs distribution, defined over any subsequence, can be characterized in terms of the functional order parameter -- the limiting distribution of one overlap. It is not difficult to show that the functional in (\ref{ch12AS2repr}) is continuous with respect to the distribution of the overlap array $(R_{\ell,\ell'})_{\ell,\ell'\geq 1}$ generated by $G_N'$ and, therefore, in the limit it can be approximated by (\ref{equpthere}), completing the proof of the lower bound and the Parisi formula.


\section{Phase transition}\label{SecLabel-PT}

We mentioned above that the physicists picture is non-trivial at low temperature. First of all, this means that the functional order parameter $\zeta$ that achieves the infimum in the Parisi formula (\ref{ch30Parisi}), which is called \emph{the Parisi measure}, is not concentrated on one point. The fact that the minimizer is unique follows from the strict convexity of the functional with respect to $\zeta$, which was conjectured in \cite{PConvex} (where a partial result was proved) and recently proved by Auffinger and Chen in \cite{AufChen2}. 

In the SK model without external field that we considered above, the phase transition occurs at $\beta = 1/\sqrt{2}$. This means that for $\beta\leq 1/\sqrt{2}$ the Parisi measure is concentrated on one point, and for $\beta>1/\sqrt{2}$ it is non-trivial. At high temperature, $\beta\leq 1/\sqrt{2}$, the Parisi formula is simplified to $\log 2 + \beta^2/2$, which corresponds to $\zeta$ concentrated at $0.$ This was first proved by Aizenman, Lebowitz and Ruelle in \cite{ALR}. This high temperature region is very well understood (see \cite{SG2, SG2-2}). The fact that the Parisi measure can not concentrate on one point for $\beta>1/\sqrt{2}$ was proved by Toninelli in \cite{TonAT} in a more general setting that includes the external field term $h\sum_{i=1}^N \sigma_i.$ In this general case, he showed that the Parisi measure can not concentrate on one point whenever
\begin{equation}
\beta^2 \e \ch^{-4}(\beta z \sqrt{2q} + h) > \frac{1}{2},
\label{ATline}
\end{equation}
where $z$ is a standard Gaussian random variable and $q$ is the largest solution of the equation
\begin{equation}
q = \e \mbox{th}^2 (\beta z \sqrt{2q} + h).
\label{qcritical}
\end{equation}
This last equation describes the critical point when we minimize the Parisi formula over all $\zeta$ concentrated on one point $q$. The Latala-Guerra lemma (see Appendix A.14 in \cite{SG2, SG2-2}) shows that if $h\not =0$ then (\ref{qcritical}) has a unique solution, and if $h=0$ then $q=0$ is a solution and for $\beta> 1/\sqrt{2}$ there is one more positive solution. The region (\ref{ATline}) was first described by de Almeida and Thouless in \cite{AT} and is known as (below) \emph{the dA-T line}. The high temperature region where the Parisi measure is concentrated on one point was predicted to be `above' the dA-T line,
\begin{equation}
\beta^2 \e \ch^{-4}(\beta z \sqrt{2q} + h) \leq \frac{1}{2}.
\label{ATline2}
\end{equation}
The first rigorous characterization of the high temperature region given by Talagrand in \cite{SG} looked formally different, but it was checked numerically that it coincides with (\ref{ATline2}) (it is still an open problem to show this analytically). Talagrand's description of the high temperature region can now also be obtained directly from the recent result of Auffinger and Chen in \cite{AufChen2}. Since they showed that the functional in the Parisi formula is strictly convex in $\zeta$, in order to decide whether the minimum over measures concentrated on one point (say, achieved on $\zeta_0$) is the global minimum over all $\zeta$, we only need to check that the derivative of the functional at $\zeta_0$ in the direction of any other measure concentrated on one point is non-negative. 

The functional order parameter also has the physical meaning of the distribution of one overlap, so it is interesting to see what can be said about this distribution at high and low temperature in the above sense. At high temperature, this is well understood and the overlap concentrates around the point $q$ defined in (\ref{qcritical}) (see Talagrand \cite{TalagrandSK} or \cite{SG2, SG2-2}). In particular, in the absence of external field, the overlap concentrates around $0$. At low temperature, the results are not as clear cut. If $\zeta$ is the Parisi measure, we saw in the discussion of the lower bound in the previous section that, at least, $\zeta$ will be the distribution of one overlap for some asymptotic Gibbs distribution corresponding to $G_N'$ in (\ref{ch12MeasureGNprime}). However, this statement does not say anything about the distribution of the overlap $R_{1,2}$ under $\e G_N^{\otimes 2}$ for all large enough  $N$. However, there is a class of models where the physical picture for the Gibbs distribution can be stated clearly in terms of the Parisi measure.  

\medskip
\noindent
\textbf{Generic mixed $p$-spin models.} Let us consider the {mixed $p$-spin model} with the Hamiltonian (\ref{mixedH}). First of all, the analogue of the Parisi formula can be proved for any such model (see \cite{PPF} or \cite{SKmodel}) and the limit of the free energy is given by (recall $\xi$ in (\ref{Covxi}))
\begin{equation}
F(\beta) = \lim_{N\to\infty} F_N(\beta) = \inf_{\zeta}\, \Bigl(\log 2 + \PP(\zeta) -\frac{1}{2}\int_{0}^{1}\! \xi''(t) \zeta(t)t\, dt\Bigr),
\label{ch30ParisiG}
\end{equation}
where the infimum is taken over all cumulative distribution functions $\zeta$ on $[0,1]$  and $\PP(\zeta) = f(0,0)$, where $f(t,x)$ is the solution of the parabolic differential equation
$$
\frac{\partial f}{\partial t} = -\frac{\xi''(t)}{2}
\Bigl(
\frac{\partial^2 f}{\partial x^2} + \zeta(t)
\Bigl(\frac{\partial f}{\partial x}\Bigr)^2
\Bigr)
$$
on $[0,1]\times \Reals$ with the boundary condition $f(1,x) = \log \ch(x)$. The Parisi measure that minimizes this functional is still unique, as shown in \cite{AufChen1}. (Various other interesting propeties of the Parisi measures in mixed $p$-spin models were obtained by Auffinger and Chen in \cite{AufChen1}.)

\medskip
\noindent
\textbf{Definition.}
We will say that a mixed $p$-spin model is \emph{generic} if the linear span of $1$ and power functions $x^p$ corresponding to $\beta_p \not = 0$ in (\ref{mixedH}) is dense in $(C[-1,1],\|\cdot\|_\infty)$. 

\medskip
\noindent
Each pure $p$-spin term contains some information about the $p$th moment of the overlap and, as a result, for generic models we can characterize the Gibbs distribution in the thermodynamic limit without the help of the perturbation term (\ref{ch31mixedHpert}) in the proof of the Ghirlanda-Guerra identities. More specifically, if $G_N$ is the Gibbs distribution corresponding to the Hamiltonian (\ref{mixedH}) of a generic mixed $p$-spin model then:
\begin{enumerate}
\item[(i)] the distribution of the overlap $R_{1,2}$ under $\e G_N^{\otimes 2}$ converges to some limit $\zeta^*$;

\item[(ii)] the distribution of the entire overlap array $(R_{\ell,\ell'})_{\ell,\ell'\geq 1}$ under $\e G_N^{\otimes \infty}$ converges;

\item[(iii)] the limit of the free energy is given by 
$$
\lim_{N\to\infty} F_N = 
\log 2 + \PP(\zeta^*) -\frac{1}{2}\int_{0}^{1}\! \xi''(t)\zeta^*(t)t\, dt.
$$
\end{enumerate}
In other words, the Parisi measure $\zeta^*$ that minimizes this functional is the limiting distribution of the overlap under $\e G_N^{\otimes 2}$. The proof of this statement can be found in Section 3.7 in \cite{SKmodel} and is based on two earlier results.

The first result is as follows. Let us denote by $\M$ the set of all limits over subsequences of the distribution of the overlap $R_{1,2}$ under $\e G_N^{\otimes 2}$. It was proved by Talagrand in \cite{PM} that, for each $p\geq 1$, the Parisi formula is differentiable with respect to  $\beta_p$ and 
\begin{equation}
\frac{\partial \PP}{\partial \beta_p} = 
\beta_p\Bigl(1-\int\! q^p \smsp d\zeta(q) \Bigr)
\label{ch36derP}
\end{equation}
for all $\zeta\in \M$ (a simplified proof was given in \cite{ECPs}). If $\beta_p \not = 0$, then this implies that all the limits $\zeta\in \M$ have the same $p$th moment and for generic models this implies that $\M = \{\zeta^*\}$ for some unique distribution $\zeta^*$ on $[-1,1]$, which is the statement (i) above. Here, we do not know yet that this distribution is on $[0,1]$. 

As a consequence of the differentiability of the Parisi formula, it was proved in Panchenko \cite{PGGmixed} that, whenever $\beta_p\not = 0$, the Ghirlanda-Guerra identities
\begin{equation}
\e  \bigl\la f R_{1,n+1}^p \bigr\ra 
=
 \frac{1}{n}\e \bigl\la f \bigr\ra \smsp \e\bigl\la R_{1,2}^p\bigr\ra 
 +
 \frac{1}{n}\sum_{\ell=2}^{n}\e \bigl\la f R_{1,\ell}^p\bigr\ra
\label{ch31GGlimitAgain}
\end{equation}
for the $p$th moment of the overlap hold in the thermodynamic limit in a strong sense -- for the Gibbs distribution $G_N$ corresponding to the original Hamiltonian (\ref{mixedH}) without the perturbation term (\ref{ch31mixedHpert}). For generic models, this implies the general Ghirlanda-Guerra identities in distribution. The Ghirlanda-Guerra uniquely determine the distribution of the entire overlap array in terms of  the distribution of one overlap. Since the limiting distribution of one overlap is unique, $\zeta^*,$ the distribution of the entire overlap array under $\e G_N^{\otimes\infty}$ also has a unique limit, which is the second statement (ii) above.  Notice also that, by Talagrand's positivity principle, the distribution $\zeta^*$ is, actually, supported on $[0,1]$.

Finally, one can show that, for generic models, the Gibbs distribution $G_N'$ in (\ref{ch12MeasureGNprime}) that appears in the cavity computation in the previous section satisfies all the same properties as $G_N$, so the cavity computation can be used to show the statement (iii) above.

\section{Some open problems}\label{SecLabel-Problems}

In conclusion, we will mention a few open problems, both in the setting of the SK model and beyond.

\medskip
\noindent
\textbf{Strong ultrametricity.} Suppose that the SK model has an external field term $h\sum_{i=1}^N\sigma_i$ with $h>0$. This term will remove the symmetry under flipping one replica $\sigma^1\to-\sigma^1$ and ensure that the overlap is non-negative in the limit (Section 14.12 of Talagrand \cite{SG2, SG2-2}). Otherwise, ultrametricity can not hold. Above, we mentioned an approach to ultrametricity suggested by Talagrand. In essence, it reduces to a technical problem of estimating a certain variational formula and relating it to the Parisi formula (see Section 15.7 in \cite{SG2, SG2-2} for details). Although this method is based on computations very specific to the SK model and it seems unlikely that it would work in other models, it is still of interest for a couple of reasons. Most importantly, it is designed to give a strong exponential control of the probability that ultrametricity is violated,
$$
\e G_N^{\otimes 3}\bigl(R_{2,3}\leq \min(R_{1,2}, R_{1,3}) - \eps\bigr) \leq Ke^{-N/K},
$$
for some constant $K$ independent of $N$. Moreover, this approach does not involve any perturbation of the Hamiltonian, which makes it appealing for aesthetic reasons. 
  
\medskip
\noindent
\textbf{Chaos problem.} Chaos problem is a general phenomenon suggested in the physics literature by Fisher, Huse \cite{FH86} and Bray, Moore \cite{BM87}. It states that any small change of parameters of the model results in complete reorganization of the Gibbs distribution. Chaos in temperature means that two Gibbs distributions $G_N$ and $G_N'$ corresponding to inverse temperature parameters $\beta$ and  $\beta'$, and defined with the same disorder $(g_{ij})$, will be concentrated on configurations almost orthogonal to each other,
$$
\e\, G_N\times G_N' \bigl((\sigma^1,\sigma^2) : |R_{1,2}| \geq \eps \bigr) \to 0,
$$ 
for any $\eps>0$, if $\beta \not = \beta'$. Chaos in external field is defined similarly for $h \not = h',$  except that two systems will have a preferred direction, so 
$$
\e\, G_N\times G_N' \bigl((\sigma^1,\sigma^2) : |R_{1,2} - q| \geq \eps \bigr) \to 0
$$ 
for any $\eps>0$ and some constant $q.$ Another type of chaos is chaos in disorder. One way to `slightly modify' the interaction parameters $g_{ij}$ is simply to assume that the interaction parameters $g_{ij}^1$ and $g_{ij}^2$ in the two models are strongly correlated, $\e g_{ij}^1 g_{ij}^2 = t,$ where $t<1$ is very close to $1$.  

Disorder chaos in the SK model was first proved by Chatterjee in Theorem 1.4 in \cite{Chatt09} (and Theorem 1.3 for a different modification of disorder), in the case when there is no external field. Disorder chaos in the presence of external field was proved by Chen in \cite{ChenChaos0}. The external field complicates the problem, because finding the constant $q$ around which the overlap concentrates is very non-trivial. 

First examples of chaos in external field and chaos in temperature for some classes of mixed $p$-spin models were given by Chen in \cite{ChenChaos}, building upon the ideas in Chen, Panchenko in \cite{ChenChaos2}. However, the most basic cases, such as the SK model, are still open.

\medskip
\noindent
\textbf{Concentration of the ground state and free energy.} Let us considers the free energy $F_N(\beta)$ in the SK model at positive temperature, or the ground state $\max_\sigma H_N(\sigma)/N$. Classical Gaussian concentration inequalities guarantee that the fluctuations of these random quantities around their expected values are of order at most $N^{-1/2}$. At high temperature, $\beta < 1/\sqrt{2}$, the fluctuations were proved to be of order $1/N$ by Aizenman, Lebowitz, Ruelle in \cite{ALR}. The physicists  predicted that fluctuations of the free energy at low temperature, $\beta>1/\sqrt{2},$ and fluctuations of the ground state should be of order $N^{-5/6}$ (see Crisanti, Paladin, Sommers, Vulpiani \cite{CPSV},  and Parisi, Rizzo \cite{ParRiz}), which is an open problem.

The best known results are due to Chatterjee. In Theorem 1.5 in \cite{Chatt09}, he showed that the fluctuations of the free energy $F_{N}(\beta)$ in the SK model are of order at most $1/\sqrt{N\log N}$, at any positive temperature $\beta>0.$ In Theorem 9.2 in \cite{Chatt08}, he showed that for a certain class of mixed mixed $p$-spin models, the fluctuations are at most of order $N^{-5/8}$. See Chatterjee \cite{Chatt14} for a detailed exposition.

\medskip
\noindent
\textbf{Properties of the Parisi measure.} We already mentioned the papers of Auffinger, Chen \cite{AufChen1, AufChen2}, where a number of interesting properties of the Parisi measure were obtained. However, most of the conjectures of the physicists (discovered numerically) are still open, and we refer to \cite{AufChen1} for an overview. One related problem that was mentioned in the previous section is to show that Talagrand's characterization of the high temperature region coincides with the dA-T line.

\medskip
\noindent
\textbf{Multi-species and bipartite SK model.} Let us consider the following modification of the SK model. Consider a finite set $\spe$. The elements of $\spe$ will be called species. Let us divide all spin indices into disjoint groups indexed by the species,
\begin{equation}
I = \bigl\{1\, ,\ldots,\, N\bigr\} = \bigcup_{s\in\spe} I_s.
\label{species}
\end{equation}
These sets will vary with $N$ and we will assume that their cardinalities $N_s =|I_s|$ satisfy 
\begin{equation}
\lim_{N\to\infty} \frac{N_s}{N} = \lambda_s \in (0,1)
\,\mbox{ for all }\, s\in\spe.
\end{equation}
For simplicity of notation, we will omit the dependence of $\lambda_s^N := N_s/N$ on $N$ and will simply write $\lambda_s$. The Hamiltonian of the multi-species SK model resembles the usual SK Hamiltonian,
\begin{equation}
H_N(\sigma) = \frac{1}{\sqrt{N}} \sum_{i,j =1}^N g_{ij}\sigma_i \sigma_j,
\label{SKHms}
\end{equation}
where the interaction parameters $(g_{ij})$ are independent Gaussian random variables, only now they are not necessarily identically distributed but, instead, satisfy
\begin{equation}
\e g_{ij}^2 = \Delta_{st}^2 \,\mbox{ if }\, i\in I_s, j\in I_t \,\mbox{ for }\, s,t\in\spe.
\label{Delta}
\end{equation}
In other words, the variance of the interaction between $i$ and $j$ depends only on the species they belong to. In the case when the matrix $\Delta^2 = (\Delta_{st}^2)_{s,t\in\spe}$ is symmetric and nonnegative definite, the analogue of the Parisi formula for the free energy was proposed by Barra, Contucci, Mingione, Tantari in \cite{MS}. Using a modification of Guerra's interpolation, they showed that their formula gives an upper bound on the free energy. The matching lower bound was proved in Panchenko \cite{multi-species} utilizing a new multi-species version of the Ghirlanda-Guerra identities to show that the overlaps within species are deterministic functions of the global overlaps. 

The proof of the lower bound does not use the condition $\Delta^2\geq 0,$ but it is used in the proof of the upper bound via the Guerra interpolation. An important open problem is to understand what happens when $\Delta^2$ is not positive definite. One example would be a bipartite SK model, where $\spe = \{1,2\}$ and $\Delta_{11}^2 = \Delta_{22}^2 = 0$, i.e. only interactions between species are present.

\medskip
\noindent
\textbf{The Edwards-Anderson model.} If one considers the Edwards-Anderson model (\ref{EA}) on the lattice, one can use the same perturbation as in the SK model to obtain the Ghirlanda-Guerra identities and all their consequences in the thermodynamic limit (see e.g. Contucci, Mingione, Starr \cite{CEA}). However, it is a very difficult basic problem to understand whether the distribution of the overlaps can be non-trivial, which means that there is replica symmetry breaking. In the SK model, this phase transition information  was extracted from the Parisi formula, and there is no analogue of the Parisi formula that is expected in the EA model.

There is another, more classical, point of view in the EA model taking into account that the model is defined on a finite block of the infinite lattice $\mathbb{Z}^d$, and that the boundary spins can be connected to their neighbours outside of the finite block. The interactions with these outside nearest neighbours can be incorporated into the model, and it becomes a very important feature of this point of view to understand the effect of the boundary conditions on the behaviour of the Gibbs distribution inside the block. The questions here are notoriously difficult and we refer to Newman, Stein \cite{NS} for up-to-date results and open problems.

\medskip
\noindent
\textbf{Diluted models.} Diluted models, such as the diluted SK model or random $K$-sat model mentioned in Section \ref{SecLabel-Beyond}, are some of the most natural spin glass models to study, now that the SK model has been understood reasonably well. 

There is the analogue of the Parisi formula originating from the work of M\'ezard and Parisi in \cite{Mezard}, which comes with detailed predictions about the structure of the Gibbs distribution that build upon the picture in the SK model that we described above. There is the analogue of the Guerra upper bound proved by Franz, Leone in \cite{FL} (see Panchenko, Talagrand \cite{PT} for another exposition). There is a natural analogue of the cavity method (see Panchenko \cite{Pspins}), so proving the matching lower bound reduces to demonstrating the structure of Gibbs distribution predicted by M\'ezard and Parisi. The most recent results in this direction in Panchenko \cite{finiteRSB} are at the same stage as the SK model was after the work of Arguin, Aizenman \cite{AA} and Panchenko \cite{PGG}, namely, the M\'ezard-Parisi picture is proved in some generic sense (sufficient for proving the formula for the free energy) under the technical condition that the overlap takes finitely many values. 

Even if one could finish the proof of the M\'ezard-Parisi formula for the free energy, it would only be a starting point in the study of most interesting questions in these models related to the phase transition. The variational problem that defines this formula is not very well understood, since it involves a new complicated order parameter. For example, even numerical study of this variational problem relies on heuristic algorithms that are not completely justified (see Chapter 22 in M\'ezard, Montanari \cite{MezMont}). 

\medskip
\noindent
\textbf{Other models.} Above, we gave only a small sample of problems and references, and one can find much more information, for example, in \cite{SG2, SG2-2, NS, MezMont}.

\end{document}